\documentclass[a4,12pt,twoside]{article}
\usepackage{amsfonts}
\usepackage{amssymb}
\usepackage{amsthm}
\usepackage{amsmath}
\usepackage[dvips]{graphicx,color}
\usepackage{latexsym,bm}

\numberwithin{equation}{section}

\def\R2{\mathbb{R}^2}

\begin{document}           

\begin{flushleft}
\renewcommand{\thefootnote}{\fnsymbol{footnote}}
\center{\bf\LARGE  Exact spectrum of the Laplacian on a domain in
the Sierpinski gasket
}  

\footnotetext{This research  was supported by the National Science
Foundation of China, Grant 10901081. }
 \vskip0.5cm \center{\large
 HUA QIU}


\end{flushleft}

\noindent{\bf Abstract.} For a certain domain $\Omega$ in the
Sierpinski gasket $\mathcal{SG}$ whose boundary is a line segment, a
complete description of the eigenvalues of the Laplacian, with an
exact count of dimensions of eigenspaces, under the Dirichlet and
Neumann boundary conditions is presented. The method developed in
this paper is a weak version of the
 spectral decimation method due to Fukushima and Shima,
since for a lot of ``bad" eigenvalues the spectral decimation method
can not be used directly. Let $\rho^0(x)$, $\rho^\Omega(x)$ be the
eigenvalue counting functions of the Laplacian associated to
$\mathcal{SG}$ and $\Omega$ respectively. We prove a comparison
between $\rho^0(x)$ and $\rho^\Omega(x)$ says that $ 0\leq
\rho^{0}(x)-\rho^\Omega(x)\leq C x^{\log2/\log 5}\log x$ for
sufficiently large $x$ for some positive constant $C$. As a
consequence, $\rho^\Omega(x)=g(\log x)x^{\log 3/\log
5}+O(x^{\log2/\log5}\log x)$ as $x\rightarrow\infty$,
 for some (right-continuous
discontinuous) $\log 5$-periodic function $g:\mathbb{R}\rightarrow
\mathbb{R}$ with $0<\inf_{\mathbb{R}}g<\sup_\mathbb{R}g<\infty$.
Moreover, we explain that the asymptotic expansion of
$\rho^\Omega(x)$ should admit a second term of the order $\log2/\log
5$, that becomes apparent from the experimental data. This is very
analogous to the conjectures of Weyl and Berry.



\hspace{0.2cm}
\\

\section{Introduction}
The study of the Laplacian on fractals was originated by S. Kusuoka
$\cite{Ku}$ and S. Goldstein $\cite{G}$. They independently
constructed the Laplacian as the generator of a diffusion process on
the Sierpinski gasket $\mathcal{SG}$. Later an analytic approach was
developed by J. Kigami $\cite{Ki1}$, who constructed the Laplacian
both as a renormalized limit of difference operators and a weak
formulation using the theory of Dirichlet forms.

Let $V_0$ be the boundary of $\mathcal{SG}$, which consists of the
three vertices of the equilateral triangle containing
$\mathcal{SG}$. Consider the following Dirichlet eigenvalue problem:
\begin{eqnarray*}\left\{
\begin{array}{l}
  -\Delta u=\lambda u \mbox{ in } \mathcal{SG}\setminus{V_0},\\
  u|_{V_0}=0,
\end{array}\right.
\end{eqnarray*}
where $\Delta$ is the standard Laplacian (with respect to the
standard self-similar measure $\mu$) on $\mathcal{SG}$.
 Physicists R. Rammal and G. Toulouse $\cite{Ra}$
found that an appropriate choice of a series of eigenvalues of
successive difference operators produces an orbit of a dynamical
system related to a quadratic polynomial, and all the eigenvalues of
$-\Delta$ on $\mathcal{SG}\setminus{V_0}$ should be obtained by
tracking back the orbits. This is the phenomenon which M. Fukushima
and T. Shima $\cite{Fu,Shi1}$ described from the mathematical point
of view, by saying that $\mathcal{SG}$ admits \emph{spectral
decimation} with respect to a quadratic polynomial. Using the
spectral decimation, all the eigenvalues and eigenfunctions of
$-\Delta$ on $\mathcal{SG}\setminus{V_0}$ have been determined
exactly. This method also works for the eigenvalue problem of
$-\Delta$ with Neumann boundary condition.

Later the theory of the Laplacian was developed for nested fractals
and p.c.f. self-similar sets by T. Lindstr{\o}m $\cite{Lin}$ and
Kigami $\cite{Ki2}$ by introducing the notion of \emph{harmonic
structure}. Every p.c.f. self-similar set is approximated by an
increasing sequence of finite graphs and the harmonic structure
determines a sequence of difference operators on the successive
graphs, which converges to the Laplacian. Then some generalizations
of the spectral decimation to a class of p.c.f. self-similar sets
were developed by Shima $\cite{Shi2}$, L. Malozemov and A. Teplyaev
$\cite{Ma}$, in which some strong symmetry conditions are supposed
to be satisfied to ensure the spectral decimation applies to the
corresponding graph sequences. Under such strong symmetry
conditions, the spectrum of the Laplacian can also be determined in
terms of the iteration of a rational function. Recently, the
spectrum of the Laplacian  on some other fractals has been analyzed
either numerically $\cite{Ada}$ or using the spectral decimation
method $\cite{Con,Dren,Te,Zhou}$ by R. S. Strichartz (with
co-authors), D. Zhou and  Teplyaev. In all the references mentioned
above, spectral decimation plays a key role in the theoretical study
of the spectrum of the Laplacian.

The Weyl asymptotic behavior of the eigenvalue counting function of
$-\Delta$ on $\mathcal{SG}\setminus{V_0}$  has also been studied by
Fukushima and Shima $\cite{Fu}$. Afterwards,  a general spectral
distribution theory on p.c.f. self-similar sets was obtained by
Kigami and M. L. Lapidus $\cite{KiL,Ki5}$. Denote by $\rho^{0}(x)$
the number of eigenvalues of $-\Delta$(taking the multiplicities
into account) on $\mathcal{SG}\setminus{V_0}$ not exceeding $x$,
with Dirichlet boundary condition at the three vertices. As proved
in $\cite{Fu,KiL}$, there exist positive constant $c,C$ such that
\begin{equation}\label{abc}
cx^{d_S/2}\leq\rho^{0}(x)\leq C x^{d_S/2}
\end{equation}
 for
sufficiently large $x$, where $d_S={\log9}/{\log5}$ is the
\emph{spectral dimension} of $\mathcal{SG}$. In particular,
$\rho^0(x)$ varies highly irregularly at $\infty$ due to the high
multiplicities of localized eigenfunctions,
\begin{equation}\label{sss}
0<\liminf_{x\rightarrow\infty}\rho^{0}(x)x^{-{d_S}/{2}}<\limsup_{x\rightarrow\infty}\rho^{0}(x)x^{-{d_S}/{2}}<\infty.
\end{equation}
Furthermore, using a refinement of the Renewal Theorem, Kigami
$\cite{Ki5}$ showed that the remainder of $\rho^0(x)$ is bounded,
\begin{equation}\label{ssssss}
\rho^{0}(x)=g(\log x)x^{d_S/2}+O(1)\quad{ as } \quad
x\rightarrow\infty, \end{equation} for some (right-continuous
discontinuous) $\log 5$-periodic function $g:\mathbb{R}\rightarrow
\mathbb{R}$ with $0<\inf_{\mathbb{R}}g<\sup_\mathbb{R}g<\infty$.
Exactly the same results hold for the eigenvalue counting function
for the Neumann Laplacian.

In this paper, we are mainly concerned with eigenvalue problems for
a domain in $\mathcal{SG}$. Although analysis on fractals has been
made possible by the definition of the Laplacian, there has been
little research into differential equations on bounded subsets of
fractals.  Recall that $\mathcal{SG}$ is the attractor of the
\emph{iterated function system} $\{F_0,F_1,F_2\}$ with
$F_ix=\frac{1}{2}(x+q_i)$ where $q_0,q_1,q_2$ are the vertices of an
equilateral triangle in the plane,
$$\mathcal{SG}=\bigcup_{i=0}^2F_i(\mathcal{SG}).$$
In Kigami's theory the boundary of $\mathcal{SG}$ consists of the
three points $q_0,q_1,q_2$, and the space of \emph{harmonic
functions} (solutions of $\Delta u=0$) is three dimensional, with
$u$ determined explicitly by its boundary values $u(q_i)$. (Note
that this boundary is not a topological boundary.) Thus the harmonic
function theory on $\mathcal{SG}$ is more closely related to the
theory of linear functions on the unit interval than to harmonic
functions on the disk. To get a richer theory we should take an open
set $\Omega$ in $\mathcal{SG}$ and restrict the Laplacian on
$\mathcal{SG}\setminus{V_0}$ to functions defined on $\Omega$. Hence
we believe it is appropriate to begin the study of differential
equations related to a bounded domain $\Omega$ in $\mathcal{SG}$.
\begin{figure}[ht]
\begin{center}
\includegraphics[width=13cm,totalheight=5.4cm]{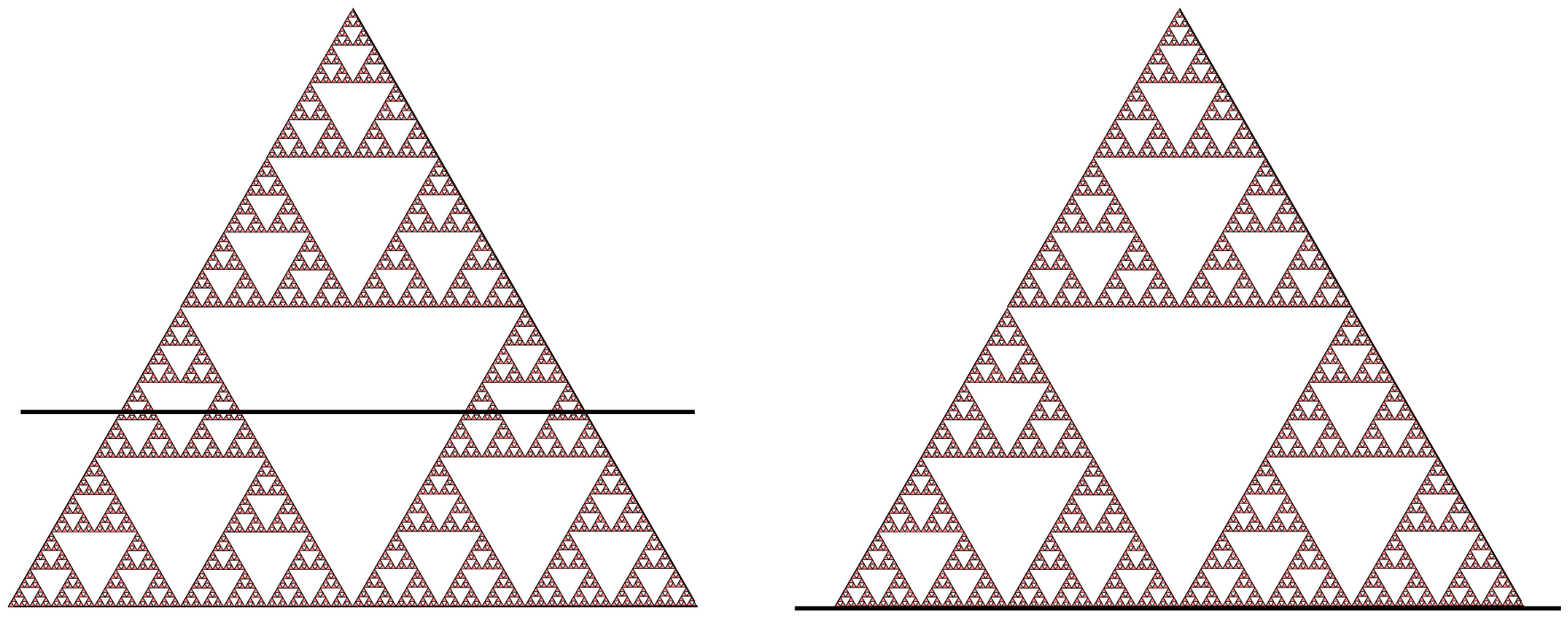}
\begin{center}
\textbf{Fig. 1.1.} \small{$\Omega_x$ and $\Omega_1$.}
\end{center}
\end{center}
\end{figure}

For simplicity, here we particularly focus on the certain domain
$\Omega_x$ which is a triangle obtained by cutting $\mathcal{SG}$
with a horizontal line at any vertical height $x$ ($0< x\leq 1$ if
we suppose that the height of $\mathcal{SG}$ is equal to $1$.) below
the top vertex $q_0$. See Fig. 1.1. An important motivation for
studying this kind of domains is that they are the simplest examples
which could serve as a testing ground for questions and conjectures
on analysis of more general fractal domains with fractal boundaries.
These domains were first introduced by  Strichartz $\cite{Str1}$ and
later studied by J. Owen and Strichartz $\cite{Owen}$, where they
gave an explicit analog of the \emph{Poisson integral formula} to
recover a harmonic function $u$ on $\Omega_x$ from its boundary
values. It is also natural to calculate an explicit \emph{Green's
function} for the Laplacian on $\Omega_x$. This was studied by Z.
Guo, R. Kogan and  Strichartz in $\cite{Guo}$ which is completely
similar to the construction of the Green's function on
$\mathcal{SG}\setminus{V_0}$ given by Kigami in
$\cite{Ki1,Ki2,Ki6}$. For some other analytic topics related to this
kind of domains, see $\cite{Hi,Jon,Ki3,Ki4}$.

In the present paper, we study the spectral properties of the
Laplacian on $\Omega_x$, which is an open problem posed in
$\cite{Owen}$. For the simplicity of description, we mainly
concentrate our attention to a particular domain $\Omega_1$ (We drop
the subscript $1$ on $\Omega$ in all that follows without causing
any confusion.) which is the complement of $\{q_0\}\cup L$, where
$L$ is the line segment joining $q_1$ and $q_2$ (in this case
$\partial \Omega=\{q_0\}\cup L$). We give a complete description of
the Dirichlet and Neumann spectra of the Laplaician on $\Omega$.

In our context, for a number of ``bad" eigenvalues (whose associated
eigenfunctions have supports touching the bottom boundary line $L$)
the spectral decimation method can not be used directly, which makes
things much more complicated. By choosing a sequence of appropriate
graph approximations, we describe a phenomenon on those eigenvalues
called \emph{weak spectral decimation} which approximates to
spectral decimation when the levels of the successive graphs go to
infinity. And we use this weak spectral decimation to replace the
role of spectral decimation in the original  Fukushima and Shima's
work $\cite{Fu}$. Actually, similarly to the standard case, weak
spectral decimation can also produce a ``weak" orbit related to the
same quadratic polynomial by an appropriate series of eigenvalues of
successive difference operators on graph approximations. We can then
trace back those ``weak" orbits to capture all the ``bad"
eigenvalues. More precisely, we classify the eigenvalues of
$-\Delta$ on $\Omega$ into three types, the \emph{localized
eigenvalues}, \emph{primitive eigenvalues} and \emph{miniaturized
eigenvalues}. The localized eigenfunctions associated to localized
eigenvalues on $\Omega$ are just a subspace of the localized
eigenfunctions on $\mathcal{SG}\setminus{V_0}$, whose supports are
disjoint from $L$. This type of eigenvalues can be dealt with in a
same way as the $\mathcal{SG}\setminus{V_0}$ case, for which the
spectral decimation can apply. The primitive and miniaturized
eigenvalues are the so-called ``bad'' eigenvalues. They are the
eigenvalues need to be paid particular attention to. We will give a
precise description of the structure of the Dirichlet and Neumann
spectra of $-\Delta$ on $\Omega$ in Section 3, before giving out the
technical proofs.

Now what happens to the asymptotic behavior of the eigenvalue
counting function $\rho^\Omega(x)$(with Dirichlet boundary condition
on $\partial\Omega$) on $\Omega$? A natural analogue of
$(\ref{abc})$ holds. Namely, there exists some positive constant
$c,C$ such that for sufficiently large $x$,
\begin{equation}\label{ssss}
 cx^{d_S/2}\leq
\rho^\Omega(x)\leq Cx^{d_S/2},
\end{equation}
 which can be proved by
first considering the asymptotic behavior of the eigenvalue counting
function for each type of eigenvalues separately, then adding up
them together. In fact, $(\ref{ssss})$ can be even easily proved
without involving the structure of the Dirichlet spectrum on
$\Omega$, as follows:
 the Dirichlet eigenvalue counting function on the top
cell $F_0(\mathcal{SG})$ is given by $\rho^{0}(x/5)$ by the
self-similarity of both the Dirichlet form and the measure $\mu$.
Therefore it follows from the minmax principle that
$\rho^{0}(x/5)\leq \rho^\Omega(x)\leq \rho^{0}(x)$, which together
with $(\ref{abc})$, also yields $(\ref{ssss})$. Moreover, the high
multiplicities of localized eigenfunctions immediately imply that
$\rho^\Omega(x)$ does not vary regularly at $\infty$, similarly to
$(\ref{sss})$. Thus,
$$
0<\liminf_{x\rightarrow\infty}\rho^\Omega(x)x^{-{d_S}/{2}}<\limsup_{x\rightarrow\infty}\rho^\Omega(x)x^{-{d_S}/{2}}<\infty.
$$

Since most eigenvalues are localized, $\rho^0(x)$ and
$\rho^\Omega(x)$ are very close. We are interested in the difference
$\rho^0(x)-\rho^\Omega(x)$. More precisely, is there some power
$\beta$ such that $\rho^0(x)-\rho^\Omega(x)\approx x^\beta$? For
this question, we have the following partial result:

\textbf{Theorem 3.10.} \emph{There exists some constant $C>0$ such
that for sufficiently large $x$,
$$
0\leq \rho^{0}(x)-\rho^\Omega(x)\leq C x^{\log2/\log 5}\log x. $$}
As a consequence, it then follows from $(\ref{ssssss})$ that
\begin{equation}\label{aaaa}
 \rho^\Omega(x)=g(\log x)x^{\log 3/\log
5}+O(x^{\log2/\log5}\log x) \quad\mbox{ as }
x\rightarrow\infty.\end{equation}

The same argument also works for the Neumann Laplacian.

Nevertheless, this should not be the entire story for the Weyl
asymptotic behavior of $\rho^\Omega(x)$. Recall that in the
classical case. Suppose $D$ is an arbitrary nonempty bounded open
set in $\mathbb{R}^n$ with smooth boundary $\partial D$, then Weyl's
classical asymptotic formula can be stated as follows:

$$\rho(x)=(2\pi)^{-n}c_n|D|_n x^{n/2}+O(x^{(n-1)/2})$$
as $x\rightarrow\infty$, where $c_n$ depends only on $n$. See
details in $\cite{Ph,Se1,Se2}$. The above remainder estimate
constitutes an important step on the way to H. Weyl's conjecture
$\cite{We}$ which states that if $\partial D$ is sufficiently
``smooth", then the asymptotic expansion of $\rho(x)$ admits a
second term, proportional to $x^{(n-1)/2}$.  Extending Weyl's
conjecture to the fractal case, M. V. Berry $\cite{Be1,Be2}$
conjectured that if $D$ has a fractal boundary $\partial D$ with
Hausdorff dimension (which later was revised into Minkowski
dimension in $\cite{Bro,La}$) $d_{\partial D}\in (n-1,n]$, then the
order of the second term should be replaced by $d_{\partial D}/2$.
See further discussion and a partial resolution of the conjectures
of Weyl and Berry in Lapidus's work $\cite{La}$. Hence it is natural
to ask that whether there is an analogue result in
$\mathcal{SG}\setminus{V_0}$ or $\Omega$ setting. For
$\mathcal{SG}\setminus{V_0}$ case,  Kigami $\cite{Ki5}$ showed that
the remainder is bounded, see $(\ref{ssssss})$. Note that this is
consistent with the fact that the boundary of $\mathcal{SG}$
consists of three points, hence has dimension zero. This was refined
by Strichartz in $\cite{Str5}$, where an exact formula was presented
with no remainder term at all, provided we restrict attention to
almost every $x$.  As for $\Omega$ case, $(\ref{aaaa})$ can be
viewed as a weak analog of the Weyl-Berry's conjecture. Moreover, We
will show that although we are unable to prove, it becomes apparent
there is a second term of order $\log2/ \log 5$ in the expansion of
the eigenvalue counting function on $\Omega$ from observing the
experimental data.

We note that our work deals with the vibrations of ``drums with
fractal membrane" since the domain itself is a fractal. The order of
the second term should has a close connection with the dimension of
the boundary $\partial \Omega$ due to Weyl-Berry's conjectures.
Moreover, when consider a more general domain $\Omega_x$, we will
meet ``drums with fractal membrane" with also fractal boundary.

The paper is organized as follows. In Section 2 we will briefly
introduce some key notions from analysis on fractals and give a
concise description of the Dirichlet and Neumann spectra of the
Laplacian for the standard $\mathcal{SG}\setminus{V_0}$ case, which
will be used in the rest of the paper.

In Section 3, we will present the exact structure of the Dirichlet
spectrum of $-\Delta$ on $\Omega$ in a self-contained and precise
way before going into the technical details. We will find an
appropriate sequence of graph approximations for the fractal domain
$\Omega$, and describe the exact structures of the discrete
Dirichlet spectra of the corresponding successive difference
operators on them. Accordingly, for each graph all the graph
eigenvalues are also divided into three types, localized, primitive
and miniaturized. By using an eigenspace dimensional counting
argument, we will show that they should make up the whole discrete
Dirichlet spectrum. We will also briefly describe how to relate the
spectra of consecutive levels and how to pass the graph
approximations to the limit by using spectral decimation for
localized eigenvalues and weak spectral decimation for other types
of eigenvalues.  We will also present analogous results for
Laplacians with Neumann boundary conditions.

In Section 4, we will describe the discrete graph primitive
Dirichlet eigenvalues on the graph approximations for each level. We
will divide our discussion into symmetric case and skew-symmetric
case. In each case, we will prove that for each level the primitive
graph eigenvalues  are exactly the total roots of a high degree
polynomial. And we will describe the weak spectral decimation
phenomenon by studying the relation between roots of consecutive
polynomials. Moreover, we will prove that for each level, the
complete discrete spectrum is made up of the three types of
eigenvalues as expected.

In Section 5, we will discuss the primitive Dirichlet  eigenvalues
of $-\Delta$ on $\Omega$ by passing the results of Section 4 on
graph approximations to the limit. Since we can only use weak
spectral decimation which is essentially based on estimates,
comparing to the $\mathcal{SG}\setminus{V_0}$ case, some trivial
results  become nontrivial and need to be proved in this section.

In Section 6, first we will prove that the whole Dirichlet spectrum
on $\Omega$ is made up of the three types of eigenvalues as
expected, following the basic idea of Fukushima and Shima's work.
Then we will give a comparison concerning the eigenvalue asymptotics
of the eigenvalue counting functions between
$\mathcal{SG}\setminus{V_0}$ case and $\Omega$ case.

In Section 7, we will give a brief discussion on how to deal with
the Neumann spectrum. We will find a similar weak spectral
decimation for primitive eigenvalues by establishing a relation
between
 symmetric (or skew-symmetric) primitive graph eigenvalues with some
high degree polynomials, but the proof is quite different from that
in the Dirichlet case.

Then in Section 8, we will list some conjectures concerning
eigenvalue asymptotics (especially the existence of the second term
of the expansion of the eigenvalue counting function), gaps in the
ratios of consecutive eigenvalues and eigenvalue clusters, which
become apparent from observing the experimental data.

We will also give a brief discussion on how to extend our method
from $\Omega$ to $\Omega_x$ with $0<x<1$ in Section 9.

The purpose of this paper is to work out the details for one
specific example. We hope this example will provide insights which
will
 inspire future
work on a more general theory.

\section{Spectral decimation on $\mathcal{SG}\setminus{V_0}$}
First we collect some key facts from analysis on $\mathcal{SG}$ that
we need to state and prove our results. These come from Kigami's
theory of analysis on fractals, and can be found in
$\cite{Ki1,Ki2,Ki6}$. An elementary exposition can be found in
$\cite{Str0,Str2}$. The fractal $\mathcal{SG}$ will be realized as
the limit of a sequence of graphs $\Gamma_0, \Gamma_1,\cdots$ with
vertices $V_0\subseteq V_1\subseteq\cdots$. The initial graph
$\Gamma_0$ is just the complete graph on $V_0=\{q_0,q_1,q_2\}$, the
vertices of an equilateral triangle in the plane, which is
considered the boundary of $\mathcal{SG}$. See Fig. 2.1. The entire
fractal is the only $0$-cell, which has $V_0$ as its boundary. At
stage $m$ of the construction, all the cells of level $m-1$ lie in
triangles whose vertices make up $V_{m-1}$. Each cell of level $m-1$
splits into three cells of level $m$, adding three new vertices to
$V_m$.
\begin{figure}[ht]
\begin{center}
\includegraphics[width=12cm,totalheight=3.8cm]{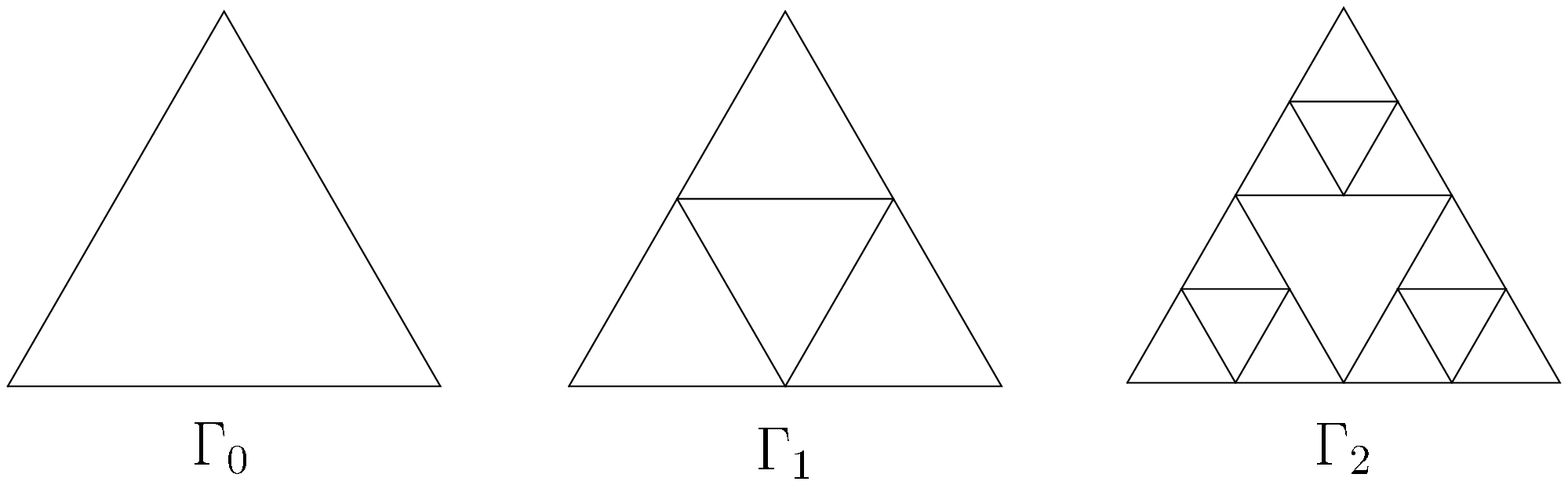}
\begin{center}
\textbf{Fig. 2.1.} \small{The first $3$ graphs, $\Gamma_0, \Gamma_1,
\Gamma_2$ in the approximations to the Sierpinski gasket.}
\end{center}
\end{center}
\end{figure}

We define the unrenormalized energy of a function $u$ on $\Gamma_m$
by
$$E_m(u):=\sum_{x\sim_m y}(u(x)-u(y))^2.$$
The energy renormalization factor is $r=\frac{3}{5}$, so the
renormalized graph energy on $\Gamma_m$ is
$$\mathcal{E}_m(u):=r^{-m}E_m(u),$$ and we can define the \emph{fractal
energy} $\mathcal{E}(u):=\lim_{m\rightarrow\infty}\mathcal{E}_m(u)$.
We define $\mathcal{F}$ as the space of continuous functions with
finite energy. Then $\mathcal{E}$ extends by polarization to a
bilinear form $\mathcal{E}(u,v)$ which serves as an inner product in
this space. The energy $\mathcal{E}$ gives rise to a natural
distance on $\mathcal{SG}$ called the \emph{effective resistance
metric} on $\mathcal{SG}$, which is defined  by
\begin{equation}\label{2}
d(x,y):=(\min\{\mathcal{E}(u): u(x)=0 \mbox{ and } u(y)=1\})^{-1}
\end{equation}
for $x,y\in \mathcal{SG}$.
 It is known that $d(x,y)$ is bounded above and below by constant
multiples of $|x-y|^{\log(5/3)/\log 2}$, where $|x-y|$ is the
Euclidean distance. Furthermore, the definition $(\ref{2})$ implies
that functions on $\mathcal{F}$ are H\"{o}lder continuous of order
$\frac{1}{2}$ in the effective resistance metric.

 We let $\mu$ denote the standard probability measure on
$\mathcal{SG}$ that assigns the measure $3^{-m}$ to each cell of $m$
level. The standard Laplacian may then be defined using the weak
formulation: $u\in dom\Delta$ with $-\Delta u=f$ if $f$ is
continuous, $u\in \mathcal{F}$, and
\begin{equation}\label{1}
\mathcal{E}(u,v)=\int_{\mathcal{SG}} fvd\mu
\end{equation}
 for all $v\in
\mathcal{F}_0$, where $\mathcal{F}_0=\{v\in \mathcal{F}:
v|_{V_0}=0\}$. There is also a pointwise formula (which is proven to
be equivalent in $\cite{Str2}$) which, for nonboundary points in
$V_*=\bigcup_m V_m$ (not in $V_0$) computes
$$\Delta u(x)=\frac{3}{2}\lim_{m\rightarrow\infty}5^m\Delta_m u(x),$$
where $\Delta_m$ is a discrete Laplacian associated to the graph
$\Gamma_m$, defined by
$$\Delta_m u(x):=\sum_{y\sim_m x}(u(y)-u(x))$$ for $x$ not on
the boundary.

The Laplacian satisfies the scaling property
$$\Delta(u\circ F_i)=\frac{1}{5}(\Delta u)\circ F_i$$
and by iteration
$$\Delta(u\circ F_w)=\frac{1}{5^m}(\Delta u)\circ F_w$$ for $F_w=F_{w_1}\circ F_{w_2}\circ\cdots\circ
F_{w_m}$.

Although there is no satisfactory analogue of gradient, there is
\emph{normal derivative} $\partial_n u(q_i)$ defined at boundary
points by
$$\partial_n u(q_i):=\lim_{m\rightarrow\infty}\sum_{y\sim_m q_i}r^{-m}(u(q_i)-u(y)),$$
the limit existing for all $u\in dom\Delta$. The definition may be
localized to boundary points of cells. For each point $x\in
V_m\setminus V_0$, there are two cells containing $x$ as a boundary
point, hence two normal derivatives at $x$. For $u\in dom\Delta$,
the normal derivatives at $x$ satisfy the \emph{matching condition}
that their sum is zero. The matching condition allows us to glue
together local solutions to $-\Delta u=f$.

The above matching condition property follows easily from a local
version of the following \emph{Gauss-Green formula}, which is an
extension of $(\ref{1})$ to the case when $v$ doesn't vanish on the
boundary:
$$\mathcal{E}(u,v)=\int_{\mathcal{SG}}(-\Delta u)vd\mu+\sum_{V_0}v\partial_n u.$$
The local version of the Gauss-Green formula is
$$\mathcal{E}_A(u,v)=\int_A(-\Delta u)vd\mu+\sum_{\partial
A}v\partial_nu,$$ where $A$ is any finite union of cells and
$\mathcal{E}_A(u,v)$ is the restriction of the energy bilinear form
$\mathcal{E}(u,v)$ to $A$, which can also be defined directly by
$$\mathcal{E}_A(u,v):=\lim_{m\rightarrow\infty}\sum_{x\sim_m y\atop in A}(u(x)-u(y))(v(x)-v(y)).$$

Now we come to a brief recap of the spectral decimation on
$\mathcal{SG}$. Our goal is to find all solutions of the eigenvalue
equation
$$-\Delta u=\lambda u\quad \mbox{ on } \mathcal{SG}\setminus{V_0}$$ as
limits of solutions of the discrete version
$$-\Delta_m u_m=\lambda_m u_m \quad \mbox{ on } V_m\setminus V_0.$$ In
the $\mathcal{SG}\setminus{V_0}$ case, we are lucky that we may take
$u_m=u|_{V_m}$, which is necessarily convenient for the spectral
decimation. We should emphasize that this is not true for $\Omega$
case.

 The method of
spectral decimation on $\mathcal{SG}$ was invented by  Fukushima and
Shima ${\cite{Fu}}$  to relate eigenfunctions and eigenvalues of the
discrete Laplacian $-\Delta_m$'s on the graph approximation
$\Gamma_m$'s for different values of $m$ to each other and the
eigenfunctions and eigenvalues of the fractal Laplacian $-\Delta$ on
$\mathcal{SG}\setminus{V_0}$. In essence, an eigenfunction on
$\Gamma_m$ with eigenvalue $\lambda_m$ can be extended to an
eigenfunction on $\Gamma_{m+1}$ with eigenvalue $\lambda_{m+1}$,
where $\lambda_m=f(\lambda_{m+1})$ for an explicit function $f$
defined by
\begin{equation}\label{40}
f(x):=x(5-x),
\end{equation} except for certain specified
\emph{forbidden eigenvalues}, and all eigenfunctions on
$\mathcal{SG}\setminus{V_0}$ arise as limits of this process
starting at some level $m_0$ which is called the \emph{generation of
birth}. This is true regardless of the boundary conditions, but if
we specify Dirichlet or Neumann boundary condition we can describe
explicitly all eigenspaces and their multiplicities.

Denote the real valued inverse functions of $f(x)$ by
$\phi_{\pm}(x)$. That is
\begin{equation}\label{4}
\phi_{\pm}(x):=\frac{5\pm\sqrt{25-4x}}{2}.
\end{equation}

We describe the procedure briefly here. First, there is a
\emph{local extension algorithm} that shows how to uniquely extend
an eigenfunction $u_m$ defined on $V_m$ to a function defined on
$V_{m+1}$ such that the $\lambda$-eigenvalue equations hold on all
points of $V_{m+1}\setminus V_m$.   For $\mathcal{SG}$, the
extension algorithm is: Suppose $u_m$ is an eigenfunction on
$\Gamma_m$ with eigenvalue $\lambda_m$. Let
$\lambda_{m+1}=\phi_{\pm}(\lambda_m)$. Consider an $m$-cell with
boundary points $x_0,x_1,x_2$ and let $y_0,y_1,y_2$ denote the
points in $V_{m+1}\setminus V_m$ in that cell, with $y_i$ opposite
$x_i$. Extend $u_m$ to a function $u_{m+1}$ on $V_{m+1}$ by defining
(for simplicity of notation, we drop the subscripts on $u$)
\begin{equation}\label{3}
u(y_i)=\frac{(4-\lambda_{m+1})((u(x_{i+1})+u(x_{i-1})))+2u(x_i)}{(2-\lambda_{m+1})(5-\lambda_{m+1})},
\quad i=0,1,2.
\end{equation}
Then we have the following proposition  taken from $\cite{Str2}$.

\textbf{Proposition 2.1.} \emph{ Suppose $\lambda_{m+1}\neq 2,5$ or
$6$, and $\lambda_m=f(\lambda_{m+1})$. If $u_m$ is a
$\lambda_m$-eigenfunction of $-\Delta_m$ and is extended to a
function $u_{m+1}$ on $V_{m+1}$ by $(\ref{3})$, then $u_{m+1}$ is a
$\lambda_{m+1}$-eigenfunction of $-\Delta_{m+1}$, Conversely, if
$u_{m+1}$ is a $\lambda_{m+1}$-eigenfunction of $-\Delta_{m+1}$ and
is restricted to a function $u_m$ on $V_m$, then $u_m$ is a
$\lambda_m$-eigenfunction of $-\Delta_m$.}

The forbidden eigenvalues $\{2,5,6\}$ are singularities of the
spectral decimation function $f$. It is ``forbidden" to decimate to
a forbidden eigenvalue. Because forbidden eigenvalues have no
predecessor, we speak of forbidden eigenvalues being ``born" at a
level of approximation $m$.

Next we want to take the limit as $m\rightarrow\infty$. We assume
that we have an infinite sequence $\{\lambda_{m}\}_{m\geq m_0}$
related by $\lambda_{m+1}=\phi_{\pm}(\lambda_m)$ with all but a
finite number of $\phi_{-}$'s. Then we may define
$$\lambda:=\frac{3}{2}\lim_{m\rightarrow\infty}5^m\lambda_m.$$
It is easy to see that the limit exists since
\begin{equation}\label{5}
\phi_{-}(x)=\frac{1}{5}x+O(x^2)
 \end{equation}
 as $x\rightarrow 0$.
We start with a $\lambda_{m_0}$-eigenfunction $u$ of $-\Delta_{m_0}$
on $V_{m_0}$, and extend $u$ to $V_*$ successively using
$(\ref{3})$, assuming that none of $\lambda_m$ is a forbidden
eigenvalue. Since $(\ref{5})$ implies $\lambda_m=O(\frac{1}{5^m})$
as $m\rightarrow\infty$, it is easy to see that $u$ is uniformly
continuous on $V_*$ and so extends to a continuous function on
$\mathcal{SG}$. Moreover, it satisfies the $\lambda$-eigenvalue
equation for $-\Delta$.

 A proof in $\cite{Fu}$ guarantees that this
spectral decimation produces all possible eigenvalues and
eigenfunctions (up to linear combination).

To describe the explicit Dirichlet and Neumann spectra, we have to
describe all possible generations of birth and values for
$\lambda_{m_0}$, and describe the multiplicity of the eigenvalue by
giving an explicit basis for the $\lambda_{m_0}$-eigenspace of
$-\Delta_{m_0}$. For each $m$, we have to add up the dimensions of
eigenspaces with generation of birth $m_0\leq m$, extended to
$\Gamma_m$ in all allowable ways. This total must be $\sharp V_m$
(Neumann) or $\sharp V_m-3$ (Dirichlet), the dimension of the space
on which the symmetric operator $-\Delta_m$ acts.  Now we give a
brief description of the structure of the Dirichlet and Neumann
spectra of $-\Delta$ on $\mathcal{SG}\setminus{V_0}$ respectively.

\textbf{Dirichlet spectrum.}

We denote by $\mathcal{D}$ the Dirichlet spectrum of $-\Delta$ on
$\mathcal{SG}\setminus{V_0}$ and by $\mathcal{D}_m$ the discrete
Dirichlet spectrum of $-\Delta_m$ on $\Gamma_m$ for $m\geq 1$. Due
to the above discussion, we only need to make clear the spectrum
$\mathcal{D}_m$ for each level $m$. There are two kinds of
eigenvalues, \emph{initial} and \emph{continued}. The continued
eigenvalues will be those that arise from eigenvalues of
$\mathcal{D}_{m-1}$ by the spectral decimation. Those that remain,
the initial eigenvalues, must be some of the forbidden eigenvalues
by Proposition 2.1.

In $\cite{Shi1}$, it is proved that $\mathcal{D}_1$ consists of two
eigenvalues $2$ and $5$ with multiplicities 1 and 2 respectively,
and for $m\geq 2$, the only possible initial eigenvalues in
$\mathcal{D}_{m}$ are the two forbidden eigenvalues $5$ and $6$ with
multiplicities $\frac{3^{m-1}+3}{2}$ and $\frac{3^m-3}{2}$
respectively. Hence we may classify eigenvalues into three series,
which we call the $2$-series, $5$-series, and $6$-series, depending
on the value of $\lambda_{m_0}$. The eigenvalues in the $2$-series
all have multiplicity $1$, while the eigenvalues in the other series
all exhibit higher multiplicity. Also, if $\lambda$ is an eigenvalue
in the $5$-series or $6$-series, then $5^m\lambda$ is also an
eigenvalue, corresponding to a generation of birth $m_0+m$, with the
same choice of $\phi_{\pm}$ relations (suitably reindexed).

\textbf{Neumann spectrum.}

We impose a Neumann condition on the graph $\Gamma_m$ by imagining
that it is embedded in a larger graph by reflecting in each boundary
vertex and imposing the $\lambda_m$-eigenvalue equation on the even
extension of $u$. This just means that we impose the equation
$$(4-\lambda_m)u(q_i)=2u(F_i^mq_{i+1})+2u(F_i^{m}q_{i-1})$$ at $q_i$
for $i=0,1,2$. Then the Neumann $\lambda_m$-eigenvalue equations
consist of exactly $\sharp V_m$ equations in $\sharp V_m$ unknowns.
Similarly to the Dirichlet case, we also only need to make clear all
the discrete spectra. The result is very similar to  the Dirichlet
spectrum, with only a few changes. We omit it.

It should be emphasized here that those eigenfunctions which are
simultaneously Dirichlet and Neumann play an important role in the
spectral analysis of $\mathcal{SG}$. Here we call them
\emph{localized eigenfunctions} since all of them have small
supports. (Here this definition of localized eigenfunctions is
slightly different from that of $\cite{Ba,Ki5,Str2}$ for the
convenience of further discussion for $\Omega$ case.) Similarly to
$\mathcal{D}$, to describe the structure of localized
eigenfunctions, we only need to make clear the structure of all
initial localized eigenvalues,  which consists of $5$-series and
$6$-series eigenvalues. In fact, the multiplicity of a $5$-series
eigenvalue with generation of birth $m$, is
$\rho_{m}(5):=\frac{3^{m-1}-1}{2}$ with an eigenfunction associated
to each $m$-level loop (a $m$-level circuit around an empty
upside-down triangle in the graph $\Gamma_m$).  The eigenfunction
$u$ associated to each loop takes value $0$ on all $m$-level points
 not lying in that loop. Moreover, the support of $u$ is exactly the union of all $m$-cells intersecting that
loop. The multiplicity of a $6$-series eigenvalue with generation of
birth $m$, is $\rho_{m}(6):=\frac{3^m-3}{2}$ with an eigenfunction
associated to each point $x$ in $V_{m-1}\setminus V_0$.  Each such
eigenfunction $u$ takes value $0$ on all points in $V_{m-1}$ except
$x$. Moreover, $u$ is supported in the union of two $(m-1)$-level
cells containing $x$. The existence of localized eigenfunctions is
unprecedented in all of smooth mathematics. However, for a class of
p.c.f. self-similar sets, including $\mathcal{SG}$, localized
eigenfunctions dominate global eigenfunctions. See more details in
$\cite{Ba,Ki5}$. See also Section 4 of Kigami's book $\cite{Ki6}$,
where most results are explained in detail.

\section{The structures of Dirichlet spectrum and Neumann spectrum on $\Omega$}
To give the readers an intuitive perception of the structure of the
spectrum of $-\Delta$ on $\Omega$ in advance, in this section we
describe all Dirichlet and Neumann eigenvalues and eigenfunctions on
$\Omega$ avoiding involving technical proofs. We will go to the
details in the remaining sections.

\subsection{Dirichlet spectrum}

We begin with the Dirichlet case. First we formulate the eigenvalue
problem of $-\Delta$ on $\Omega$ with Dirichlet boundary
condition(for short, the \emph{Dirichlet Laplacian}).

\textbf{Definition 3.1.} \emph{Let $\mathcal{F}_\Omega:=\{u\in
\mathcal{F}: u|_{\partial\Omega}=0\}$. The Dirichlet Laplacian
$\Delta_D$ on $\Omega$ with domain $\mathcal{D}[\Delta_D]$ is
formulated as follows: for $u\in \mathcal{F}_\Omega$ and $f\in
L^2(\Omega,\mu|_\Omega)$,
$$u\in\mathcal{D}[\Delta_D]\mbox{ and } -\Delta_D u=f \mbox {\quad\quad if and only if\quad\quad } \mathcal{E}(u,v)=\int_\Omega fvd\mu\mbox{ for any } v\in \mathcal{F}_\Omega.$$
}

If we replace $\Omega$ by $\mathcal{SG}\setminus V_0$ in the above
definition, then we get the standard Dirichlet Laplacian which is
introduced in $\cite{Ki6}$.

\textbf{Definition 3.2.} \emph{For $\lambda\in \mathbb{R}$ and
$u\in\mathcal{D}[\Delta_D]$ if
$$-\Delta_D u=\lambda u,$$
then $\lambda$ is called an  eigenvalue of $-\Delta_D$  on $\Omega$
(or, a Dirichlet eigenvalue of $-\Delta$ on $\Omega$), and $u$ is
called an associated (Dirichlet) eigenfunction.}

 Let  $\mathcal{S}$ denote the
 spectrum of $-\Delta_D$ on $\Omega$ ($\mathcal{S}$ is also called the Dirichlet spectrum of $-\Delta$ on $\Omega$).  We will consider three kinds of
Dirichlet eigenfunctions, \emph{localized}, \emph{primitive}, and
\emph{miniaturized} eigenfunctions.

In the following, we will always use $u$ to denote an eigenfunction
of $-\Delta_D$ on $\Omega$ and $\lambda$ to denote the associated
eigenvalue of $u$.

\textbf{Definition 3.3.} \emph{$u$ is called a localized
eigenfunction if it is a localized eigenfunction on
$\mathcal{SG}\setminus{V_0}$ whose support is disjoint from $L$ (the
line segment joining $q_1$ and $q_2$). }

The associated eigenvalue $\lambda$ is called a \emph{localized
eigenvalue}. Denote by  $\mathcal{L}$ the set consisting of all such
 eigenvalues. Obviously, all the eigenvalues in $\mathcal{L}$ have
generation of birth $m_0\geq 3$ (the ones with $m_0=2$ all have
supports intersecting $L$) and $\lambda_{m_0}=5$ or $6$.

Comparing to the $\mathcal{SG}\setminus{V_0}$ case, instead of the
eigenfunctions associated to the $2$-series eigenvalues, there is
also a type  of global eigenfunctions in $\Omega$ case,  which will
be sorted into symmetric and skew-symmetric parts according to the
reflection symmetry fixing $q_0$.

\textbf{Definition 3.4.} \emph{$u$ is called a symmetric primitive
eigenfunction if it is symmetric under the reflection symmetry
fixing $q_0$ and also local symmetric in each cell
$F_w(\mathcal{SG})$ under the reflection symmetry fixing $F_wq_0$
with word $w$ taking symbols only from $\{1,2\}$.}

 Fig. 3.1. gives a symbolic picture of the above mentioned symmetries,
indicated by dotted lines. The associated eigenvalue $\lambda$ is
called a \emph{symmetric primitive eigenvalue}. Denote by
$\mathcal{P}^+$ the set consisting of all such eigenvalues.

Similarly,

\textbf{Definition 3.5.} \emph{If $u$ is skew-symmetric under the
reflection symmetry fixing $q_0$,  but still local symmetric in
small cells, then it is called a skew-symmetric eigenfunction. }

The associated eigenvalue $\lambda$ is called a \emph{skew-symmetric
eigenvalue}. Denote by  $\mathcal{P}^-$ the set consisting of all
such eigenvalues.

Both the symmetric and skew-symmetric primitive eigenfunctions are
called \emph{primitive eigenfunctions}. All the associated
eigenvalues are called \emph{primitive eigenvalues}. Let
$\mathcal{P}$  denote the set consisting of all of them. Namely,
$$\mathcal{P}=\mathcal{P}^+\cup\mathcal{P}^-.$$

The primitive eigenfunction $u$ (either the symmetric or
skew-symmetric case) is uniquely determined by the values denoted by
$(b_0,b_1,b_2,\cdots)$ of $u$ on vertices $(q_0, F_1q_0,
F_{1}^2q_0,\cdots)$ by using the eigenfunction extension algorithm
described in (\ref{3}). Due to the Dirichlet boundary condition,
$\forall\lambda\in\mathcal{P}$, for the associated eigenfunction $u$
of $\lambda$, we always have $b_0=0$ and
$\lim_{m\rightarrow\infty}b_m=0$. We call $(q_0, F_1q_0,
F_{1}^2q_0,\cdots)$ a \emph{skeleton} of $\Omega$ since it plays a
critical role in the study of primitive eigenfunctions.

\textbf{Theorem 3.1.} \emph{All the  primitive eigenvalues are of
 multiplicity $1$.}

This theorem will be proved in Section 6.

\begin{figure}[ht]
\begin{center}
\includegraphics[width=6.2cm,totalheight=5.6cm]{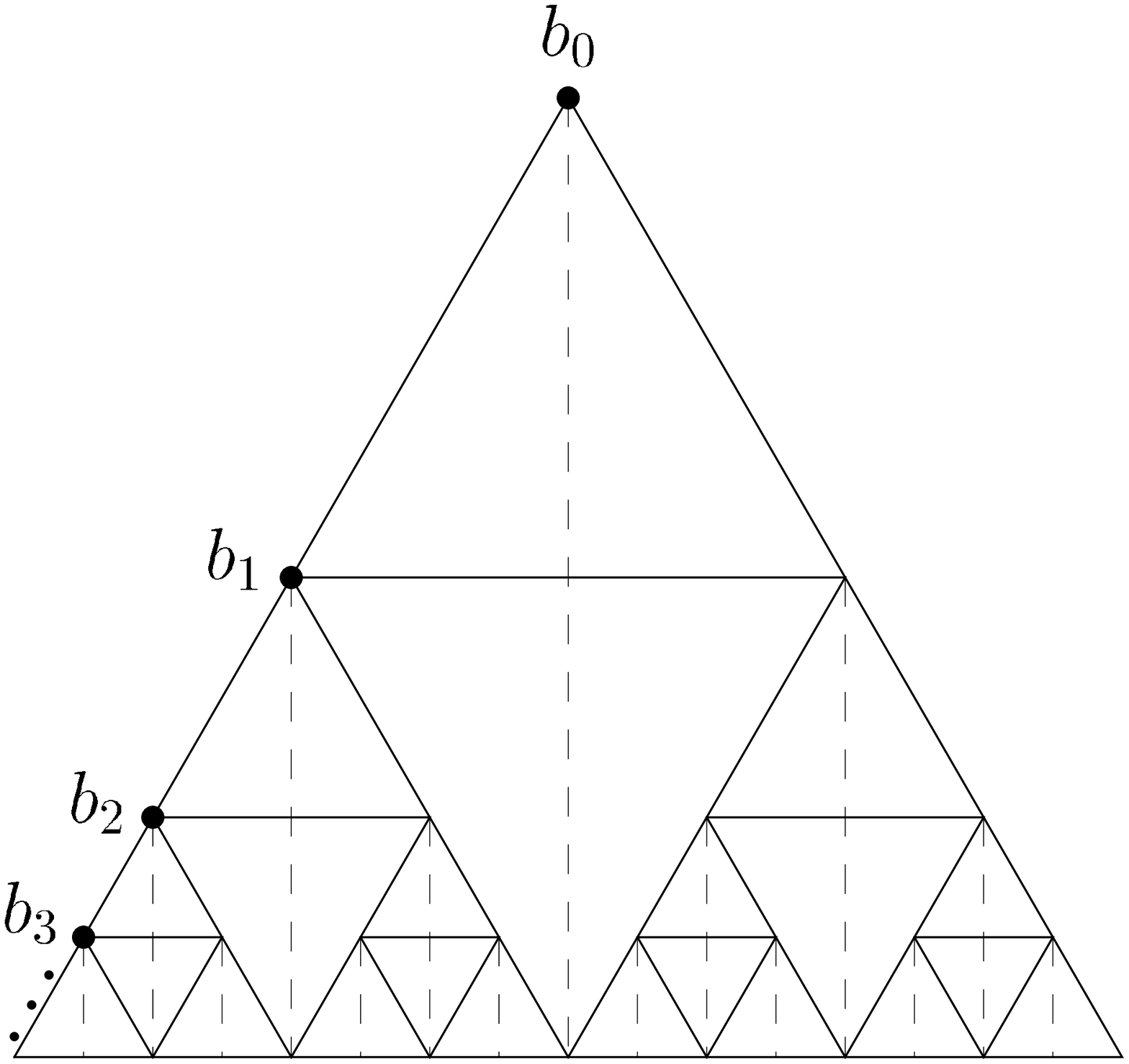}
\begin{center}
\textbf{Fig. 3.1.} \small{The first $4$ level symmetries and the
skeleton of $\Omega$.}
\end{center}
\end{center}
\end{figure}

The following argument will show that there is another type of
eigenfunctions. For each skew-symmetric eigenvalue $\lambda\in
\mathcal{P}^-$, there is a family of eigenfunctions with eigenvalue
$5^k\lambda$ and multiplicity $2^k$ for $k=1,2,3,\cdots$. To get
such an eigenfunction, just take the $\lambda$-eigenfunction $u$,
contract it $k$ times, place it in any one of the $2^k$ bottom cells
of level $k$, and take value $0$ elsewhere. See Fig 3.2. The reason
we can do this is that on the boundary point $q_0$, $u(q_0)=0$ and
$\partial_n u(q_0)=0$ which make the matching condition holds
automatically.

\textbf{Definition 3.6.} \emph{We call all the above obtained
eigenfunctions miniaturized eigenfunctions. }

 If
$u$ is a miniaturized eigenfunction obtained by contracting a
skew-symmetric primitive eigenfunction $k$ times, then we call $u$ a
\emph{$k$-contracted miniaturized eigenfunction.} Let $\mathcal{M}$
denote all the eigenvalues associated to them. Obviously,
$\mathcal{M}$ is determined by $\mathcal{P}^-$.
\begin{figure}[ht]
\begin{center}
\includegraphics[width=9cm,totalheight=4.2cm]{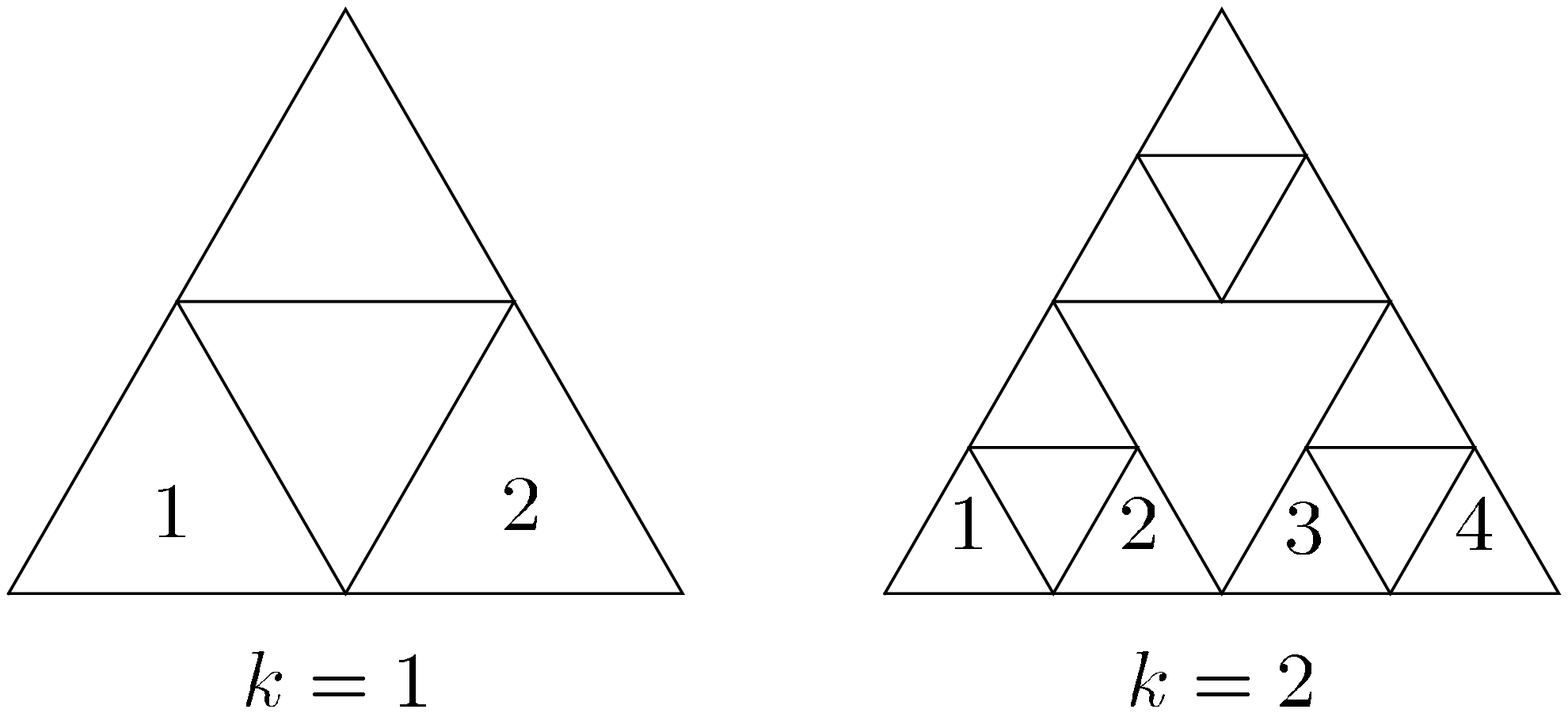}
\begin{center}
\textbf{Fig. 3.2.} \small{The first $2$ level miniaturized
eigenfunctions.}
\end{center}
\end{center}
\end{figure}

We will prove in Section 6 that all the eigenfunctions of
$-\Delta_D$ on $\Omega$  fall into one of these three types, and
there are no coincidences of eigenvalues among different types.

For the purpose of studying the exact structure of the spectrum
$\mathcal{S}$, the first thing we should consider is to describe the
Dirichlet spectra of $-\Delta_m$'s on the associated graph
approximations. Similarly to the $\mathcal{SG}\setminus{V_0}$ case,
the fractal domain $\Omega$ can be realized as the limit of a
sequence of graphs $\Omega_m$. More precisely, $\forall m\geq 1$,
let $V_m^\Omega$ be a subset of $V_m$ with all vertices lying along
$L$ removed. Let $\Omega_m$ be the subgraph of $\Gamma_m$ restricted
to $V_m^\Omega$. Denote by $\partial\Omega_m$ the boundary of the
finite graph $\Omega_m$. It is easy to find that
$V_m^\Omega\setminus\partial\Omega_m$ and $\partial\Omega_m$
approximate to $\Omega$ and $\partial\Omega$ as $m$ goes to infinity
respectively. See Fig. 3.3.

A routing argument shows that the Dirichlet Laplaician $\Delta_D$
could be viewed as the limit of suitably scaled graph Laplaicians
$\Delta_m$ on $\Omega_m$, as is done in
$\cite{Ki1,Ki2,Ki6,Str0,Str2}$ for the standard
$\mathcal{SG}\setminus{V_0}$ case. Hence, there is also a pointwise
formula which, for nonboundary points in $V_*\cup \Omega$, computes
$$\Delta u(x)=\frac{3}{2}\lim_{m\rightarrow\infty}5^m\Delta_m u(x),$$
where $\Delta_m$ is a discrete Laplacian associated to the graph
$\Omega_m$, defined by
$$\Delta_m u(x)=\sum_{y\sim_m x}(u(y)-u(x))$$ for $x$ in $V_m^\Omega\setminus\partial\Omega_m$.

We denote by $\mathcal{S}_m$ the discrete Dirichlet spectrum of
$-\Delta_m$ on $\Omega_m$ for $m\in \mathbb{N}$. On $\Omega_m$ the
Dirichlet $\lambda_m$-eigenvalue equations consist of exactly
$\sharp (V_m^\Omega\setminus\partial\Omega_m)$ equations in $\sharp
(V_m^\Omega\setminus\partial\Omega_m)$ unknowns.
 We start from m=2 since there is no Dirichlet $\lambda_1$-eigenvalue
equation. For simplicity, let $a_m=\sharp
(V_m^\Omega\setminus\partial\Omega_m)$.  It is easy to check that
$a_2=5$, $a_3=24$, and more generally,

\textbf{Proposition 3.1.} \emph{$a_m=\frac{3^{m+1}-1}{2}-2^{m+1}$,
$\forall m\in \mathbb{N}$.}

\emph{Proof.} Notice that $a_m=a_{m-1}+3^m+2^{m-1}-3\cdot 2^{m-1}$,
where $3^m=\sharp(V_m\setminus V_{m-1})$, $2^{m-1}$ is the number of
points lying on the bottom boundary of $\Omega_{m-1}$, and $3\cdot
2^{m-1}$ is the number of points in $V_m\setminus V_{m-1}$ lying on
$L$ or $\partial\Omega_m$. $\Box$
\begin{figure}[ht]
\begin{center}
\includegraphics[width=14.5cm,totalheight=4.6cm]{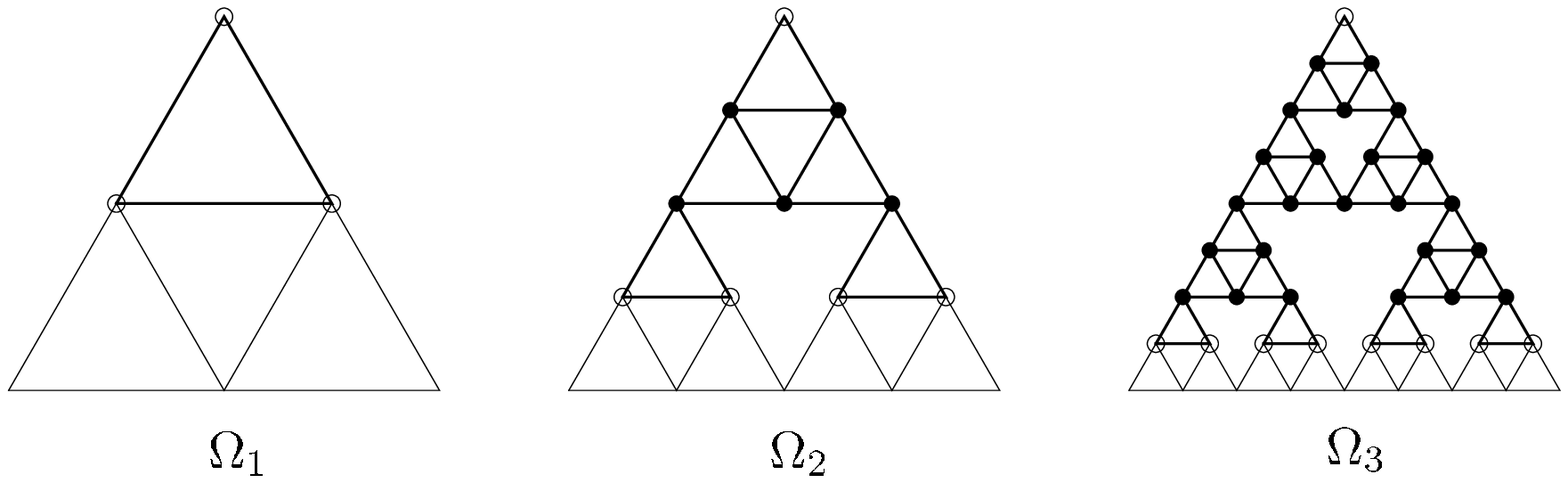}
\begin{center}
\textbf{Fig. 3.3.} \small{The first $3$ graphs, $\Omega_1, \Omega_2,
\Omega_3$ in the approximations to $\Omega$ with  inside points and
boundary points represented by dots and circles respectively.}
\end{center}
\end{center}
\end{figure}

Due to different types of eigenvalues of $-\Delta_D$ on $\Omega$, we
should consider the associated different types of graph Dirichlet
eigenvalues of $-\Delta_m$.  We now describe how to define
$\mathcal{L}_m$, $\mathcal{P}_m$ and $\mathcal{M}_m$ respectively.
For simplicity, in the rest of this subsection, without causing any
confusion, we omit the term ``Dirichlet" for graph eigenfunctions
and eigenvalues.

In the following, we will always use $u_m$ to denote an
eigenfunction of $-\Delta_m$ on $\Omega_m$ satisfying the Dirichlet
boundary condition, and $\lambda_m$ to denote the associated
eigenvalue of $u_m$.

In fact, by the spectral decimation recipe, each localized
eigenfunction $u$ of $-\Delta_D$ on $\Omega$ whose generation of
birth $m_0\leq m$ can be restricted to $\Omega_m$ to get a graph
eigenfunction $u_m$ of $-\Delta_m$, with the Dirichlet boundary
condition of $u_m$ on $\partial\Omega_m$ holding automatically.

\textbf{Definition 3.7.} \emph{Let $u$ be a localized eigenfunction
of $-\Delta_D$ on $\Omega$ with generation of birth $m_0\leq m$,
then its restricted graph function $u_m$ on $\Omega_m$ is called a
$m$-level localized graph eigenfunctions of $-\Delta_m$ on
$\Omega_m$. }

All the associated eigenvalues are called $m$-level \emph{localized
graph eigenvalues}. We use $\mathcal{L}_m$ to denote the set
consisting of all these type eigenvalues.

 We
can not imitate the above process to get the $m$-level primitive
graph  eigenfunctions since the Dirichlet boundary condition would
be destroyed if we do the similar restriction. But we can define
$m$-level primitive graph eigenfunctions on $\Omega_m$ directly in
the following way.

\textbf{Definition 3.8.}\emph{ A Dirichlet eigenfunction $u_m$ on
$\Omega_m$ is called a $m$-level symmetric primitive graph
eigenfunction if it is symmetric under the reflection symmetry
fixing $q_0$ and also local symmetric in $F_w(\mathcal{SG})\cap
V_m^\Omega$ under the reflection symmetry fixing $F_wq_0$ for each
word $w$ taking symbols only from $\{1,2\}$.}

The associated eigenvalue $\lambda_m$ is called a $m$-level
\emph{symmetric primitive graph eigenvalue}. Denote by
$\mathcal{P}_m^+$ the set of all this type of  eigenvalues.

Similarly,

\textbf{Definition 3.9.} \emph{If $u_m$ is skew-symmetric under the
reflection symmetry fixing $q_0$,  but still local symmetric in
small cells, then it is called a $m$-level skew-symmetric primitive
graph eigenfunction. }

$\mathcal{P}_m^-$ can be defined in a similar way. Let
$\mathcal{P}_m$ denote all the $m$-level primitive graph
eigenvalues. Namely,
$$\mathcal{P}_m=\mathcal{P}_m^+\cup\mathcal{P}_m^-.$$

Similarly to the limit case, the primitive graph eigenfunction $u_m$
(either the symmetric or skew-symmetric case) is uniquely determined
by the values denoted by $(b_0,b_1,b_2,\cdots,b_m)$ of $u_m$ on
vertex points $(q_0, F_1q_0, F_{1}^2q_0,\cdots,F_{1}^mq_0)$ by using
the eigenfunction extension algorithm  (\ref{3}). Due to the
Dirichlet boundary condition, we always have $b_0=b_m=0$ for a
$m$-level primitive graph  eigenfunction $u_m$. We call $(q_0,
F_1q_0, F_{1}^2q_0,\cdots,F_{1}^mq_0)$ a \emph{skeleton} of
$\Omega_m$. It also plays a critical role in the study of $m$-level
primitive graph eigenfunctions.

 Miniaturized graph eigenfunctions on
$\Omega_m$ can be defined in a similar way by using miniaturization
of skew-symmetric primitive graph eigenfunctions whose level
strictly less than $m$.

\textbf{Definition 3.10.} \emph{For a Dirichlet eigenfunction $u_m$
on $\Omega_m$, if there exists an integer $m'<m$ and a $m'$-level
skew-symmetric primitive graph eigenfunction $u_{m'}$ such that
after contracting $u_{m'}$ $m-m'$ times, placing it in one of the
$2^{m-m'}$ bottom copies of $\Omega_{m'}$ in $\Omega_m$, and taking
value $0$ elsewhere, one can obtain $u_m$, then $u_m$ is called a
$m$-level miniaturized graph eigenfunction.}

The associated eigenvalue $\lambda_m$ is called a $m$-level
\emph{miniaturized graph eigenvalue}. $m'$ is called the \emph{type}
of $\lambda_m$. Denote by $\mathcal{M}_m$ the set of all such
eigenvalues. Obviously, $\mathcal{M}_m$ is determined by all
$\mathcal{P}_k^-$'s with $k<m$.

It is not difficult to make clear all the localized graph
eigenvalues in $\mathcal{L}_m$, since they are almost the same as
the $\mathcal{SG}\setminus{V_0}$ case. There are two kinds of
eigenvalues in $\mathcal{L}_m$, \emph{initial} and \emph{continued.}

\textbf{Theorem 3.2.} \emph{Let $m\geq 2$, then
$\sharp\mathcal{L}_m=\frac{3^{m+1}-1}{2}-(m-2)\cdot2^m-26\cdot
2^{m-3}$. Moreover, the initial eigenvalues in $\mathcal{L}_m$ are
$5$ and $6$ with multiplicity
$\rho_m^\Omega(5):=\frac{3^{m-1}+1}{2}-2^{m-1}$ and
$\rho_m^\Omega(6):=\frac{3^m-1}{2}-2^m$ respectively.}

\emph{Proof.} Similarly to the $\mathcal{SG}\setminus{V_0}$ case,
the initial eigenvalues are $5$ and $6$. For the $6$-eigenfunctions
of $-\Delta_m$ on $\Omega_m$, comparing to the $6$-eigenfunctions of
$-\Delta_m$ on $\Gamma_m$, the only difference is those
eigenfunctions whose support intersecting the boundary $\partial
\Omega_m$ should be removed. A similar analysis shows that they are
indexed by points in $V_{m-1}^\Omega\setminus\partial\Omega_{m-1}$.
Hence the multiplicity of $6$ is
$$\rho_m^\Omega(6)=a_{m-1}=\frac{3^m-1}{2}-2^m.$$ Similarly, the
$5$-eigenfunctions of $-\Delta_m$ on $\Omega_m$  are indexed by
$m$-level loops except those loops touching $\partial\Omega_m$.
Hence the multiplicity of $5$ is
$$\rho_m^\Omega(5)=\rho_m(5)-(1+2+2^2+\cdots+2^{m-2})=\frac{3^{m-1}+1}{2}-2^{m-1}.$$

The continued eigenvalues will be those that arise from eigenvalues
of $\mathcal{L}_{m-1}$ by the spectral decimation. Note that every
eigenvalue $\lambda_{m-1}$ of $-\Delta_{m-1}$ bifurcates into two
choices of $\lambda_m$  of $-\Delta_m$ by $(\ref{4})$, except
$\lambda_{m-1}=6$, which just yields the single choice $\lambda_m=3$
since the other is a forbidden eigenvalue $2$. We know that
$\rho_{m-1}^\Omega(6)$ of $\mathcal{L}_{m-1}$ correspond to
eigenvalue $6$ of $-\Delta_{m-1}$, while the remaining
$\sharp\mathcal{L}_{m-1}-\rho_{m-1}^\Omega(6)$ of them correspond to
other eigenvalues, leading to a space of continued eigenfunctions of
dimension
$2\cdot(\sharp\mathcal{L}_{m-1}-\rho_{m-1}^\Omega(6))+\rho_{m-1}^\Omega(6)=2\cdot\sharp\mathcal{L}_{m-1}-\frac{3^{m-1}-1}{2}+2^{m-1}$.
If we add to this $\rho_m^\Omega(6)=\frac{3^m-1}{2}-2^m$ and
$\rho_m^\Omega(5)=\frac{3^{m-1}+1}{2}-2^{m-1}$, we should obtain
$\sharp \mathcal{L}_m$. Hence we have
\begin{eqnarray*}
\sharp
\mathcal{L}_m&=&2\cdot\sharp\mathcal{L}_{m-1}-\frac{3^{m-1}-1}{2}+2^{m-1}+\frac{3^m-1}{2}-2^m+\frac{3^{m-1}+1}{2}-2^{m-1}\\
&=&2\cdot\sharp\mathcal{L}_{m-1}+\frac{3^m+1}{2}-2^m.
\end{eqnarray*}
Combining this with $\sharp\mathcal{L}_2=0$, we can easily get
$$\sharp\mathcal{L}_m=\frac{3^{m+1}-1}{2}-(m-2)\cdot2^m-26\cdot 2^{m-3}\quad\mbox{ for } m\geq 2. \Box$$

As for primitive graph eigenvalues $\mathcal{P}_m$, things become
more complicated. We consider $\mathcal{P}_m^+$ and
$\mathcal{P}_m^-$ respectively. We will show in the next section the
spectral decimation recipe for this type of eigenvalues can not be
used directly. In fact there is even not an analytic relation
between elements in $\mathcal{P}_m^+$ (or $\mathcal{P}_m^-$) and
elements in $\mathcal{P}_{m+1}^+$ (or $\mathcal{P}_{m+1}^-$). A
rough but intuitive explanation of why does this ``bad" thing happen
is that the Dirichlet boundary condition will be destroyed when we
use the eigenfunction extension algorithm $(\ref{3})$ to extend a
$\lambda_m$-eigenfunction $u_m$ from $\Omega_m$ to $\Omega_{m+1}$ or
restrict a $\lambda_{m+1}$-eigenfunction $u_{m+1}$ from
$\Omega_{m+1}$ to $\Omega_m$. However, a weak but useful relation
between $\mathcal{P}_m^+$ (or $\mathcal{P}_m^-$) and
$\mathcal{P}_{m+1}^+$ (or $\mathcal{P}_{m+1}^-$) will be found in
the next section, which will take the place of spectral decimation
in the further discussion. Let $\phi_{\pm}(x)$ be the same functions
as defined in $(\ref{4})$. We will prove that:

\textbf{Theorem 3.3.} \emph{For each $m\geq 2$, $\mathcal{P}_m^+$
consists of $r_m:=2^m+2^{m-2}-2$ distinct eigenvalues with
multiplicity 1, between $0$ and $6$ strictly, denoted by
$\lambda_{m,1},\lambda_{m,2},\cdots,\lambda_{m,r_m}$ in increasing
order, satisfying
$$0<\lambda_{m,1}<\lambda_{m,2}<\cdots<\lambda_{m,r_m-1}<5<\lambda_{m,r_m}<6.$$
 Moreover,
$r_{m+1}=2r_m+2$ and
\begin{eqnarray*}
0&<&\lambda_{m+1,1}<\phi_{-}(\lambda_{m,1}),\\
\phi_{-}({\lambda_{m,k-1}})&<&\lambda_{m+1,k}<\phi_{-}(\lambda_{m,k}),
\quad\forall 2\leq k\leq r_m,\\
\phi_{-}({\lambda_{m,r_m}})&<&\lambda_{m+1,r_m+1}<\phi_{+}(\lambda_{m,r_m}),\\
\phi_{+}({\lambda_{m,2r_{m}+2-k}})&<&\lambda_{m+1,k}<\phi_{+}(\lambda_{m,2r_{m}+1-k}),
\quad\forall r_m+2\leq k\leq 2r_m,\\
\phi_{+}(\lambda_{m,1})&<&\lambda_{m+1,2r_m+1}<5,\\
5&<&\lambda_{m+1,2r_m+2}<6.
\end{eqnarray*}
Similar properties hold for $\mathcal{P}_m^-$ with $r_m$ replaced by
$s_m:=2^m-2$.}

 In order to study the relation between
$\mathcal{P}_m^+$ (or $\mathcal{P}_m^-$) and $\mathcal{P}_{m+1}^+$
(or $\mathcal{P}_{m+1}^-$), we introduce the following notations. In
symmetric case, let $\widetilde{\phi}_-({\lambda_{m,1}})$ denote the
$(m+1)$-level eigenvalue between $0$ and $\phi_-(\lambda_{m,1})$.
Let $\widetilde{\phi}_{-}(\lambda_{m,k})$ denote the $(m+1)$-level
eigenvalue between $\phi_{-}(\lambda_{m,k-1})$ and
$\phi_{-}(\lambda_{m,k})$ for  $2\leq k\leq r_m$. Let
$\widetilde{\phi}_{+}(\lambda_{m,k})$ denote the $(m+1)$-level
eigenvalue between $\phi_{+}(\lambda_{m,k})$ and
$\phi_{+}(\lambda_{m,k-1})$ for $2\leq k\leq r_m$. Let
$\widetilde{\phi}_+(\lambda_{m,1})$ denote the $(m+1)$-level
eigenvalue between $\phi_+(\lambda_{m,1})$ and $5$. Call this kind
of $(m+1)$-level eigenvalues \emph{continued eigenvalues}. There
remain two other $(m+1)$-level eigenvalues: one is between
$\phi_-(\lambda_{m,r_m})$ and $\phi_{+}(\lambda_{m,r_m})$, the other
 is between $5$ and $6$. Call these two $(m+1)$-level eigenvalues
\emph{initial eigenvalues} with generation of birth $m+1$.  For the
$2$ level, all $r_2=3$ symmetric primitive  eigenvalues
$\lambda_{2,1}, \lambda_{2,2}$ and $\lambda_{2,3}$ are called
initial eigenvalues with generation of birth $2$. We define the
similar notations for skew-symmetric case in an obvious way with
$r_m$ replaced by $s_m$. From this point of view, the continued
primitive eigenvalues in $\mathcal{P}_{m+1}^+$ (or
$\mathcal{P}_{m+1}^-$) will be those arise from eigenvalues in
$\mathcal{P}_m^+$ (or $\mathcal{P}_m^-$) by a
$\widetilde{\phi}_{\pm}$ bifurcation similar (but never equal) to
 $\phi_{\pm}$ bifurcation. We call this phenomenon \emph{weak
spectral decimation}, which will be proved playing a critical role
in the study of the exact structure of primitive eigenvalues on
$\Omega$ in stead of spectral decimation.

We should emphasize here that $\widetilde{\phi}_\pm$ is not a real
function relation. (It is just a notation for simplicity.) See the
following diagram for the relation between $\mathcal{P}_m^+$  and
$\mathcal{P}_{m+1}^+$. The skew-symmetric case is similar.

\begin{flushleft}
\scriptsize{
\begin{tabular}{lllllllllllll}
     &  &$\quad\quad\quad\quad\lambda_{m,1}$ & &$\quad\quad\quad\quad\lambda_{m,2}$ && $\cdots$ & $\quad\quad\quad\quad\lambda_{m,r_m}$\\
    &&  $\quad\swarrow$ & $\searrow$ & $\quad\swarrow$ & $\searrow$ & $\cdots$ & $\quad\swarrow$ & $\searrow$\\
    &  & $\widetilde{\phi}_-(\lambda_{m,1})$ & $\widetilde{\phi}_+(\lambda_{m,1})$ & $\widetilde{\phi}_-(\lambda_{m,2})$ & $\widetilde{\phi}_+(\lambda_{m,2})$ & $\cdots$ & $\widetilde{\phi}_-(\lambda_{m,2})$ & $\widetilde{\phi}_+(\lambda_{m,2})$ & $\lambda_{m+1,r_m+1}$ & $\lambda_{m+1,r_{m+1}}$ \\
\end{tabular}}
\end{flushleft}

The structure of $\mathcal{M}_m$ depends on the structure of all
$\mathcal{P}_k^-$'s with $k<m$ by the definition of $\mathcal{M}_m$.
In fact, it is easy to check that

\textbf{Theorem 3.4.} \emph{Let $m\geq 2$, then for each eigenvalue
$\lambda_m$ in $\mathcal{M}_m$, it has multiplicity $2^{m-m'}$,
where $m'$ is the type of $\lambda_m$. Moreover, $\sharp
\mathcal{M}_m=(m-3)\cdot 2^m+4.$}

\emph{Proof.} Let $\lambda_{m'}$ be the graph eigenvalue associated
to $u_{m'}$ as used in Definition 3.10. Then obviously,
$\lambda_{m}=5^{m-m'}\lambda_{m'}$ and $\lambda_{m}$ has
multiplicity $2^{m-m'}$. Hence by using Theorem 3.3,
$$\sharp \mathcal{M}_m=\sum_{m'=2}^{m-1}2^{m-m'}\sharp\mathcal{P}_{m'}^-=\sum_{m'=2}^{m-1}2^{m-m'}(2^{m'}-2)=(m-3)\cdot 2^m+4.\quad\Box$$

From Theorem 3.2, Theorem 3.3 and Theorem 3.4, it is easy to check
that $\sharp \mathcal{L}_m$, $\sharp \mathcal{P}_m$ and $\sharp
\mathcal{M}_m$ add up to $\sharp(V_m^\Omega\setminus\partial
\Omega_m)$ since
\begin{equation}\label{11}
\sharp\mathcal{L}_m+\sharp \mathcal{P}_m+\sharp
\mathcal{M}_m=\frac{3^{m+1}-1}{2}-(m-2)\cdot2^m-26\cdot
2^{m-3}+r_m+s_m+(m-3)\cdot 2^m+4=a_m.
\end{equation}
This means there is no more Dirichlet eigenvalues except in the case
of localized, primitive and miniaturized types. Or, more precisely,
we will prove in the next section that

\textbf{Theorem 3.5.}\emph{ Let $m\geq 2$, then all the above
mentioned three types of Dirichlet eigenfunctions of $-\Delta_m$ on
$\Omega_m$ are linearly independent. Moreover, the Dirichlet
spectrum $\mathcal{S}_m$ of $-\Delta_m$ on $\Omega_m$ satisfies
$$\mathcal{S}_m=\mathcal{L}_m\cup \mathcal{P}_m^{+}\cup \mathcal{P}_m^{-}\cup
\mathcal{M}_m,$$ where the union is disjoint.}

 Hence we have the complete Dirichlet
spectrum $\mathcal{S}_m$ of $-\Delta_m$ on $\Omega_m$.  In Table
3.1, we list the eigenspace dimensions of all different types of
eigenvalues in $\mathcal{S}_m$ for level $m=2,3,4,5$.
\begin{flushleft}
\begin{center}
\scriptsize{
\begin{tabular}{cccccccccc}
  \hline
level & $\sharp \mathcal{L}_m$  & $\sharp \mathcal{P}^+_m$ & $\sharp \mathcal{P}^-_m$ & $\sharp \mathcal{M}_m$ & $\sharp \mathcal{S}_m$\\
\hline
  m &$\frac{3^{m+1}-1}{2}-(m-2)\cdot2^m-26\cdot 2^{m-3}$ & $2^m+2^{m-2}-2$& $2^m-2$ & $(m-3)\cdot 2^m+4$&$\frac{3^{m+1}-1}{2}-2^{m+1}$ \\
  2 & 0 & 3 & 2 & 0 & 5 \\
  3 & 6 & 8 & 6 & 4 & 24\\
  4 & 37 & 18 & 14 & 20 & 89\\
  5 & 164 & 38 & 30 & 68 & 300\\
 \hline
\end{tabular}}
\begin{center}
\textbf{Table 3.1.} \small{Eigenspace dimensions of different types
of eigenvalues in $\mathcal{S}_m$.}
\end{center}
\end{center}
\end{flushleft}

Next we want to pass the approximations to the limit.

For $\mathcal{L}$ case, assume that  $\{\lambda_m\}_{m\geq m_0}$ is
an infinite sequence of localized graph  eigenvalues related by
$\phi_{\pm}$ relations, with all but a finite number of
$\phi_{-}$'s. Then we define
$$\lambda=\frac{3}{2}\lim_{m\rightarrow\infty}5^m\lambda_m.$$ By
successively using the eigenfunction extension algorithm $(\ref{3})$
from a $\lambda_{m_0}$-eigenfunction $u_{m_0}$ of $-\Delta_{m_0}$ on
$\Omega_{m_0}$, one can extend $u_{m_0}$ to a localized
eigenfunction $u$ of $-\Delta_D$ on $\Omega$ associated to
$\lambda$. This method generates all the localized  eigenvalues
$\mathcal{L}$ as for the $\mathcal{SG}\setminus{V_0}$ case.

For $\mathcal{P}^+$  case, for an infinite sequence of
$\mathcal{P}^+$ type graph eigenvalues $\{\lambda_m\}_{m\geq m_0}$
related by $\widetilde{\phi}_{\pm}$ relations, with all but a finite
number of $\widetilde{\phi}_{-}$'s, we define
$$\lambda=\frac{3}{2}\lim_{m\rightarrow\infty}5^m\lambda_m.$$ We
will also show the existence of the limit $\lambda$. As pointed out
before, now we could not use the eigenfunction extension algorithm
directly. However, we will prove that $\lambda$ is still a
$\mathcal{P}^+$ type eigenvalue of $-\Delta_D$ on $\Omega$ by a
nonconstructive method. That is in Section 5, we will prove

\textbf{Theorem 3.6.} \emph{For each sequence of symmetric primitive
eigenvalues $\{\lambda_m\}_{m\geq m_0}$ related by
$\widetilde{\phi}_{\pm}$ relations, with all but a finite number of
$\widetilde{\phi}_{-}$'s, the limit
$\lambda=\frac{3}{2}\lim_{m\rightarrow\infty}5^m\lambda_m$ exists.
Moreover, $\lambda\in \mathcal{P}^+$.}

Furthermore, all $\mathcal{P}^+$ type eigenvalues come in this way.
 This will be done in Section 6.

\textbf{Theorem 3.7.} \emph{For each element $\lambda\in
\mathcal{P}^+$, there is uniquely a sequence of symmetric primitive
eigenvalues $\{\lambda_m\}_{m\geq m_0}$ as described in Theorem 3.6
such that
$\lambda=\frac{3}{2}\lim_{m\rightarrow\infty}5^m\lambda_m$.}

The $\mathcal{P}^-$ and $\mathcal{M}$ cases are completely similar
to the $\mathcal{P}^+$ case.

Now we could consider the complete  spectrum $\mathcal{S}$ of
$-\Delta_D$ on $\Omega$. That is, we will prove

\textbf{Theorem 3.8. } \emph{$\mathcal{S}=\mathcal{L}\cup
\mathcal{P}\cup \mathcal{M}$ where the union is disjoint.}

The following is a description of the multiplicity of each of the
above mentioned three types of eigenvalues.

\textbf{Theorem 3.9.} \emph{If $\lambda\in \mathcal{L}$ with
generation of birth of $m_0$, then  $\lambda$ is either a $5$-series
eigenvalue with multiplicity $\frac{3^{{m_0}-1}+1}{2}-2^{{m_0}-1}$
or a $6$-series eigenvalue with multiplicity
$\frac{3^{m_0}-1}{2}-2^{m_0}$; If $\lambda\in \mathcal{P}$, then the
multiplicity of $\lambda$ is $1$; If $\lambda\in \mathcal{M}$ is a
$k$-contracted miniaturized eigenvalue, then the multiplicity of
$\lambda$ is $2^{k}$.}

\emph{Proof.} It follows directly from Theorem 3.1, Theorem 3.2 and
Definition 3.6. $\Box$

Now we describe the Weyl's eigenvalue asymptotics on $\Omega$. As
introduced in the introduction section, let $\rho^{0}(x)$ and
$\rho^\Omega(x)$ be the Dirichlet eigenvalue counting functions for
$\mathcal{SG}\setminus{V_0}$ case and $\Omega$ case respectively.
Then we will prove in Section 6 the following comparison between
$\rho^{0}(x)$ and $\rho^\Omega(x)$.

\textbf{Theorem 3.10.} \emph{There exists a positive constant $C$ such
that for sufficiently large $x$, we have
\begin{equation}\label{sssss}
0\leq \rho^{0}(x)-\rho^\Omega(x)\leq C x^{\log2/\log 5}\log x.
\end{equation}}

In Section 8, we present the eigenvalues and their multiplicities in
$\mathcal{S}_m$ for level $m=2,3,4,5$(see Table 8.1, 8.2, 8.3 and
8.4).

We will prove Theorem 3.3 and Theorem 3.5 in Section 4, Theorem 3.6
in Section 5, Theorem 3.1, Theorem 3.7, Theorem 3.8 and Theorem 3.10
in Section 6.

\subsection{Neumann spectrum}

Before going to the Neumann spectrum of $-\Delta$ on $\Omega$, we
should give a precise description of the non-positive self-adjoint
operator $\Delta_N$, the Neumann Laplacian, on $\Omega$ under
consideration.

For $m\in \mathbb{N}$, define the renormalization energy on
$\Omega_m$ by
$$\mathcal{E}_m^\Omega(u,v):=r^{-m}\sum_{x,y\in\Omega_m,x\sim_my}(u(x)-u(y))(v(x)-v(y)),$$ for $u,v\in \mathbb{R}^{\Omega_m}$, where
$r=\frac{3}{5}$ is the energy renormalization factor. Then it is
easy to see that $\{\mathcal{E}_m^\Omega\}_{m\in \mathbb{N}}$ forms
a compatible sequence of discrete Dirichlet forms in the sense of
Kigami $\cite{Ki6}$(Definition 2.2.1), and hence we can define a
resistance form $(\mathcal{E}^\Omega, \mathcal{F}^\Omega)$ on
$\Omega_*:=\bigcup_{m\in \mathbb{N}}\Omega_m$ by
$$\mathcal{F}^\Omega:=\{u\in \mathbb{R}^{\Omega_*}|\lim_{m\rightarrow\infty}\mathcal{E}_m^\Omega(u|_{\Omega_m},v|_{\Omega_m})<\infty\},$$
$$\mathcal{E}^\Omega(u,v):=\lim_{m\rightarrow\infty}\mathcal{E}_m^\Omega(u|_{\Omega_m},v|_{\Omega_m}),\quad u,v\in\mathcal{F}^\Omega.$$
Let $R^\Omega$ be the resistance metric on $\Omega_*$ associated
with the resistance form $(\mathcal{E}^\Omega,\mathcal{F}^\Omega)$,
and let $\widetilde{\Omega}$ be the $R^\Omega$-completion of
$\Omega_*$. Then each $u\in \mathcal{F}^\Omega$ is naturally
identified with its unique extension to a continuous function on
$\widetilde{\Omega}$ by virtue of its $1/2$-H\"{o}lder continuity
with respect to $R^\Omega$. Kigami$\cite{Ki6}$(Definition 2.3.10)
assures that the resulting bilinear form
$(\widetilde{\mathcal{E}},\widetilde{\mathcal{F}})$, defined for
functions on $\widetilde{\Omega}$, is a resistance form on
$\widetilde{\Omega}$ with resistance metric (the completion of)
$R^\Omega$, which was actually done by Kigami and
Takahashi$\cite{Ki4}$. On the other hand, it is not difficult to
show that $\widetilde{\Omega}$ can be identified as
$(\mathcal{SG}\setminus L)\cup \widetilde{L}$, where $L$ denotes the
bottom line segment of $\mathcal{SG}$ and $\widetilde{L}$ is the
Cantor set naturally appearing as the ``boundary" at the bottom of
the graphs $\Omega_m$. Roughly speaking, $(\mathcal{SG}\setminus
L)\cup\widetilde{L}$ is obtained from $\mathcal{SG}$ by
distinguishing the left and right sides of each dyadic rational on
the bottom line $L$ to regard $L$ as a Cantor set $\widetilde{L}$.

From here on we will identify $\widetilde{\Omega}$ and
$(\mathcal{SG}\setminus L)\cup \widetilde{L}=\Omega\cup\{q_0\}\cup
\widetilde{L}$. Call $\{q_0\}\cup \widetilde{L}$ the boundary of
$\widetilde{\Omega}$, denoted by
$\widetilde{\Omega}\setminus\Omega$. Define a Borel measure
$\widetilde{\mu}$ on $\widetilde{\Omega}$ by
$\widetilde{\mu}(A)=\mu(A\cap\Omega)$ for each Borel subset $A$ of
$\widetilde{\Omega}$. Now we define the Neumann Laplacian as the
non-positive self-adjoint operator  associated with the Dirichlet
form $(\widetilde{\mathcal{E}},\widetilde{\mathcal{F}})$ on
$L^2(\widetilde{\Omega},\widetilde{\mu})$  in the following way.

\textbf{Definition 3.11.} \emph{ The Neumann Laplacian $\Delta_N$ on
$\Omega$ with domain $\mathcal{D}[\Delta_N]$ is formulated as
follows: for $u\in \widetilde{\mathcal{F}}$ and $f\in
L^2(\widetilde{\Omega},\widetilde{\mu})$,
$$u\in\mathcal{D}[\Delta_N]\mbox{ and } -\Delta_N u=f \mbox {\quad\quad if and only if\quad\quad } \widetilde{\mathcal{E}}(u,v)=\int_\Omega fvd\mu\mbox{ for any } v\in \widetilde{\mathcal{F}}.$$
}

\textbf{Definition 3.12.} \emph{For $\lambda\in \mathbb{R}$ and
$u\in\mathcal{D}[\Delta_N]$ if
$$-\Delta_N u=\lambda u,$$
then $\lambda$ is called an eigenvalue of $-\Delta_N$ on $\Omega$
(or, a Neumann eigenvalue of $-\Delta$ on $\Omega$), and $u$ is
called an associated (Neumann) eigenfunction.}

Let $\mathcal{S}^N$ denote the spectrum of $-\Delta_N$ on $\Omega$
($\mathcal{S}^N$ is also called the Neumann spectrum of $-\Delta$ on
$\Omega$).

Similarly to the Dirichlet case, the Neumann Laplacian $\Delta_N$ on
$\Omega$ could also be realized by the limit of graph Laplacians
$\Delta_m$ on $\Omega_m$. Hence it is natural to believe that the
discrete Neumann spectrum of $-\Delta_m$ on $\Omega_m$ should
converge to the spectrum of the Neumann Laplacian on $\Omega$. Thus
we need to analyze the discrete Neumann spectra first. We denote
$\mathcal{S}^{N}_m$ the Neumann spectrum of $-\Delta_m$ on
$\Omega_m$ for $m\in \mathbb{N}$.

To study the Neumann spectrum we  impose a Neumann condition on the
graph $\Omega_m$ by extending functions from $\Omega_m$ by even
reflection, and imposing the pointwise eigenvalue equation at the
boundary points in $\partial\Omega_{m}$, which now have $4$
neighbors. Then the Neumann $\lambda_m$-eigenvalue equations consist
of exactly $\sharp V_m^\Omega$ equations in $\sharp V_m^\Omega$
unknowns. It is even convenient to allow $m=1$, in which case there
are three equations associated to the boundary $\partial\Omega_1$
and no others. In particular, on $\Omega_1$ we find eigenvalues
$\lambda_1=0$ corresponding to the constant function, and
$\lambda_1=6$ corresponding to the two dimensional space of
functions satisfying $u(q_0)+u(F_1q_0)+u(F_2q_0)=0$ which can be
split into an one dimensional symmetric space and an one dimensional
skew-symmetric space under the reflection symmetry fixing $q_0$. For
simplicity, let $b_m=\sharp V_m^\Omega$. It is easy to check that
$b_1=3$, $b_2=10$, and more generally,

\textbf{Proposition. 3.2.} $b_m=\frac{3^{m+1}+1}{2}-2^m, \forall
m\in \mathbb{N}.$

\emph{Proof.}
$b_m=a_m+\sharp\partial\Omega_m=\frac{3^{m+1}-1}{2}-2^{m+1}+2^m+1=\frac{3^{m+1}+1}{2}-2^m.$
$\Box$

Similarly to the Dirichlet case,  $\mathcal{S}^{N}$ will also
consist of three types of  eigenvalues, localized, primitive and
miniaturized, with obvious modifications, denoted by
$\mathcal{L}^N$, $\mathcal{P}^N$ and $\mathcal{M}^N$ respectively.
And correspondingly, $\mathcal{S}^{N}_m$ will consist of three types
of graph Neumann eigenvalues, denoted by $\mathcal{L}^N_m$,
$\mathcal{P}^N_m$ and $\mathcal{M}^N_m$ respectively. Moreover,
$\mathcal{P}^N$($\mathcal{P}^N_m$) can also be split into symmetric
part $\mathcal{P}^{+,N}$($\mathcal{P}_m^{+,N}$) and skew-symmetric
part $\mathcal{P}^{-,N}$($\mathcal{P}^{-,N}_m$).

The structure of localized (graph) Neumann eigenvalues is very
similar to the Dirichlet case, with only a few changes:

\textbf{Theorem 3.11.} \emph{Let $m\geq 1$, then
$\sharp\mathcal{L}_m^N=\frac{3^{m+1}-1}{2}-2^{m+1}-(m-1)\cdot 2^m$.
Moreover, the initial eigenvalues in $\mathcal{L}_m^N$ are $5$ and
$6$ with multiplicity
$\rho_m^{\Omega,N}(5):=\frac{3^{m-1}+1}{2}-2^{m-1}$ and
$\rho_m^{\Omega,N}(6):=\frac{3^m+1}{2}-2^m$ respectively.}

\emph{Proof.} Comparing to the Dirichlet case(see Theorem 3.2), the
$6$-series has multiplicity increasing by $1$, namely the
eigenfunction associated to $q_0$, while the $5$-series is
unchanged. Hence
$\rho_m^{\Omega,N}(6)=\rho_m^{\Omega}(6)+1=\frac{3^m+1}{2}-2^m$  and
$\rho_m^{\Omega,N}(5)=\rho_m^\Omega(5)=\frac{3^{m-1}+1}{2}-2^{m-1}$,
$\forall m\geq 1$. A similar discussion shows that
$$\sharp\mathcal{L}_m^N=2\cdot\sharp\mathcal{L}_{m-1}^N-\rho_{m-1}^{\Omega,N}(6)+\rho_{m}^{\Omega,N}(6)+\rho_m^{\Omega,N}(5).$$
Hence we have
$$\sharp\mathcal{L}_m^N=2\cdot\sharp\mathcal{L}_{m-1}^N+\frac{3^m+1}{2}-2^m,$$
which yields that
$$\sharp\mathcal{L}_m^N=\frac{3^{m+1}-1}{2}-2^{m+1}-(m-1)\cdot 2^m\quad\mbox{ for } m\geq
1,$$ since $\sharp\mathcal{L}_1^N=0$. $\Box$

The structure of primitive (graph) Neumann eigenvalues
$\mathcal{P}^N$($\mathcal{P}_m^N$) is also similar to the Dirichlet
case. We consider the symmetric and skew-symmetric case
respectively. We will prove that:

\textbf{Theorem 3.12.} \emph{For each $m\geq 1$,
$\mathcal{P}_m^{+,N}$ consists of $2^m$ distinct eigenvalues with
multiplicity 1, between $0$ and $6$ with $0,6$ included, denoted by
$\lambda_{m,1},\lambda_{m,2},\cdots,\lambda_{m,r_m}$ in increasing
order, satisfying
$$\lambda_{m,1}=0<\lambda_{m,2}<\cdots<\lambda_{m,2^m-1}<5<\lambda_{m,2^m}=6.$$
Moreover,
\begin{eqnarray*}
\phi_{-}({\lambda_{m,k-1}})&<&\lambda_{m+1,k}<\phi_{-}(\lambda_{m,k}),
\quad\forall 2\leq k\leq 2^m,\\
\phi_{+}({\lambda_{m,2^{m+1}-k+1}})&<&\lambda_{m+1,k}<\phi_{+}(\lambda_{m,2^{m+1}-k}),
\quad\forall 2^m+1\leq k\leq 2^{m+1}-1.
\end{eqnarray*}
Similar properties hold for $\mathcal{P}_m^{-,N}$ with $2^m$
replaced by $2^m-1$, and $\lambda_{m,1}>0$ in that case.}

For symmetric case, there is  a weak spectral decimation which
relates $\mathcal{P}_m^{+,N}$ and $\mathcal{P}_{m+1}^{+,N}$ by
introducing the following notations. Let
$\widetilde{\phi}_{-}(\lambda_{m,k})$ denote the $(m+1)$-level
eigenvalue between $\phi_{-}(\lambda_{m,k-1})$ and
$\phi_{-}(\lambda_{m,k})$ for $2\leq k\leq 2^m$. Let
$\widetilde{\phi}_{+}(\lambda_{m,k})$ denote the $(m+1)$-level
eigenvalue between $\phi_{+}(\lambda_{m,k})$ and
$\phi_{+}(\lambda_{m,k-1})$ for  $2\leq k\leq 2^m$.  Call this kind
of $(m+1)$-level eigenvalues \emph{continued eigenvalues}. Hence the
continued primitive eigenvalues in $\mathcal{P}_{m+1}^{+,N}$ will be
those arise from eigenvalues of $\mathcal{P}_m^{+,N}\setminus\{0\}$
by a $\widetilde{\phi}_{\pm}$ bifurcation similar (but never equal)
to $\phi_{\pm}$ bifurcation. There remain two other $(m+1)$-level
eigenvalues: one is $0$, called \emph{zero eigenvalue}, the other is
$6$, called \emph{initial eigenvalue} with generation of birth
$m+1$. See the following diagram of eigenvalues in
$\mathcal{P}_{m}^{+,N}$.
\begin{flushleft}
\scriptsize{
\begin{tabular}{lllllllllllll}
 &$\mathcal{P}_{1}^{+,N}:$ & $0$ && $\quad\quad\quad 6$\\
 &&$\downarrow$ && $\swarrow$ & $\searrow$\\
 &$\mathcal{P}_{2}^{+,N}:$ & $0$ & $\quad\quad\quad\quad\widetilde{\phi}_-(6)$& & $\quad\quad\quad\quad\widetilde{\phi}_+(6)$& & $\quad\quad\quad 6$ &\\
 & &$\downarrow$ & $\quad\quad\swarrow$ & $\searrow$& $\quad\quad\swarrow$ & $\searrow$& $\quad\swarrow$ & $\searrow$\\
 &$\mathcal{P}_{3}^{+,N}$: & $0$ & $\widetilde{\phi}_-\widetilde{\phi}_-(6)$ & $\widetilde{\phi}_+\widetilde{\phi}_-(6)$ & $\widetilde{\phi}_-\widetilde{\phi}_+(6)$ & $\widetilde{\phi}_+\widetilde{\phi}_+(6)$ & $\widetilde{\phi}_-(6)$ & $\widetilde{\phi}_+(6)$ & $6$\\
&$\quad\vdots$ & $\vdots$ &&& $\vdots$\\
\end{tabular}}
\end{flushleft}

  For skew-symmetric
 case, the only difference is that there is no zero eigenvalue in $\mathcal{P}_{m}^{-,N}$. $\mathcal{P}_m^{-,N}$ consists of
 $2^m-1$ distinct eigenvalues between $0$ and $6$, including
 $6$, where $6$ is an initial eigenvalue with generation of birth
 $m$ and the others are continued eigenvalues arise from previous
 level eigenvalues by a similar weak bifurcation.
  We have now the following
 decimation diagram of eigenvalues in $\mathcal{P}_m^{-,N}$.
\begin{flushleft}
\scriptsize{
\begin{tabular}{lllllllllllll}
 &$\mathcal{P}_{1}^{-,N}:$  && $\quad\quad\quad 6$\\
 & && $\swarrow$ & $\searrow$\\
 &$\mathcal{P}_{2}^{-,N}:$  & $\quad\quad\quad\quad\widetilde{\phi}_-(6)$& & $\quad\quad\quad\quad\widetilde{\phi}_+(6)$& & $\quad\quad\quad 6$ &\\
 &  & $\quad\quad\swarrow$ & $\searrow$& $\quad\quad\swarrow$ & $\searrow$& $\quad\swarrow$ & $\searrow$\\
 &$\mathcal{P}_{3}^{-,N}$:  & $\widetilde{\phi}_-\widetilde{\phi}_-(6)$ & $\widetilde{\phi}_+\widetilde{\phi}_-(6)$ & $\widetilde{\phi}_-\widetilde{\phi}_+(6)$ & $\widetilde{\phi}_+\widetilde{\phi}_+(6)$ & $\widetilde{\phi}_-(6)$ & $\widetilde{\phi}_+(6)$ & $6$\\
&$\quad\vdots$ & $\vdots$ &&& $\vdots$\\
\end{tabular}}
\end{flushleft}

As for the miniaturized Neumann eigenvalues, the structure of
$\mathcal{M}_m^{N}$ depends on the structure of all
$\mathcal{P}_k^{-,N}$'s with $k<m$ in a completely same way as the
Dirichlet case. In fact, it is easy to check that

\textbf{Theorem 3.13.} \emph{Let $m\geq 1$, then for each eigenvalue
$\lambda_m$ in $\mathcal{M}_m^N$, it has multiplicity $2^{m-m'}$,
where $m'$ is the type of $\lambda_m$. Moreover, $\sharp
\mathcal{M}_m^N=(m-2)\cdot 2^m+2.$}

\emph{Proof.} Let $\lambda_{m'}$ be the graph eigenvalue associated
to $u_{m'}$, as defined in the Neumann version Definition 3.10. Then
obviously, $\lambda_{m}=5^{m-m'}\lambda_{m'}$ and $\lambda_{m}$ has
multiplicity $2^{m-m'}$. Hence by using the fact that
$\sharp\mathcal{P}_k^{-,N}=2^k-1$ for all $k\geq 1$, $$\sharp
\mathcal{M}_m^N=\sum_{k=1}^{m-1}2^{m-k}\sharp\mathcal{P}_{k}^{-,N}=\sum_{k=1}^{m-1}2^{m-k}(2^k-1)=(m-2)\cdot
2^m+2\quad\mbox{ for } m\geq 1.\quad\Box$$

It is easy to check $\sharp \mathcal{L}_m^N$, $\sharp
\mathcal{P}_m^N$ and $\sharp \mathcal{M}_m^N$ add up to $\sharp
V_m^\Omega$, since
\begin{equation}\label{123}
\sharp\mathcal{L}_m^N+\sharp \mathcal{P}_m^N+\sharp
\mathcal{M}_m^N=\frac{3^{m+1}-1}{2}-2^{m+1}-(m-1)\cdot
2^m+2^m+2^m-1+(m-2)\cdot 2^m+2=b_m.
\end{equation}

A similar argument like the Dirichlet case will yields:

\textbf{Theorem 3.14.}\emph{ Let $m\geq 1$, then all the above
mentioned three types of Neumann eigenfunctions of $-\Delta_m$ on
$\Omega_m$ are linearly independent. Moreover, the spectrum
$\mathcal{S}_m^N$ of $-\Delta_m$ on $\Omega_m$ satisfies
$$\mathcal{S}_m^N=\mathcal{L}_m^N\cup \mathcal{P}_m^{+,N}\cup \mathcal{P}_m^{-,N}\cup
\mathcal{M}_m^N,$$ where the union is disjoint.}

Hence we have the complete Neumann spectrum of $-\Delta_m$.

In Table 3.2, we list the eigenspace dimensions of all different
types of eigenvalues in $\mathcal{S}^N_m$ for level $m=1,2,3,4,5$.
\begin{flushleft}
\begin{center}
\scriptsize{
\begin{tabular}{cccccccccc}
  \hline
level & $\sharp \mathcal{L}^N_m$  & $\sharp \mathcal{P}^{+,N}_m$ & $\sharp \mathcal{P}^{-,N}_m$ & $\sharp \mathcal{M}^N_m$ & $\sharp \mathcal{S}^N_m$\\
\hline
  m &$\frac{3^{m+1}-1}{2}-2^{m+1}-(m-1)\cdot
2^m$ & $2^m$& $2^m-1$ & $(m-2)\cdot 2^m+2$&$\frac{3^{m+1}+1}{2}-2^m$ \\
  1 & 0 & 2 & 1 & 0 & 3\\
  2 & 1 & 4 & 3 & 2 & 10 \\
  3 & 8 & 8 & 7 & 10 & 33\\
  4 & 41 & 16 & 15 & 34 & 106\\
  5 & 172 & 32 & 31 & 98 & 333\\
 \hline
\end{tabular}}
\begin{center}
\textbf{Table 3.2.} \small{Eigenspace dimensions of different types
of eigenvalues in $\mathcal{S}^N_m$.}
\end{center}
\end{center}
\end{flushleft}

Then a similar discussion on how to pass the approximations to the
limit leads to the  spectrum $\mathcal{S}^{N}$ of $-\Delta_N$ on
$\Omega$.

The counterpart of Theorem 3.6 becomes:

 \textbf{Theorem 3.15.}
\emph{For each sequence of symmetric primitive Neumann eigenvalues
$\{\lambda_m\}_{m\geq m_0}$ related by $\widetilde{\phi}_{\pm}$
relations, with all but a finite number of $\widetilde{\phi}_{-}$'s,
the limit
$\lambda:=\frac{3}{2}\lim_{m\rightarrow\infty}5^m\lambda_m$ exists.
Moreover, $\lambda\in \mathcal{P}^{+,N}$.}

Similarly, we will have

\textbf{Theorem 3.16. } \emph{$\mathcal{S}^N=\mathcal{L}^N\cup
\mathcal{P}^N\cup \mathcal{M}^N$ where the union is disjoint.}

The following is a description of the multiplicity of each of the
above mentioned three types of eigenvalues.

\textbf{Theorem 3.17.} \emph{If $\lambda\in \mathcal{L}^N$ with
generation of birth of $m_0$, then  $\lambda$ is either a $5$-series
eigenvalue with multiplicity $\frac{3^{{m_0}-1}+1}{2}-2^{{m_0}-1}$
or a $6$-series eigenvalue with multiplicity
$\frac{3^{m_0}+1}{2}-2^{m_0}$; If $\lambda\in \mathcal{P}^N$, then
the multiplicity of $\lambda$ is $1$; If $\lambda\in \mathcal{M}^N$
is a $k$-contracted miniaturized eigenvalue, then the multiplicity
of $\lambda$ is $2^{k}$.}

\emph{Proof.} It follows directly from Theorem 3.11, Theorem 3.12
and the Neumann version of Definition 3.6. $\Box$

As for the Weyl's eigenvalue asymptotics on $\Omega$. A same
argument as Theorem 3.10 works also for the Neumann Laplacian.

We will prove Theorem 3.12 and Theorem 3.14 in Section 7. We will
also make some comments in Section 7 on the proof of Theorem 3.15,
since comparing to its Dirichlet counterpart Theorem 3.6, there is
no direct analogue of Green's function for the Neumann Laplacian,
which is used in the proof of Theorem 3.6. Other proofs are omitted
since they can be easily modified suitably from the Dirichlet case.

\section{Primitive graph Dirichlet eigenvalues of $-\Delta_m$}
\begin{figure}[ht]
\begin{center}
\includegraphics[width=13.3cm,totalheight=6.2cm]{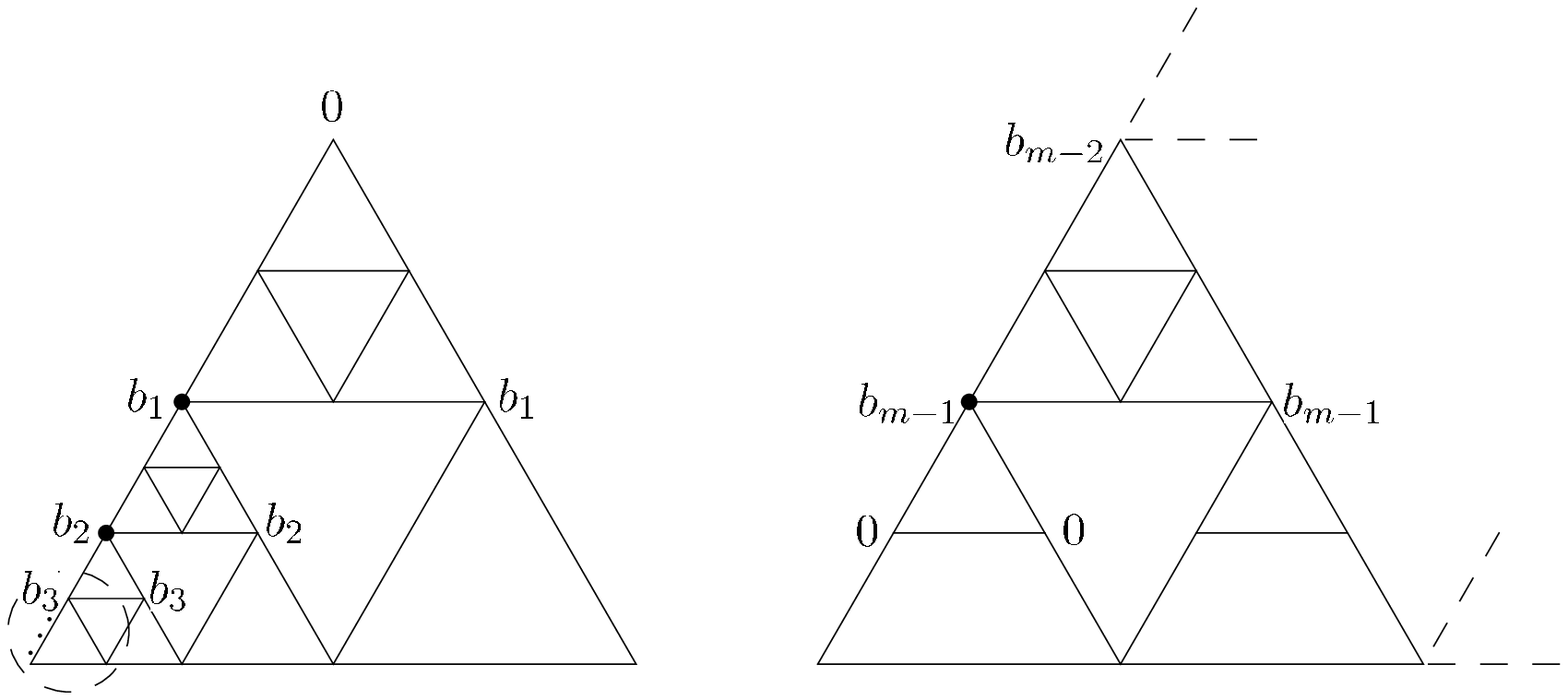}
\begin{center}
\textbf{Fig. 4.1.} \small{The values of the
$\lambda_m$-eigenfunction $u_m$ on the skeleton of $\Omega_m$ with
$\lambda_m\in\mathcal{P}_m^+$.}
\end{center}
\end{center}
\end{figure}
In this section, we work with $m$-level graph approximation
$\Omega_m$, $m=2,3,4\cdots$.  As introduced in Section 3, we use
$\mathcal{P}_m$ to denote the totality of the primitive graph
Dirichlet eigenvalues of the discrete Laplacian $-\Delta_m$ on
$\Omega_m$. Throughout this section, for simplicity, we omit the
terms ``graph" and ``Dirichlet" without causing any confusion. The
main object  in this section is to prove Theorem 3.3 and Theorem
3.5.

Let $f$ and $\phi_{\pm}(x)$ be the same functions as defined in
$(\ref{40})$ and $(\ref{4})$ respectively. In the following we use
$f^{(n)}$ to denote the $n$'th iteration of $f$, $n\geq 1$. In
particular, $f^{(0)}$ is the identity map of $\mathbb{R}$. If
$w=f^{(n)}(x)$, $w$ is called a \emph{successor} of $x$ of order $n$
with respect to $f$, and $x$ is called a \emph{predecessor} of $w$
of order $n$ with respect to $f$.

We begin with $\mathcal{P}_m^+$, the symmetric eigenvalues in
$\mathcal{P}_m$. Let $u_m$ be a $\lambda_m$-eigenfunction of
$-\Delta_m$ with $\lambda_m\in\mathcal{P}_m^+$. Denote by
$(b_0,b_1,b_2,\cdots,b_m)$ the values of $u_m$ on the skeleton
$(q_0,F_1q_0,F_1^2q_0,\cdots,F_1^mq_0)$ of $\Omega_m$ where
$b_0=b_m=0$ by the Dirichlet boundary condition. See Fig. 4.1.
 Write
$\lambda_{i}^{(m)}$ the successor of $\lambda_m$ of order $(m-i)$
with $2\leq i\leq m$ for simplicity. Assume that none of
$\lambda_{i}^{(m)}$'s is equal to $2$ or $5$. (Later we will show
this assumption automatically holds for any
$\lambda_m\in\mathcal{P}_m^+$.) The eigenfunction extension
algorithm $(\ref{3})$ gives the value of $u_m$ on the four
$(i+1)$-level neighbors of $F_1^iq_0$ for each $1\leq i\leq m-1$,
shown in Fig. 4.2.
\begin{figure}[ht]
\begin{center}
\includegraphics[width=6cm,totalheight=5.6cm]{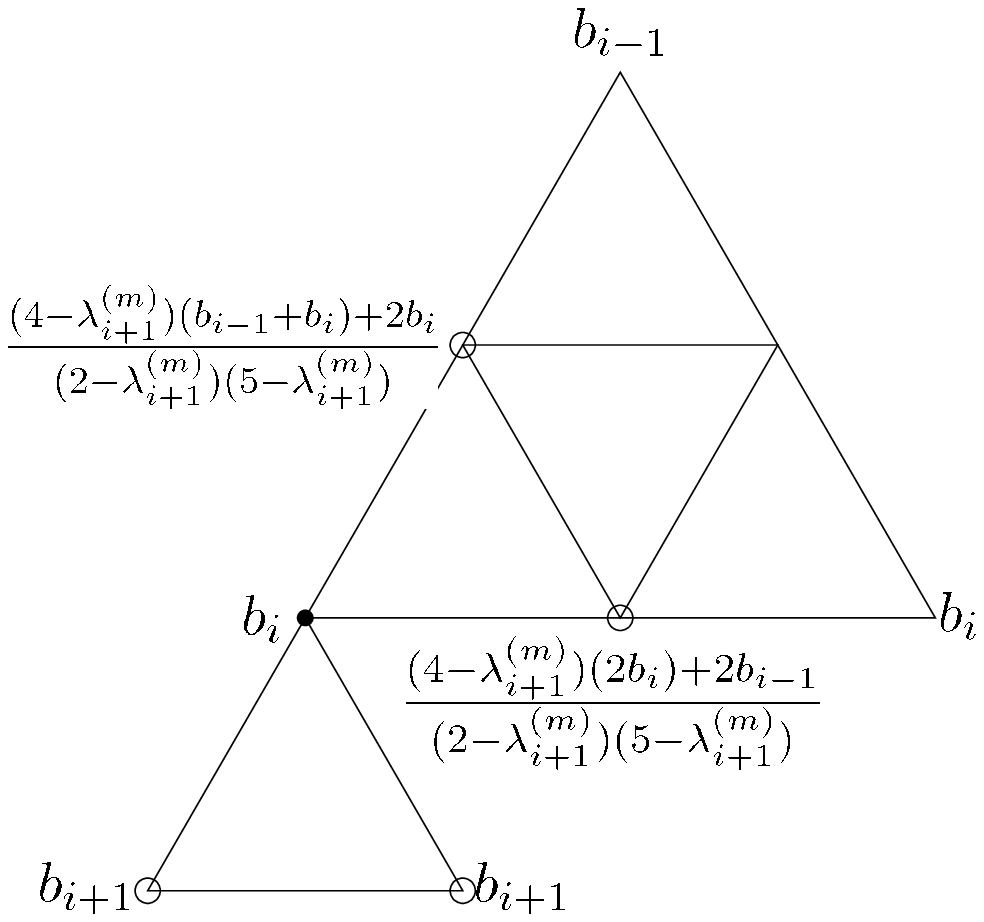}
\begin{center}
\textbf{Fig. 4.2.} \small{Values of $u_m$ on neighbors of
$F_1^iq_0$.}
\end{center}
\end{center}
\end{figure}
 Hence the $\lambda_{i+1}^{(m)}$-eigenvalue equation at the vertex $F_1^iq_0$
gives
\begin{equation}\label{91}
(4-\lambda_{i+1}^{(m)})b_i=2b_{i+1}+\frac{(14-3\lambda_{i+1}^{(m)})b_i+(6-\lambda_{i+1}^{(m)})b_{i-1}}{(2-\lambda_{i+1}^{(m)})(5-\lambda_{i+1}^{(m)})},\quad\forall
1\leq i\leq m-1,
\end{equation}
 which can be modified into
\begin{equation}\label{9}
l(\lambda_{i+1}^{(m)})b_{i-1}+s(\lambda_{i+1}^{(m)})b_i+r(\lambda_{i+1}^{(m)})b_{i+1}=0,\quad
\forall 1\leq i\leq m-1,
\end{equation}
 with $l(x):=x-6$,
$s(x):=(2-x)(4-x)(5-x)-(14-3x)$ and $r(x):=-2(2-x)(5-x)$. Still from
the eigenfunction extension algorithm, $u_m$ is uniquely determined
by $(b_1,b_2,\cdots,b_{m-1})$. Here $(b_1,b_2,\cdots,b_{m-1})$ can
be viewed as a non-zero vector solution of either of the above two
systems of equations consisting of $m-1$ equations in $m-1$
unknowns. Hence the determinants of them should both be equal to
$0$. For simplicity, we are interested in the second determinant,
although comparing to the first one, it brings the possibility that
$\lambda_{i}^{(m)}$ ($2\leq i\leq m$) could be $2$ or $5$, which
should be removed.

The determinant associated to system $(\ref{9})$ is a tridiagonal
determinant,
$$\left|
  \begin{array}{ccccc}
    s(\lambda_2^{(m)}) & r(\lambda_2^{(m)})   \\
    l(\lambda_3^{(m)}) & s(\lambda_3^{(m)}) & r(\lambda_3^{(m)})   \\
     & \ddots & \ddots & \ddots &  \\
     &  & l(\lambda_{m-1}^{(m)}) & s(\lambda_{m-1}^{(m)}) & r(\lambda_{m-1}^{(m)}) \\
     &  &  & l(\lambda_m^{(m)}) & s(\lambda_{m}^{(m)}) \\
  \end{array}
\right|.$$ Hence $\lambda_{m}$ should be a solution of the following
equation
\begin{equation}\label{6}
q_m(x):=\left|
  \begin{array}{ccccc}
    s(f^{(m-2)}(x)) & r(f^{(m-2)}(x))  \\
    l(f^{(m-3)}(x)) & s(f^{(m-3)}(x)) & r(f^{(m-3)}(x))   \\
     & \ddots & \ddots & \ddots &  \\
     &  & l(f(x)) & s(f(x)) & r(f(x)) \\
     &  &  & l(x) & s(x) \\
  \end{array}
\right|=0.
\end{equation}

Conversely, if $\lambda_m$ is a root of the polynomial $q_m(x)$ and
none of $f^{(i)}(\lambda_m)$'s with $0\leq i\leq m-2$ is equal to
$2$ or $5$, then $\lambda_m\in\mathcal{P}_m^+$. Hence we are
particular interested in all the root $x$'s of the polynomial
$q_m(x)$ excluding those satisfying $f^{(i)}(x)=2$ or $5$ for some
$0\leq i\leq m-2$.

We list some useful facts about the polynomial $q_m(x)$.

\textbf{Proposition 4.1.} \emph{Let $m\geq 2$, then}

\emph{(1) $q_m(0)>0$;}

\emph{(2) $q_m(5)>0$;}

\emph{(3) $q_m(6)<0$;}

\emph{(4) $q_m(\phi_{-}^{(m-1)}(5))<0$;}

\emph{(5)  $q_m(\phi_{-}^{(m-1)}(2))>0$;}

\emph{(6) $q_{m+2}(\phi_{-}^{(m-1)}(3))<0$ and $q_3(3)>0$.}

\emph{Proof.} (1) We will prove a stronger result,
\begin{equation}\label{7}
q_{m+1}(0)>20q_m(0)>0
\end{equation}
 for $m\geq 2$. This can
be proved by induction. It is easy to check that $q_2(0)=26>0$ and
$q_3(0)=556>20 q_2(0)$ by a direct computation.   If we assume
$q_m(0)>20q_{m-1}(0)>0$, then the expansion along the first row of
$q_{m+1}(0)$ yields that
$$q_{m+1}(0)=26q_{m}(0)-6\cdot 20q_{m-1}(0)>26q_{m}(0)-6q_{m}(0)=20q_{m}(0)>0.$$

(2) It is easy to compute that $q_2(5)=1>0$ and $q_3(5)=6>0$. For
$m\geq 4$, $q_m(5)=q_{m-1}(0)-20q_{m-2}(0)>0$ by using $(\ref{7})$.

(3) It is easy to compute that $q_2(6)=-4<0$, $q_3(6)=-3392\leq
q_2(6)<0$ and
$$q_{m}(6)=s(f^{(m-2)}(6))\cdot q_{m-1}(6)-r(f^{(m-2)}(6))\cdot l(f^{(m-3)}(6))\cdot q_{m-2}(6)$$
for $m\geq 4$ by the expansion along the first row of $q_m(6)$.

Consider a polynomial defined by $g_1(x):=s(f(x))-r(f(x)) l(x)$, it
is easy to check that $g_1(x)\geq 1$ whenever $x\leq -6$. In fact,
we can write $g_1(x)=(2-f(x))(5-f(x))(4-f(x)+2(x-6))-(14-3f(x))$ by
substituting the expressions for $s(f(x))$, $r(f(x))$ and $l(x)$.
Noticing that $4-f(x)+2(x-6)=x^2-3x-8\geq 46$ and $f(x)<0$ whenever
$x\leq -6$, we have $g_1(x)\geq
46(2-f(x))(5-f(x))-(14-3f(x))=46(f(x))^2-319f(x)+446$. Moreover,
since $f(x)\leq-66$ whenever $x\leq -6$, we finally have $g_1(x)\geq
46(-66)^2-319\cdot(-66)+446\geq 1.$

Then we can prove $q_m(6)\leq q_{m-1}(6)<0$ by induction. Suppose
$q_{m-1}(6)\leq q_{m-2}(6)<0$. (This is true for $m=4$.) Write
$q_m(6)=aq_{m-1}(6)+bq_{m-2}(6)$ with $a=s(f^{(m-2)}(6))$ and
$b=-r(f^{(m-2)}(6))\cdot l(f^{(m-3)}(6)).$ Noticing that $m\geq 4$,
we have $f^{(m-3)}(6)\leq -6$ and $f^{(m-2)}(6)<0$. Hence
$a+b=g_1(f^{(m-3)}(6))\geq 1$ and $b<0$. So by the induction
assumption, we have
$$q_{m}(6)\leq aq_{m-1}(6)+bq_{m-1}(6)=(a+b)q_{m-1}(6)\leq q_{m-1}(6)<0.$$

Hence we always have $q_m(6)<0$ for $m\geq 2$.

(4) For simplicity, denote $\alpha_m=q_m(\phi_-^{(m-1)}(5))$. By
direct computation, we have $\alpha_2=-4<0$ and $\alpha_3\approx
-92.10<0$. We will prove a stronger result, $\alpha_{m+1}\leq
10\alpha_m<0$, $\forall m\geq 2$. It holds for $m=2$. In order to
use the induction, we assume $\alpha_{m+1}\leq 10\alpha_m<0$. An
expansion of $\alpha_{m+2}$ along the last row yields that
$$\alpha_{m+2}=s(\phi_{-}^{(m+1)}(5))\alpha_{m+1}-r(\phi_{-}^{(m)}(5))l(\phi_{-}^{(m+1)}(5))\alpha_{m}.$$
 Since $2-\phi_{-}^{(m)}(5)>0$, $5-\phi_{-}^{(m)}(5)>0$ and
$\phi_{-}^{(m+1)}(5)-6<0$, we have
\begin{eqnarray*}
\alpha_{m+2}&=&s(\phi_{-}^{(m+1)}(5))\alpha_{m+1}-\frac{1}{10}r(\phi_{-}^{(m)}(5))l(\phi_{-}^{(m+1)}(5))\cdot(10\alpha_{m})\\
&\leq&
s(\phi_{-}^{(m+1)}(5))\alpha_{m+1}-\frac{1}{10}r(\phi_{-}^{(m)}(5))l(\phi_{-}^{(m+1)}(5))\alpha_{m+1}\\
&=&
[s(\phi_{-}^{(m+1)}(5))-\frac{1}{10}r(\phi_{-}^{(m)}(5))l(\phi_{-}^{(m+1)}(5))]\alpha_{m+1}.
\end{eqnarray*}

Consider a polynomial
$$g_2(x):=s(x)-\frac{1}{10}r(f(x))l(x)=14+9x-\frac{172}{5}x^2+\frac{87}{5}x^3-\frac{16}{5}x^4+\frac{1}{5}x^5.$$
It is easy to compute that
$$g_2'(x)=9-\frac{344}{5}x+\frac{261}{5}x^2-\frac{64}{5}x^3+x^4\geq 9-\frac{344}{5}(\phi_{-}^{(3)}(5))-\frac{64}{5}(\phi_{-}^{(3)}(5))^3\approx 4.91>0$$
whenever $0\leq x\leq \phi_-^{(3)}(5)$. Hence $g_2(x)$ is monotone
increasing in the interval $[0,\phi_{-}^{(3)}(5)]$.

Since $0<\phi_-^{(m+1)}(5)\leq\phi_-^{(3)}(5)$, we have
$g_2(\phi_{-}^{(m+1)}(5))\geq g_2(0)\geq 10$. Hence
$\alpha_{m+2}\leq g_2(\phi_{-}^{(m+1)}(5))\alpha_{m+1}\leq
10\alpha_{m+1}<0$.

The proofs of (5) and (6) are similar to that of (4). $\Box$

 Now we discuss the possibility  of the roots of $q_m(x)$
satisfying $f^{(i)}(x)=2$ or $5$ for some $0\leq i\leq m-2$. The
following well-known basic algebra lemma should be useful.

\textbf{Lemma 4.1.} \emph{Let $g, h$ be two polynomials whose
coefficients all belong to $\mathbb{Q}$ (the field of rational
numbers), i.e., $g,h\in \mathbb{Q}[x]$. If $g$ is irreducible in
$\mathbb{Q}[x]$ and $g, h$ have a common root in $\mathbb{R}$, then
$g$ divides $h$ in $\mathbb{Q}[x]$, i.e., all real roots of $g$
belong to those of $h$.}

 \textbf{Lemma 4.2.} \emph{Let $x$ be a predecessor of $2$ of order
$i$ with $0\leq i\leq m-3$. Then $q_m(x)=0$. Let $x$ be a
predecessor of $2$ of order $m-2$. Then $q_m(x)\neq 0$.}

\emph{Proof.} Firstly, let $m\geq 3$ and $x$ be a predecessor of $2$
of order $i$ with $0\leq i\leq m-3$. Then $f^{(i)}(x)=2$ and
$f^{(i+1)}(x)=6$. Substituting them into $(\ref{6})$, noticing
$s(2)=r(6)=-8$, $s(6)=l(2)=-4$ and $r(2)=l(6)=0$, we get
\begin{eqnarray*}
q_m(x)&=&\left|
  \begin{array}{ccccccc}
      \ddots & \ddots &\ddots  \\
   &l(f^{(i+1)}(x)) &s(f^{(i+1)}(x)) &r(f^{(i+1)}(x))  \\
         & &  l(f^{(i)}(x)) & s(f^{(i)}(x)) & r(f^{(i)}(x)) \\
     &&  & \ddots &\ddots & \ddots \\
  \end{array}
\right|\\&=&\left|
  \begin{array}{ccccccc}
      \ddots &\ddots &\ddots \\
   &0 & -4 &-8  &  \\
          & & -4 & -8 & 0 \\
   &  &  & \ddots &\ddots & \ddots \\
  \end{array}
\right|=0.
\end{eqnarray*}

Secondly, let $x$ be a predecessor of $2$ of order $m-2$. Then
$f^{(m-2)}(x)=2$. If $m=2$, then $x=2$. It is easy to check that
$x=2$ is not a root of $q_2(x)$. If $m\geq 3$, suppose $x$ is a root
of $q_m(x)$, then using Lemma 4.1, all roots of $f^{(m-2)}(x)-2$ are
roots of $q_m(x)$. Noticing that $\phi_-^{(m-2)}(2)$ is also a root
of $f^{(m-2)}(x)-2$, we have $q_m(\phi_-^{(m-2)}(2))=0$. But
\begin{eqnarray*}
q_m({\phi_{-}^{(m-2)}(2))}&=&\left|
  \begin{array}{ccccccc}
      s(2) &  r(2)  \\
      l(\phi_-(2)) & s(\phi_-(2)) & r(\phi_-(2))\\
     &\ddots &\ddots &\ddots\\
    && l(\phi_-^{(m-2)}(2))& s(\phi_-^{(m-2)}(2))\\
  \end{array}
\right|\\&=&s(2)\cdot\left|
  \begin{array}{ccccccc}
      s(\phi_-(2)) &  r(\phi_-(2)) \\
           \ddots &\ddots &\ddots\\
   & l(\phi_-^{(m-2)}(2))& s(\phi_-^{(m-2)}(2))\\
  \end{array}
\right|\\&=&(-8)\cdot q_{m-1}(\phi_{-}^{(m-2)}(2)),
\end{eqnarray*}
since $r(2)=0$. By using Proposition 4.1(5), we get
$q_m({\phi_{-}^{(m-2)}(2))}<0$ which contradicts to
$q_m({\phi_{-}^{(m-2)}(2))}=0$. Hence $q_m(x)\neq 0$. $\Box$

\textbf{Lemma 4.3.} \emph{Let $x$ be a predecessor of $5$ of order
$i$ with $0\leq i\leq m-2$. Then $q_m(x)\neq 0$.}

\emph{Proof.} Let $x$ be a predecessor of $5$ of order $i$ with
$0\leq i\leq m-2$. Then $f^{(i)}(x)=5$. Hence if $x$ is a root of
$q_m(x)$, then using Lemma 4.1, all roots of $f^{(i)}(x)-5$ are
roots of $q_m(x)$. Noticing that $\phi_{-}^{(i)}(5)$ is also a root
of $f^{(i)}(x)-5$, we have $q_m(\phi_{-}^{(i)}(5))=0$. But
$q_m(\phi_{-}^{(0)}(5))=q_m(5)>0$ by Proposition 4.1(2). More
generally for $0<i\leq m-2$,
\begin{eqnarray*}
q_m({\phi_{-}^{(i)}(5))}&=&\left|
  \begin{array}{ccccccc}
      s(f^{(m-2-i)}(5)) &  r(f^{(m-2-i)}(5))  \\
     \ddots &\ddots &\ddots\\
   &   l(5) & s(5)& r(5)&\\
   && l(\phi_-(5)) & s(\phi_-(5)) & r(\phi_{-}(5))\\
           & & &\ddots &\ddots &\ddots \\
   &  && & l(\phi_{-}^{(i)}(5)) & s(\phi_{-}^{(i)}(5)) \\
  \end{array}
\right|\\&=&\left|
  \begin{array}{ccccccc}
      s(f^{(m-2-i)}(5)) &  r(f^{(m-2-i)}(5)) \\
     \ddots &\ddots &\ddots\\
   &  l(5) & s(5)& 0&\\
   && l(\phi_-(5)) & s(\phi_-(5)) & r(\phi_{-}(5))\\
          &  & &\ddots &\ddots &\ddots \\
   &  &&  & l(\phi_{-}^{(i)}(5)) & s(\phi_{-}^{(i)}(5)) \\
  \end{array}
\right|,
\end{eqnarray*}
since $r(5)=0$. Thus
\begin{eqnarray*} q_m({\phi_{-}^{(i)}(5))}&=&\left|
  \begin{array}{ccccccccccc}
      s(f^{(m-2-i)}(5)) &  r(f^{(m-2-i)}(5))  \\
     \ddots &\ddots &\ddots\\
   &&   l(5) & s(5)
     \end{array}
\right|\\&&\cdot \left|\begin{array}{cccccccccc}
    s(\phi_-(5)) &  r(\phi_-(5))  \\
     \ddots &\ddots &\ddots\\
   &  & l(\phi_-^{(i)}(5))& s(\phi_-^{(i)}(5))
     \end{array}
\right|\\
&=&q_{m-i}(5)\cdot q_{i+1}(\phi_-^{(i)}(5))<0
\end{eqnarray*}
by the $2$'nd and $4$'th statements in Proposition 4.1. Hence
$\forall 0\leq i\leq m-2$, we have proved
$q_m({\phi_{-}^{(i)}(5))}\neq0$ which yields a contradiction to
$q_m({\phi_{-}^{(i)}(5))}=0$. So $q_m(x)\neq 0$. $\Box$

From Lemma 4.2 and Lemma 4.3, for $m\geq 3$, the total unwanted
roots of $q_m(x)$ are those predecessors of $2$ of order $i$ with
$0\leq i\leq m-3$. $q_2(x)$ does not have any unwanted root. Hence
to exclude them out, we define
$$p_m(x):=\frac{q_m(x)}{(x-2)(f(x)-2)\cdots (f^{(m-3)}(x)-2)}\quad \mbox{ for } m\geq
3,$$ and
$$p_2(x):=q_2(x)=s(x).$$ $p_m(x)$ is still a polynomial from Lemma 4.2, although it looks
like a rational function.

 Now we can say if $\lambda_m$ is a root
of the polynomial $p_m(x)$, then $\lambda_m\in\mathcal{P}_m^+$. Note
that the degree of the polynomial $q_m(x)$ is
$3+3\cdot2+\cdots+3\cdot 2^{m-2}=3(2^{m-1}-1)$ and  the number of
all the unwanted roots of $q_m(x)$ is $1+2+\cdots+2^{m-3}=2^{m-2}-1$
for $m\geq 3$ and $0$ for $m=2$. Hence it is easy to check that the
degree of $p_m(x)$ is $r_m:=2^m+2^{m-2}-2$.

The following is a list of some useful facts about the polynomial
$p_m(x)$.

\textbf{Proposition 4.2.} \emph{(1) $(-1)^mp_m(0)> 0$, $\forall
m\geq 2$;}

\emph{(2) $p_2(5)>0$ and  $(-1)^{m-1}p_m(5)>0$, $\forall m\geq 3$;}

\emph{(3) $p_2(6)<0$ and  $(-1)^{m}p_m(6)>0$, $\forall m\geq 3$.}

\emph{Proof.} It can be checked by a direct computation when $m=2$.
When $m\geq 3$, noticing that by the definition of $p_m(x)$,
$$p_{m}(0)=\frac{q_{m}(0)}{(-2)^{m-2}},\quad p_{m}(5)=\frac{q_{m}(5)}{3\cdot (-2)^{m-3}}$$ and $$ p_{m}(6)=\frac{q_{m}(6)}{(6-2)(f(6)-2)\cdots(f^{(m-3)}(6)-2)}.$$
 Using Proposition 4.1(1)-(3), we get
the desired result. $\Box$

We now present a more precise result about the distribution of the
roots of $p_m(x)$ and show an useful relation between roots of two
consecutive polynomials.

\textbf{Lemma 4.4.} \emph{Let $m\geq 2$. Then $p_m(x)$ has $r_m$
distinct real roots, denoted by $\lambda_{m,1},$ $\lambda_{m,2},$
$\cdots,$ $\lambda_{m,r_m}$ in increasing order, satisfying
$$0<\lambda_{m,1}<\lambda_{m,2}<\cdots<\lambda_{m,r_m-1}<5<\lambda_{m,r_m}<6.$$
Moreover, $(-1)^{m+k-1}p_{m+1}(\phi_{-}(\lambda_{m,k}))>0$ and
$(-1)^{m+k}p_{m+1}(\phi_{+}(\lambda_{m,k}))>0$, $\forall 1\leq k\leq
r_m$.}

\emph{Proof.} We prove it by using the induction on $m$.

When $m=2$, $p_2(x)=s(x)$ has $3$ distinct roots:
$\lambda_{2,1}\approx 1.0646$, $\lambda_{2,2}\approx 4.4626$ and
$\lambda_{2,3}\approx 5.4728$ by a direct computation.

Let $\lambda$ be one of $\lambda_{2,k}$'s, then $p_2(\lambda)=0$,
i.e., $s(\lambda)=0$, and
$p_3(\phi_-(\lambda))=\frac{q_3(\phi_-(\lambda))}{\phi_-(\lambda)-2}=\frac{2(\phi_-(\lambda)-6)(2-\lambda)(5-\lambda)}{\phi_-(\lambda)-2}$
by using $s(\lambda)=0$. Since $0<\lambda<6$, we have
$\phi_-(\lambda)-2<0$ and $\phi_-(\lambda)-6<0$. Hence
$p_3(\phi_-(\lambda))\sim (2-\lambda)(5-\lambda)$ where ``$\sim$"
means both sides of ``$\sim$" have the same signs. Similarly,
$p_3(\phi_+(\lambda))=\frac{2(\phi_+(\lambda)-6)(2-\lambda)(5-\lambda)}{\phi_+(\lambda)-2}$
and $p_3(\phi_+(\lambda))\sim -(2-\lambda)(5-\lambda)$.

Hence $0<\lambda_{2,1}<2$ yields that $p_3(\phi_-(\lambda_{2,1}))>0$
and $p_3(\phi_+(\lambda_{2,1}))<0$; $2<\lambda_{2,2}<5$ yields that
$p_3(\phi_-(\lambda_{2,2}))<0$ and $p_3(\phi_+(\lambda_{2,2}))>0$;
$\lambda_{2,1}>5$ yields that $p_3(\phi_-(\lambda_{2,3}))>0$ and
$p_3(\phi_+(\lambda_{2,3}))<0$. So our lemma holds for $m=2$.

We now assume our lemma holds for $m$, and prove it for $m+1$.

Noticing that from Proposition 4.2, we have $p_{m+1}(0)\sim
(-1)^{m-1}$, $p_{m+1}(5)\sim (-1)^{m}$ and $p_{m+1}(6)\sim
(-1)^{m-1}$. Hence if we write
\begin{equation}\label{8}
0,\phi_-(\lambda_{m,1}),\phi_-(\lambda_{m,2}),\cdots,\phi_{-}(\lambda_{m,r_m}),\phi_{+}(\lambda_{m,r_m}),\cdots,\phi_{+}(\lambda_{m,2}),\phi_{+}(\lambda_{m,1}),5,6
\end{equation}
in increasing order, then the values of $p_{m+1}$ on them have
alternating signs by the induction assumption. Hence there exist at
least $2r_m+2=r_{m+1}$ distinct roots of $p_{m+1}(x)$, with each
located strictly between each two consecutive  points in
$(\ref{8})$. Moreover, these are the totality of the roots of
$p_{m+1}(x)$ since the degree of $p_{m+1}(x)$ is also $r_{m+1}$.
Hence we can write them in increasing order:
$$0<\lambda_{m+1,1}<\lambda_{m+1,2}<\cdots<\lambda_{m+1,r_{m+1}-1}<5<\lambda_{m+1,r_{m+1}}<6.$$

Now we study the signs of $p_{m+2}(\phi_{\pm}(\lambda_{m+1,k}))$'s.
Let $\lambda$ be one of $\lambda_{m+1,k}$'s, then
$p_{m+1}(\lambda)=0$. Moreover,
\begin{eqnarray*}
p_{m+2}(\phi_-(\lambda))&=&\frac{q_{m+2}(\phi_-(\lambda))}{(\phi_-(\lambda)-2)(\lambda-2)\cdots(f^{(m-2)}(\lambda)-2)}\\
&=&\frac{s(\phi_-(\lambda))q_{m+1}(\lambda)+2(\phi_-(\lambda)-6)(2-\lambda)(5-\lambda)q_m(f(\lambda))}{(\phi_-(\lambda)-2)(\lambda-2)\cdots(f^{(m-2)}(\lambda)-2)}
\end{eqnarray*}
by using the expansion of $q_{m+2}(\phi_-(\lambda))$ along the last
row. Since $p_{m+1}(\lambda)=0$, we have $q_{m+1}(\lambda)=0$. Hence
$$p_{m+2}(\phi_-(\lambda))=\frac{2(\phi_-(\lambda)-6)(2-\lambda)(5-\lambda)q_m(f(\lambda))}{(\phi_-(\lambda)-2)(\lambda-2)\cdots(f^{(m-2)}(\lambda)-2)}=\frac{-2(\phi_-(\lambda)-6)(5-\lambda)p_m(f(\lambda))}{\phi_-(\lambda)-2}.$$
Since $0<\lambda<6$, we have $\phi_-(\lambda)-2<0$ and
$\phi_-(\lambda)-6<0$, hence
$$p_{m+2}(\phi_-(\lambda))\sim(\lambda-5)p_m(f(\lambda)).$$
Similarly,
$$p_{m+2}(\phi_+(\lambda))=\frac{-2(\phi_+(\lambda)-6)(5-\lambda)p_m(f(\lambda))}{\phi_+(\lambda)-2}$$
and
$$p_{m+2}(\phi_+(\lambda))\sim(5-\lambda)p_m(f(\lambda)).$$

When $\lambda=\lambda_{m+1,1}$, we have
$0<\lambda<\phi_-(\lambda_{m,1})$, hence
$0<f(\lambda)<\lambda_{m,1}$. Noticing that $\lambda_{m,1}$ is the
least root of $p_m(x)$ and $\lambda_{m,1}>0$ by the induction
assumption, we have $p_m(f(\lambda))\sim p_m(0)\sim (-1)^m$. Hence
$p_{m+2}(\phi_-(\lambda_{m+1,1}))\sim(-1)^{m+1}$ and
$p_{m+2}(\phi_+(\lambda_{m+1,1}))\sim(-1)^{m}$ since
$\lambda_{m+1,1}<5$.

When $\lambda=\lambda_{m+1,k}$ with $2\leq k\leq r_m$, we have
$\phi_-(\lambda_{m,k-1})<\lambda<\phi_-(\lambda_{m,k})$, hence
$\lambda_{m,k-1}<f(\lambda)<\lambda_{m,k}$. Noticing that
$p_m(\lambda_{m,k-1})=0$ and $p_{m}(0)\sim (-1)^m$, we have
$p_m(f(\lambda))\sim(-1)^{m+k+1}$ by using the induction assumption.
Hence $p_{m+2}(\phi_-(\lambda_{m+1,k}))\sim(-1)^{m+k}$ and
$p_{m+2}(\phi_+(\lambda_{m+1,k}))\sim(-1)^{m+k+1}$ since
$\lambda_{m+1,k}<5$.

When $\lambda=\lambda_{m+1,r_m+1}$, we have
$\phi_-(\lambda_{m,r_m})<\lambda<\phi_+(\lambda_{m,r_m})$, hence
$f(\lambda)>\lambda_{m,r_m}$. Noticing that $\lambda_{m,r_m}$ is the
last root of $p_m(x)$ and $p_{m}(0)\sim (-1)^m$, we have
$p_m(f(\lambda))\sim(-1)^{m+r_m}$ by using the induction assumption.
Hence $p_{m+2}(\phi_-(\lambda_{m+1,r_m+1}))\sim(-1)^{m+1+r_m}$ and
$p_{m+2}(\phi_+(\lambda_{m+1,r_m+1}))\sim(-1)^{m+r_m}$ since
$\lambda_{m+1,r_m+1}<5$.

When $\lambda=\lambda_{m+1,k}$ with $r_m+2\leq k\leq 2r_m$, we have
$\phi_+(\lambda_{m,r_{m+1}-k})<\lambda<\phi_+(\lambda_{m,r_{m+1}-k-1})$,
hence $\lambda_{m,r_{m+1}-k-1}<f(\lambda)<\lambda_{m,r_{m+1}-k}$.
Noticing that $p_m(\lambda_{m,r_{m+1}-k-1})=0$ and $p_{m}(0)\sim
(-1)^m$, we have
$p_m(f(\lambda))\sim(-1)^{m+r_{m+1}-k-1}\sim(-1)^{m+k-1}$ by using
the induction assumption. Hence
$p_{m+2}(\phi_-(\lambda_{m+1,k}))\sim(-1)^{m+k}$ and
$p_{m+2}(\phi_+(\lambda_{m+1,k}))\sim(-1)^{m+k-1}$ since
$\lambda_{m+1,k}<5$.

When $\lambda=\lambda_{m+1,2r_m+1}$, we have
$\phi_+(\lambda_{m,1})<\lambda<5$, hence $f(\lambda)<\lambda_{m,1}$.
So we have $p_m(f(\lambda))\sim p_m(0)\sim(-1)^{m}$. Hence
$p_{m+2}(\phi_-(\lambda_{m+1,2r_m+1}))\sim(-1)^{m+1}$ and
$p_{m+2}(\phi_+(\lambda_{m+1,2r_m+1}))\sim(-1)^{m}$ since
$\lambda_{m+1,2r_m+1}<5$.

When $\lambda=\lambda_{m+1,2r_m+2}$, we have $5<\lambda<6$, hence
$f(\lambda)<0$. So we have $p_m(f(\lambda))\sim p_m(0)\sim(-1)^{m}$.
But now $\lambda>5$, hence
$p_{m+2}(\phi_-(\lambda_{m+1,2r_m+2}))\sim(-1)^{m}$ and
$p_{m+2}(\phi_+(\lambda_{m+1,2r_m+2}))\sim(-1)^{m-1}$.

Hence we have proved
$(-1)^{m+1+k-1}p_{m+2}(\phi_{-}(\lambda_{m+1,k}))>0$ and
$(-1)^{m+1+k}p_{m+2}(\phi_{+}(\lambda_{m+1,k}))>0$, $\forall 1\leq
k\leq r_{m+1}$. So our lemma holds for $m+1$. $\Box$

Thus by Lemma 4.4, in particular the proof of Lemma 4.4 and the fact
that each root of $p_m(x)$ belongs to $\mathcal{P}_m^+$,  we have
the following result:

\textbf{Lemma 4.5.} \emph{For each $m\geq 2$, $\mathcal{P}_m^{+}$
consists of at least $r_m$ distinct eigenvalues satisfying
\begin{equation}0<\lambda_{m,1}<\lambda_{m,2}<\cdots<\lambda_{m,r_m-1}<5<\lambda_{m,r_m}<6.
\end{equation}
Moreover, \begin{eqnarray}\label{10}
0&<&\lambda_{m+1,1}<\phi_{-}(\lambda_{m,1}),\nonumber\\
\phi_{-}({\lambda_{m,k-1}})&<&\lambda_{m+1,k}<\phi_{-}(\lambda_{m,k}),
\quad\forall 2\leq k\leq r_m,\nonumber\\
\phi_{-}({\lambda_{m,r_m}})&<&\lambda_{m+1,r_m+1}<\phi_{+}(\lambda_{m,r_m}),\nonumber\\
\phi_{+}({\lambda_{m,2r_{m}+2-k}})&<&\lambda_{m+1,k}<\phi_{+}(\lambda_{m,2r_{m}+1-k}),
\quad\forall r_m+2\leq k\leq 2r_m,\\
\phi_{+}(\lambda_{m,1})&<&\lambda_{m+1,2r_m+1}<5,\nonumber\\
5&<&\lambda_{m+1,2r_m+2}<6.\nonumber
\end{eqnarray}}

\textbf{Remark.} \emph{The third inequality in $(\ref{10})$ can be
refined into $2<\lambda_{m+1,r_m+1}<\phi_{+}(\lambda_{m,r_m})$. See
details in Theorem A in Appendix.}

Moreover, we have

\textbf{Lemma 4.6.} \emph{Let $\lambda_m$ be a root of $p_m(x)$,
$u_m$  a primitive $\lambda_m$-eigenfunction on $\Omega_m$, and
$(b_0,b_1,\cdots,b_m)$  the values of $u_m$ on the skeleton of
$\Omega_m$. Then $b_1\neq 0$ and $b_{m-1}\neq 0$.}

\emph{Proof.} Without loss of generality, assume $m\geq 3$. We still
use $\lambda_{i}^{(m)}$ to denote the successor of $\lambda_m$ of
order $(m-i)$ with $2\leq i\leq m$. From the definition of $p_m(x)$,
$\lambda_{i}^{(m)}\neq 2$ or $5$, for each $2\leq i\leq m$. From the
discussion in the beginning of this section, the vector
$(b_1,b_2,\cdots,b_{m-1})$ can be viewed as a non-zero vector
solution of system $(\ref{9})$ of equations.

Suppose $b_{m-1}=0$. Then $(b_1,b_2,\cdots,b_{m-2})$ can be viewed
as a non-zero vector solution of the system of equations consisting
of the first $(m-2)$ equations of $(\ref{9})$ in $(m-2)$ unknowns.
Hence the determinant of this system $q_{m-1}(\lambda_{m-1}^{(m)})$
should be equal to $0$. Thus $\lambda_{m-1}^{(m)}$ is a root of
$p_{m-1}(x)$ since its all successors
$\lambda_{2}^{(m)},\cdots,\lambda_{m-1}^{(m)}$ do not take value
from $\{2,5\}$ obviously. Then Lemma 4.4 says that neither of
$\phi_{\pm}(\lambda_{m-1}^{(m)})$ should be a root of $p_m(x)$. This
contradicts to $p_m(\lambda_m)=0$ since $\lambda_m$ is equal to
either of $\phi_{\pm}(\lambda_{m-1}^{(m)})$. Hence $b_{m-1}\neq 0$.

On the other hand, if $b_1=0$, then by substituting it into
$(\ref{9})$, noticing that none of $\lambda_{i}^{(m)}$'s is equal to
$2$ or $5$, we can get $b_2=0,\cdots, b_{m-1}=0$ successively, which
contradicts to $b_{m-1}\neq 0$. Hence $b_1\neq 0$. $\Box$

Next we give a brief discussion of the skew-symmetric case. It is
very similar to the symmetric case. Let $u_m$ be a
$\lambda_m$-eigenfunction of $-\Delta_m$ with
$\lambda_m\in\mathcal{P}_m^-$. Denote by $(b_0,b_1,b_2,\cdots,b_m)$
the values of $u_m$ on the skeleton of $\Omega_m$ where $b_0=b_m=0$
by the Dirichlet boundary condition.
 Write
$\lambda_{i}^{(m)}$ the successor of $\lambda_m$ of order $(m-i)$
with $2\leq i\leq m$.  Comparing to the symmetric case, the
eigenvalue equations at the vertex $F_1^iq_0$'s are unchanged except
the one at $F_1 q_0$, since now the values of $u_m$ on the four
$2$-level neighbors of $F_1 q_0$ are modified as shown in Fig. 4.3.
\begin{figure}[ht]
\begin{center}
\includegraphics[width=6cm,totalheight=5.8cm]{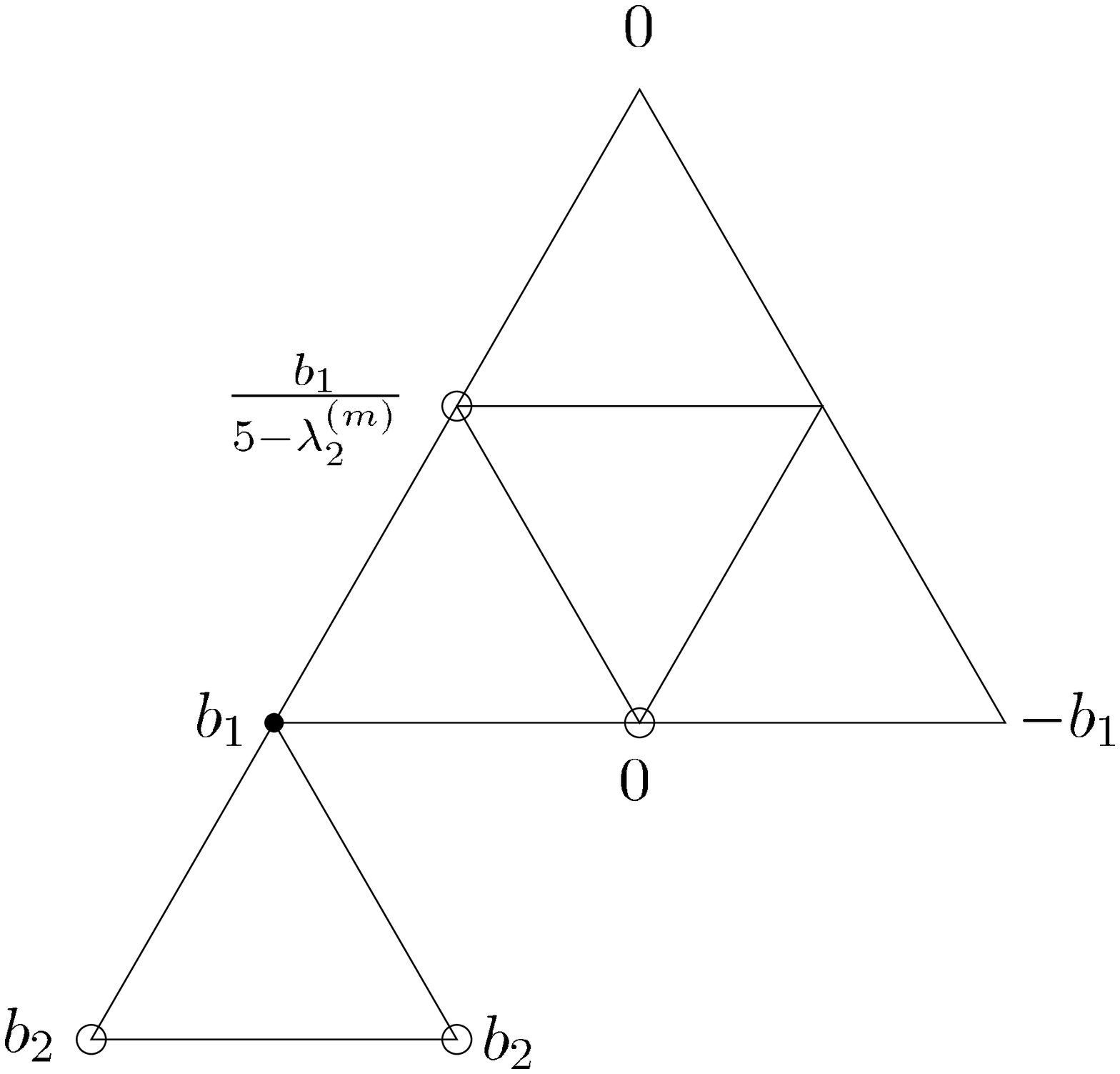}
\begin{center}
\textbf{Fig. 4.3.} \small{Values of $u_m$ on neighbors of $F_1q_0$.}
\end{center}
\end{center}
\end{figure}
 Hence we still have the same modified
eigenvalue equation
$$
l(\lambda_{i+1}^{(m)})b_{i-1}+s(\lambda_{i+1}^{(m)})b_i+r(\lambda_{i+1}^{(m)})b_{i+1}=0,\quad
\forall 2\leq i\leq m-1,
$$
while the first equation in $(\ref{9})$ is replaced by
$$\widetilde{s}(\lambda_{2}^{(m)})b_1+\widetilde{r}(\lambda_{2}^{(m)})b_2=0,$$
with $\widetilde{s}(x):=(4-x)(5-x)-1$ and
$\widetilde{r}(x):=-2(5-x)$. Now we assume $\lambda_{2}^{(m)}\neq 5$
and none of $\lambda_{i}^{(m)}$'s is equal to $2$ or $5$ for $3\leq
i\leq m$. (Later we will show this assumption automatically holds
for any $\lambda_m\in\mathcal{P}_m^-$.) Then by the eigenfunction
extension algorithm, $u_m$ is unique and determined by its values on
the skeleton of $\Omega_m$. Using similar discussion, $\lambda_m$
should be a solution of the following equation
\begin{equation}\label{12}
\widetilde{q}_m(x):=\left|
  \begin{array}{ccccc}
    \widetilde{s}(f^{(m-2)}(x)) & \widetilde{r}(f^{(m-2)}(x))  \\
    l(f^{(m-3)}(x)) & s(f^{(m-3)}(x)) & r(f^{(m-3)}(x))   \\
     & \ddots & \ddots & \ddots &  \\
     &  & l(f(x)) & s(f(x)) & r(f(x)) \\
     &  &  & l(x) & s(x) \\
  \end{array}
\right|=0,
\end{equation}
instead of $q_m(x)=0$ in the symmetric case. Hence if $\lambda_m$ is
a root of $\widetilde{q}_m(x)$, $f^{(m-2)}(\lambda_m)\neq 5$ and
none of $f^{(i)}(\lambda_m)$'s is  equal to $2$ or $5$ for $0\leq
i\leq m-3$, then $\lambda_{m}\in \mathcal{P}_m^-$. Similarly to
Proposition 4.1, we have

\textbf{Proposition 4.3.} \emph{(1) $\widetilde{q}_m(0)>0$, $\forall
m\geq 2$; }

 \emph{(2) $\widetilde{q}_2(5)<0$ and
$\widetilde{q}_m(5)>0$, $\forall m\geq 3$;}

 \emph{(3) $\widetilde{q}_2(6)>0$ and
$\widetilde{q}_m(6)<0$, $\forall m\geq 3$.}

\emph{Proof.} The first two statements follow from a very similar
argument in the proof of Proposition 4.1. We only need to prove the
third one.

It is easy to check that $\widetilde{q}_2(6)=1>0$ and
$\widetilde{q}_3(6)=-436<0$ by a direct computation. For $m\geq 4$,
an expansion along the first row yields that
$$\widetilde{q}_m(6)=\widetilde{s}(f^{(m-2)}(6))q_{m-1}(6)+2(5-f^{(m-2)}(6))(f^{(m-3)}(6)-6)q_{m-2}(6).$$
Recall that in the proof of Proposition 4.1(3), we have proved that
$q_{m-1}(6)\leq q_{m-2}(6)<0$. Hence
$$\widetilde{q}_m(6)\leq(\widetilde{s}(f^{(m-2)}(6))+2(5-f^{(m-2)}(6))(f^{(m-3)}(6)-6))q_{m-1}(6),$$
noticing that $f^{(m-2)}(6)<f^{(m-3)}(6)\leq -6$. An easy calculus
shows that
$\widetilde{s}(f^{(m-2)}(6))+2(5-f^{(m-2)}(6))(f^{(m-3)}(6)-6)\geq
1$, hence $\widetilde{q}_m(6)\leq q_{m-1}(6)<0$.  $\Box$

Similarly to the symmetric case, the following two lemmas focus on
the possibility of the roots of $\widetilde{q}_{m}(x)$ satisfying
$f^{(m-2)}(x)=5$, or $f^{(i)}(x)=2$ or $5$ for some $0\leq i\leq
m-3$.

\textbf{Lemma 4.7.} \emph{Let $x$ be a predecessor of $2$ of order
$i$ with $0\leq i\leq m-3$. Then $\widetilde{q}_m(x)=0$.}

\emph{Proof.}  If $0\leq i<m-3$, the proof is the same as that of
Lemma 4.2. So we only need to check the $i=m-3$ case. In this case,
$f^{(m-3)}(x)=2$ and $f^{(m-2)}(x)=6$. Substituting them into
$(\ref{12})$, noticing $s(2)=-8$, $l(2)=-4$, $r(2)=0$,
$\widetilde{s}(6)=1$ and $\widetilde{r}(6)=2$, we get
$$
\widetilde{q}_m(x)=\left|
  \begin{array}{ccccccc}
    \widetilde{s}(f^{(m-2)}(x)) &\widetilde{r}(f^{(m-2)}(x))  \\
           l(f^{(m-3)}(x)) & s(f^{(m-3)}(x)) & r(f^{(m-3)}(x)) \\
     &   \ddots &\ddots & \ddots \\
  \end{array}
\right|=\left|
  \begin{array}{ccccccc}
        1 & 2   \\
           -4 & -8 & 0 \\
   &    \ddots &\ddots & \ddots \\
  \end{array}
\right|=0.\quad \Box
$$

\textbf{Lemma 4.8.} \emph{Let $x$ be a predecessor of $5$ of order
$i$ with $0\leq i\leq m-2$. Then $\widetilde{q}_m(x)\neq 0$.}

\emph{Proof.} Similarly to the proof of Lemma 4.3, we only need to
prove $\widetilde{q}_m(\phi_-^{(i)}(5))\neq 0$. A similar argument
yields that $\widetilde{q}_m(\phi_-^{(0)}(5))=\widetilde{q}_m(5)$
and $\widetilde{q}_m(\phi_-^{(i)}(5))=\widetilde{q}_{m-i}(5)\cdot
q_{i+1}(\phi_-^{(i)}(5))$ for $0<i\leq m-2$. Combined with
Proposition 4.1(4) and Proposition 4.3(2), it follows the desired
result. $\Box$

Hence the total unwanted roots of $\widetilde{q}_m(x)$ consist of
those predecessors of $2$ of order $i$ with $0\leq i\leq m-3$ for
$m\geq 3$ and $\widetilde{q}_2(x)$ does not have any unwanted root.
This is exactly the same as the symmetric case. To exclude them out,
we define
$$\widetilde{p}_m(x):=\frac{\widetilde{q}_m(x)}{(x-2)(f(x)-2)\cdots (f^{(m-3)}(x)-2)},\quad \mbox{ for } m\geq
3,$$ and
$$\widetilde{p}_2(x):=\widetilde{q}_2(x)=\widetilde{s}(x).$$ These polynomials play a very similar role to $p_m(x)$'s in
the symmetric case. It is easy to check that the degree of
$\widetilde{p}_m(x)$ is $s_m:=2^m-2$, since the degree of the
polynomial $\widetilde{q}_m(x)$ is $3+3\cdot2+\cdots+3\cdot
2^{m-3}+2\cdot 2^{m-2}=3(2^{m-2}-1)+2^{m-1}$, and the number of all
the unwanted roots of $\widetilde{q}_m(x)$ is
$1+2+\cdots+2^{m-3}=2^{m-2}-1$ for $m\geq 3$ and $0$ for $m=2$.
 The following is a list of some facts about
$\widetilde{p}_m(x)$ similar to Proposition 4.2, which can be easily
get from Proposition 4.3.

 \textbf{Proposition 4.4.} \emph{Let $m\geq 2$, then}

 \emph{(1)
$(-1)^m\widetilde{p}_m(0)> 0$;}

\emph{(2) $(-1)^{m-1}\widetilde{p}_m(5)>0$;}

\emph{(3) $(-1)^{m}\widetilde{p}_m(6)>0$.}

 Then following a similar argument, the results in
Lemma 4.5 and Lemma 4.6 still hold with $\mathcal{P}_m^+$, $p_m(x)$
and $r_m$ replaced by $\mathcal{P}_m^-$, $\widetilde{p}_m(x)$ and
$s_m$ respectively.

Hence we have found $r_m$ distinct eigenvalues in $\mathcal{P}_m^+$
and $s_m$ distinct eigenvalues in $\mathcal{P}_m^-$. We will show
these eigenvalues are the totality of $\mathcal{P}_m$. To prove
this, the following lemma is needed.

\textbf{Lemma 4.9.} \emph{Let $\mathcal{P}_m^{+,*}$ and
$\mathcal{P}_m^{-,*}$ be the sets of total roots of $p_m(x)$ and
$\widetilde{p}_m(x)$ respectively. Let $\mathcal{M}_m^{*}$ be the
set of miniaturized eigenvalues generated by $\mathcal{P}_k^{-,*}$
with $2\leq k<m$. Let $\mathcal{L}_m$ denote the set of $m$-level
localized eigenvalues. Then all eigenfunctions associated to these
eigenvalues are linearly independent.}

\emph{Proof.} Without loss of generality, assume $m\geq 3$. It is
easy to check that for each $m$-level localized eigenfunction
$u_m^{\mathcal{L}}$, it must be $0$ on $\partial\Omega_{m-1}$. Lemma
4.6 says that each $m$-level symmetric primitive
$\lambda_m$-eigenfunction $u_m^{\mathcal{P,+}}$ with
$\lambda_m\in\mathcal{P}_m^{+,*}$ must be a non-zero constant on
$\partial\Omega_{m-1}\setminus\{q_0\}$ and be a non-zero constant on
$\partial\Omega_{1}\setminus\{q_0\}$. The skew-symmetric analog of
Lemma 4.6 says that each $m$-level skew-symmetric primitive
$\lambda_m$-eigenfunction $u_m^{\mathcal{P,-}}$ with
$\lambda_m\in\mathcal{P}_m^{-,*}$ must be a non-zero constant on
each symmetric part of $\partial\Omega_{m-1}\setminus\{q_0\}$ under
the symmetry fixing $q_0$, and take non-zero value on $F_1q_0$ and
$F_2q_0$ only different in signs. From the construction of the
miniaturized eigenfunctions, for each $m$-level miniaturized
eigenfunction $u_m^{\mathcal{M}}$ with eigenvalue in
$\mathcal{M}_m^{*}$, $u_m^{\mathcal{M}}$ must take non-zero value on
a subset of $\partial\Omega_{m-1}\setminus\{q_0\}$ and be $0$ on
$\partial\Omega_1$.
 These observations implies the linearly independence of
eigenfunctions among different types. $\Box$

Hence we have

\textbf{Lemma 4.10.}\emph{ For each $m\geq 2$,
$\mathcal{P}_m^{+,*}=\mathcal{P}_m^+$ and
$\mathcal{P}_m^{-,*}=\mathcal{P}_m^-$.}

\emph{Proof.} Lemma 4.5 and its skew-symmetric analog say that
$\sharp\mathcal{P}_m^{+,*}=r_m$ and $\sharp\mathcal{P}_m^{-,*}=s_m$.
Accordingly, a similar argument in the proof of Theorem 3.4 shows
that $\sharp\mathcal{M}_m^{*}=(m-3)\cdot 2^m+4$.

Then it is easy to check that $\sharp \mathcal{L}_m$, $\sharp
\mathcal{P}_m^{+,*}$, $\sharp \mathcal{P}_m^{-,*}$ and $\sharp
\mathcal{M}_m^{*}$ add up to $\sharp(V_m^\Omega\setminus\partial
\Omega_m)$ with a suitable modification of the eigenspace
dimensional counting formula $(\ref{11})$. Using Lemma 4.9, we then
get the desired result. $\Box$

By this lemma, it is easy to see that the assumptions we made before
on symmetric and skew-symmetric $m$-level primitive eigenvalues hold
automatically.

Next we will prove each primitive eigenvalue $\lambda\in
\mathcal{P}_m$ has multiplicity $1$. For this purpose, we need the
following lemma.

\textbf{Lemma 4.11.} \emph{For each $m\geq 2$, $\mathcal{P}_m^+\cap
\mathcal{P}_{m+1}^+=\emptyset$.}

\emph{Proof.} For $m=2$, it can be checked by a direct computation.
In order to use the induction, we assume that $\mathcal{P}_m^+\cap
\mathcal{P}_{m+1}^+=\emptyset$ and will prove
$\mathcal{P}_{m+1}^+\cap \mathcal{P}_{m+2}^+=\emptyset$.

Suppose there is a $\lambda\in\mathcal{P}_{m+1}^+\cap
\mathcal{P}_{m+2}^+$. Then $p_{m+1}(\lambda)=p_{m+2}(\lambda)=0$
(hence $q_{m+1}(\lambda)=q_{m+2}(\lambda)=0$). Moreover, none of
$f^{(i)}(\lambda)$ ($0\leq i\leq m$) is equal to $2$ or $5$.

The expansion along the first row of $q_{m+2}(\lambda)$ gives
$$q_{m+2}(\lambda)=s(f^{(m)}(\lambda))q_{m+1}(\lambda)-r(f^{(m)}(\lambda))l(f^{(m-1)}(\lambda))q_m(\lambda).$$
Noticing $q_{m+1}(\lambda)=q_{m+2}(\lambda)=0$, we have
$$r(f^{(m)}(\lambda))l(f^{(m-1)}(\lambda))q_m(\lambda)=-2(2-f^{(m)}(\lambda))(5-f^{(m)}(\lambda))(f^{(m-1)}(\lambda)-6)q_m(\lambda)=0.$$
Hence $q_m(\lambda)=0$ or $f^{(m-1)}(\lambda)=6$, since
$f^{(m)}(\lambda)\neq 2$ or $5$.

If $q_m(\lambda)=0$, then $\lambda\in \mathcal{P}_m^+$, hence
$\mathcal{P}_m^+\cap \mathcal{P}_{m+1}^+\neq\emptyset$. This
contradicts to our induction assumption.

Hence we have $f^{(m-1)}(\lambda)=6$, i.e., $f^{(m-2)}(\lambda)=3$.
Noticing that $\lambda$ is also a root of $p_{m+1}(x)$,  Lemma 4.1
says that $\phi_-^{(m-2)}(3)$ is a root of $p_{m+1}(x)$. Hence
$p_{m+1}(\phi_-^{(m-2)}(3))=0$, which contradicts to Proposition
4.1(6). Hence such $\lambda$ can not exist. So we get the desired
result. $\Box$

Then we can prove:

\textbf{Lemma 4.12.} \emph{For each $m\geq 2$,
$\mathcal{P}_m^+\cap\mathcal{P}_m^-=\emptyset$.}

\emph{Proof.} For $m=2$ or $3$, it can be checked by a direct
computation. Let $m\geq 4$. Suppose there is an eigenvalue
$\lambda_m\in \mathcal{P}_m^+\cap\mathcal{P}_m^-$. Then by Lemma
4.10, $p_m(\lambda_m)=\widetilde{p}_m(\lambda_m)=0$. For each $2\leq
i\leq m$, denote by $\lambda_{i}^{(m)}$ the successor of $\lambda_m$
of order $(m-i)$. Obviously we have
$q_m(\lambda_m)=\widetilde{q}_m(\lambda_m)=0$ and
$\lambda_{i}^{(m)}\neq 2$ or $5$ for $2\leq i\leq m$. Furthermore,
by Lemma 4.11, we have $p_{m-1}(\lambda_m)\neq 0$, hence
$q_{m-1}(\lambda_m)\neq 0$.

Using the expansions of $q_m(\lambda_m)$ and
$\widetilde{q}_{m}(\lambda_m)$ along their first rows respectively,
we have
$$s(\lambda_2^{(m)})q_{m-1}(\lambda_m)-r(\lambda_2^{(m)})l(\lambda_3^{(m)})q_{m-2}(\lambda_m)=0$$
and
$$\widetilde{s}(\lambda_2^{(m)})q_{m-1}(\lambda_m)-\widetilde{r}(\lambda_2^{(m)})l(\lambda_3^{(m)})q_{m-2}(\lambda_m)=0.$$
Hence, the vector $(q_{m-1}(\lambda_m), q_{m-2}(\lambda_m))$ can be
viewed as a non-zero solution of the system of linear equations,
$$\left\{\begin{array}{c}
    s(\lambda_{2}^{(m)})x-r(\lambda_{2}^{(m)})l(\lambda_{3}^{(m)})y=0 \\
    \widetilde{s}(\lambda_{2}^{(m)})x-\widetilde{r}(\lambda_{2}^{(m)})l(\lambda_{3}^{(m)})y=0.
  \end{array}\right.
$$
Thus $$\left|
        \begin{array}{cc}
          s(\lambda_{2}^{(m)}) & 2(2-\lambda_{2}^{(m)})(5-\lambda_{2}^{(m)})(\lambda_{3}^{(m)}-6) \\
          \widetilde{s}(\lambda_{2}^{(m)}) & 2(5-\lambda_{2}^{(m)})(\lambda_{3}^{(m)}-6) \\
        \end{array}
      \right|=0.
$$ Since $\lambda_2^{(m)}\neq 5$, we have $\lambda_3^{(m)}=6$ or
$s(\lambda_{2}^{(m)})=(2-\lambda_{2}^{(m)})\widetilde{s}(\lambda_{2}^{(m)})$.
By Substituting the expressions for  $s(x)$ and $\widetilde{s}(x)$,
we get
 $\lambda_{2}^{(m)}=6$ or $\lambda_{3}^{(m)}=6$. Hence we have
 $\lambda_3^{(m)}=3$ or  $\lambda_4^{(m)}=3$, i.e.,
 $f^{(m-3)}(\lambda_m)=3$, or $f^{(m-4)}(\lambda_m)=3$.

 Noticing
 that $\lambda_m$ is a root of $q_m(x)$, by using Lemma 4.1, we can see that either $\phi_{-}^{(m-3)}(3)$ or
 $\phi_{-}^{(m-4)}(3)$ is a root of $q_m(x)$, i.e.,
 $q_m(\phi_{-}^{(m-3)}(3))=0$ or $q_m(\phi_{-}^{(m-4)}(3))=0$.
An expansion of $q_m(\phi_{-}^{(m-4)}(3))$ along the first row
yields that
$$q_m(\phi_{-}^{(m-4)}(3))=s(f^{(2)}(3))q_{m-1}(\phi_{-}^{(m-4)}(3))=848q_{m-1}(\phi_{-}^{(m-4)}(3))$$
since $l(f(3))=0$. Hence we have either $q_m(\phi_{-}^{(m-3)}(3))=0$
or $q_{m-1}(\phi_{-}^{(m-4)}(3))=0$. By Proposition 4.1(6), this is
impossible. Hence such $\lambda_m$ can not exist. So
$\mathcal{P}_m^+\cap\mathcal{P}_m^-=\emptyset$. $\Box$

\emph{Proof of  Theorem 3.3 and Theorem 3.5.}

It is an immediate consequence, by using Lemma 4.5 and its
skew-symmetric analog, Lemma 4.9, Lemma 4.10, Lemma 4.12 and the
eigenspace dimension counting formula $(\ref{11})$. $\Box$

\section{Primitive Dirichlet eigenvalues of $-\Delta$}

Having found the primitive Dirichlet eigenvalues and eigenfunctions
for $-\Delta_m$, it is natural to believe that the primitive
Dirichlet eigenvalues of $-\Delta$ could be obtained in the limit as
$m$ goes to infinity. This is true for the spectrum for
$\mathcal{SG}\setminus{V_0}$ case, benefiting from the spectral
decimation method and the eigenfunction extension algorithm
$(\ref{3})$. Our goal in this section is to extend this recipe to
$\Omega$ case by instead using the weak spectral decimation
introduced in Section 3. Comparing to the
$\mathcal{SG}\setminus{V_0}$ case, our method is more based on
estimates. We focus on the symmetric case, since the skew-symmetric
case can be got by using a similar discussion. We will prove Theorem
3.6 in this section.

We use the $\widetilde{\phi}_{\pm}$ notations introduced in Section
3. Recall that if $\alpha_{m},\beta_m$ are two consecutive
eigenvalues in $\mathcal{P}_m^+$ with $\alpha_{m}<\beta_{m}$, then
we always have
\begin{equation}\label{51}
\phi_-(\alpha_m)<\widetilde{\phi}_-(\beta_m)<\phi_-({\beta_m})\mbox{
and }
\phi_+(\beta_m)<\widetilde{\phi}_+(\beta_m)<\phi_+({\alpha_m}),
\end{equation} and if $\beta_m$ is the least eigenvalue in $\mathcal{P}_m^+$,
then instead we have
$$
0<\widetilde{\phi}_-(\beta_m)<\phi_-(\beta_m) \mbox{ and }
\phi_+(\beta_m)<\widetilde{\phi}_+(\beta_m)<5.$$

Let $m_0\geq 2$, $\lambda_{m_0}$ be a $m_0$-level symmetric
primitive eigenvalue, $\{\lambda_m\}_{m\geq m_0}$ be an infinite
sequence related by $\lambda_{m+1}=\widetilde{\phi}_-(\lambda_m)$ or
$\widetilde{\phi}_+(\lambda_m)$, $\forall m\geq m_0$, assuming that
there are only a finite number of $\widetilde{\phi}_+$ relations.
Call the minimum value $m_1$, such that $\forall m\geq m_1$,
$\lambda_{m+1}=\widetilde{\phi}_-(\lambda_m)$, the \emph{generation
of fixation} of the sequence $\{\lambda_{m}\}_{m\geq m_0}$. In all
that follows in this section, we always use $\{\lambda_{m}\}_{m\geq
m_0}$ as such a sequence without specifical declaration.

The first fact about this sequence is:

\textbf{Lemma 5.1.} \emph{$\lim_{m\rightarrow \infty}5^m
\lambda_{m}$ exists.}

\emph{Proof.} Without loss of generality, assume $\lambda_{m_1}<5$,
otherwise, we could choose $\widetilde{m}_1=m_1+1$ and use
$\widetilde{m}_1$ to replace $m_1$ in the following proof.

Let $m\geq m_1$, then
$\frac{\lambda_{m+1}}{\lambda_m}=\frac{\widetilde{\phi}_-(\lambda_m)}{\lambda_m}\leq\frac{\phi_-(\lambda_m)}{\lambda_m}=\frac{\phi_-(\lambda_m)}{\phi_-(\lambda_m)(5-\phi_-(\lambda_m))}=\frac{1}{5-\phi_-(\lambda_m)}.$
Since $0<\lambda_m<5$, we have $0<\phi_-(\lambda_m)<2$, hence
$\frac{1}{5-\phi_-(\lambda_m)}<\frac{1}{3}$. Thus $\sum_{m\geq
m_1}\lambda_m<\infty.$

Furthermore, $\frac{5^{m+1}\lambda_{m+1}}{5^m
\lambda_m}=5\frac{\lambda_{m+1}}{\lambda_m}\leq\frac{5}{5-\phi_-(\lambda_m)}=1+\frac{\phi_-(\lambda_m)}{5-\phi_-(\lambda_m)}$.
Noticing that $\sum_{m\geq
m_1}\frac{\phi_-(\lambda_m)}{5-\phi_-(\lambda_m)}\leq\frac{1}{3}\sum_{m\geq
m_1}\phi_-(\lambda_m)\leq \frac{1}{3}\sum_{m\geq
m_1}\lambda_m<\infty$ since $\phi_-'(x)<1$ whenever $0<x<5$, we get
that $\Pi_{m\geq m_1}\frac{5^{m+1}\lambda_{m+1}}{5^m\lambda_m}$
converges. Hence $\lim_{m\rightarrow\infty}5^m\lambda_m$ exists.
$\Box$

The following is an estimate of the difference between
$\widetilde{\phi}_-(\lambda_m)$ and $\phi_-(\lambda_m)$ for
$\lambda_m$ in the sequence $\{\lambda_{m}\}_{m\geq m_0}$.

\textbf{Proposition 5.1.} \emph{
$$\sum_{m\geq m_1}5^m(\widetilde{\phi}_{-}(\lambda_m)-\phi_-(\lambda_m))<\infty.$$
In particular,
$\lim_{m\rightarrow\infty}5^m(\widetilde{\phi}_-(\lambda_m)-\phi_-(\lambda_m))=0$.}

\emph{Proof.} Without loss of generality, assume $\lambda_{m_1}<5$.
 From
Lemma 5.1, we have  $\sum_{m\geq
m_1}(5^{m+1}\lambda_{m+1}-5^m\lambda_m)<\infty$. Hence
\begin{eqnarray}\label{52}
&&\sum_{m\geq m_1} 5^m(5\lambda_{m+1}-\lambda_m)\nonumber\\
=&& \sum_{m\geq
m_1}5^m(5\widetilde{\phi}_-(\lambda_m)-\phi_-(\lambda_m)(5-\phi_-(\lambda_m)))\nonumber\\
=&&\sum_{m\geq
m_1}(5^{m+1}(\widetilde{\phi}_-(\lambda_m)-\phi_-(\lambda_m))+5^m(\phi_-(\lambda_m))^2)<\infty.
\end{eqnarray}

Since $0<\phi_-'(x)<1$ whenever $0<x<5$, we have
$5^m(\phi_-(\lambda_m))^2\leq 5^m\lambda_m^2$. Still from Lemma 5.1,
we have $\lambda_m=O(\frac{1}{5^m})$, hence
$5^m(\phi_-(\lambda_m))^2\leq \frac{c}{5^m}$ for some constant $c$.
Thus $\sum_{m\geq m_1}5^m(\phi_-(\lambda_m))^2<\infty$. Combining
this with $(\ref{52})$, we get $\sum_{m\geq
m_1}5^m(\widetilde{\phi}_-(\lambda_m)-\phi_-(\lambda_m))<\infty.$
$\Box$

To reveal some further properties of the limit
$\lim_{m\rightarrow\infty}5^m\lambda_m$, the following lemma is
required, which is a generalization of formula $(\ref{51})$.

 \textbf{Lemma 5.2.} \emph{Let $m\geq 2$. $\alpha_{m}, \beta_{m}$ be two
consecutive  eigenvalues in $\mathcal{P}_m^+$ with
$\alpha_m<\beta_m$. Then $\forall l\in \mathbb{N}$,
\begin{equation}\label{54}
\phi_-^{(l)}(\alpha_m)<\widetilde{\phi}_-^{(l)}(\beta_m).
\end{equation}}

\emph{Proof.} First we need to prove the following relation.
\begin{equation}\label{53}
p_{m+l}(\phi_-^{(l)}(\alpha_m))\sim(-1)^{l-1}p_{m+1}(\phi_-(\alpha_m)),
\quad \forall l\in\mathbb{N}.
\end{equation}

In fact, when $l\geq 3$, using the Laplace theorem to expand the
determinant $q_{m+l}(\phi_-^{(l)}(\alpha_m))$ according to the last
$(l-1)$ rows, we have
$$q_{m+l}(\phi_-^{(l)}(\alpha_m))=q_{l}(\phi_-^{(l)}(\alpha_m)) q_{m+1}(\phi_-(\alpha_m))-l(\phi_-^{(2)}(\alpha_m))q_{l-1}(\phi_-^{(l)}(\alpha_m))r(\phi_-(\alpha_m))q_{m}(\alpha_m).$$
Since $q_m(\alpha_m)=0$, we have
$$q_{m+l}(\phi_-^{(l)}(\alpha_m))=q_{l}(\phi_-^{(l)}(\alpha_m))
q_{m+1}(\phi_-(\alpha_m)).$$ This equality also holds for $l=2$ by
instead using an expansion along the last row of
$q_{m+2}(\phi_-^{(2)}(\alpha_m))$. Hence for each $l\geq 2$, we
always have
$q_{m+l}(\phi_-^{(l)}(\alpha_m))=q_{l}(\phi_-^{(l)}(\alpha_m))
q_{m+1}(\phi_-(\alpha_m))$. Then from Lemma B in Appendix, we have
$q_l(\phi_-^{(l)}(\alpha_m))>0$, hence
$q_{m+l}(\phi_-^{(l)}(\alpha_m))\sim q_{m+1}(\phi_-(\alpha_m))$. By
the relation between $p_{m+l}(x)$ and $q_{m+l}(x)$, we can easily
get $(\ref{53})$.

Now we prove $(\ref{54})$. When $l=1$, $(\ref{54})$ follows from
$(\ref{51})$ directly. In order to use the induction, assuming
$(\ref{54})$ holds for $l$, we turn to prove$$
\phi_-^{(l+1)}(\alpha_m)<\widetilde{\phi}_-^{(l+1)}(\beta_m).
$$
Suppose $\alpha_m$ and $\beta_m$ are  the $k$'th and $(k+1)$'th
eigenvalues in $\mathcal{P}_m^+$ respectively. Recall that in Lemma
4.4, we have proved that $p_{m+1}(\phi_-(\alpha_m))\sim
(-1)^{m+k-1}$. Combining this with $(\ref{53})$, we have
\begin{equation}\label{55}
p_{m+l+1}(\phi_-^{(l+1)}(\alpha_m))\sim (-1)^{m+k+l-1}.
\end{equation}

On the other hand, if we denote
$\alpha_{m+l}=\widetilde{\phi}_-^{(l)}(\alpha_m)$ and
$\beta_{m+l}=\widetilde{\phi}_-^{(l)}(\beta_m)$, then it is easy to
see that $\alpha_{m+l}$ and $\beta_{m+l}$ are the $k$'th and
$(k+1)$'th eigenvalues in $\mathcal{P}_{m+l}^+$ respectively. Lemma
4.4 says that
\begin{equation}\label{56}
p_{m+l+1}(\phi_-(\alpha_{m+l}))\sim (-1)^{m+l+k-1}
\end{equation}
and
\begin{equation}\label{57}
p_{m+l+1}(\phi_-(\beta_{m+l}))\sim (-1)^{m+l+k}.
\end{equation}
Furthermore, if we denote
$\beta_{m+l+1}=\widetilde{\phi}_-^{(l+1)}(\beta_m)$, then
$\beta_{m+l+1}$ is the only root of $p_{m+l+1}(x)$ located between
$\phi_-(\alpha_{m+l})$ and $\phi_-(\beta_{m+l})$, i.e.,
\begin{equation}\label{58}
\phi_-(\alpha_{m+l})<\beta_{m+l+1}<\phi_-(\beta_{m+l}).
\end{equation}

Noticing that from the induction assumption, we have
$\phi_-^{(l+1)}(\alpha_m)<\phi_-(\beta_{m+l})$ since
$\beta_{m+l}=\widetilde{\phi}_{-}^{(l)}(\beta_m)$. Moreover,
$(\ref{55})$ and $(\ref{57})$ say that there exists at least one
root of $p_{m+l+1}(x)$, denoted by $\beta_{m+l+1}^*$,  between
$\phi_-^{(l+1)}(\alpha_m)$ and $\phi_-(\beta_{m+l})$, i.e.,
\begin{equation}\label{510}
\phi_-^{(l+1)}(\alpha_m)<\beta_{m+l+1}^*<\phi_-(\beta_{m+l}).
\end{equation}

Since
$\phi_-(\alpha_{m+l})=\phi_-(\widetilde{\phi}_-^{(l)}(\alpha_m))<\phi_-^{(l+1)}(\alpha_m)$,
we have
$$\phi_-(\alpha_{m+l})<\phi_-^{(l+1)}(\alpha_m)<\beta_{m+l+1}^*<\phi_-(\beta_{m+l}).$$
Combing this with $(\ref{58})$, from the uniqueness of
$\beta_{m+l+1}$, we have $\beta_{m+l+1}=\beta_{m+l+1}^*$. Hence
substituting it into $(\ref{510})$, we finally get
$\phi_-^{(l+1)}(\alpha_m)<\beta_{m+l+1}$, i.e.,
$\phi_-^{(l+1)}(\alpha_m)<\widetilde{\phi}_-^{(l+1)}(\beta_m)$,
which is the desired result. $\Box$

The following is an application of Lemma 5.2.

\textbf{Lemma 5.3.} \emph{Let $m_1\geq 2$, $\alpha_{m_1},
\beta_{m_1}$ be two consecutive  eigenvalues in
$\mathcal{P}_{m_1}^+$ with $\alpha_{m_1}<\beta_{m_1}$.
$\{\alpha_m\}_{m\geq m_1}$ is an infinite sequence related by
$\alpha_{m+1}=\widetilde{\phi}_-(\alpha_m),\forall m\geq m_1$;
$\{\beta_m\}_{m\geq m_1}$ is an infinite sequence related by
$\beta_{m+1}=\widetilde{\phi}_-(\beta_m),\forall m\geq m_1$. Then
$\forall m\geq m_1$, $\alpha_m<\beta_m$. Moreover,
$$\lim_{m\rightarrow\infty} 5^m\alpha_m<\lim_{m\rightarrow\infty} 5^m\beta_m.$$}

\textbf{Remark.} \emph{In $\mathcal{SG}\setminus{V_0}$ case, this is
a direct result since $\phi_-(x)$ is a definite strictly increasing
continuous function.}

\emph{Proof of Lemma 5.3.} Let $m> m_1$. Since
$\alpha_m=\widetilde{\phi}_-^{(m-m_1)}(\alpha_{m_1})$ and
$\beta_m=\widetilde{\phi}_-^{(m-m_1)}(\beta_{m_1})$, we have
\begin{equation}\label{511}
\alpha_m<\phi_-^{(m-m_1)}(\alpha_{m_1})<\widetilde{\phi}_-^{(m-m_1)}(\beta_{m_1})=\beta_m
\end{equation}
by Lemma 5.2. Hence $\forall m>m_1$, $\alpha_m<\beta_m$.

Now we prove $\lim_{m\rightarrow\infty}
5^m\alpha_m<\lim_{m\rightarrow\infty} 5^m\beta_m.$

Let $m> m_1$. Then from $(\ref{511})$, we have
$$\alpha_m<\phi_-^{(m-m_1-1)}(\widetilde{\phi}_-(\alpha_{m_1}))<\phi_-^{(m-m_1)}(\alpha_{m_1})<\beta_m.$$
Hence
$\beta_m-\alpha_m>\phi_-^{(m-m_1-1)}(\phi_-(\alpha_{m_1}))-\phi_-^{(m-m_1-1)}(\widetilde{\phi}_-(\alpha_{m_1})).$
Since $\phi_-'(x)\geq \frac{1}{5}$ whenever $0<x<5$, and
$0<\widetilde{\phi}_-(\alpha_{m_1})<\phi_-(\alpha_{m_1})<5$, we have
$$\beta_m-\alpha_m>\frac{1}{5^{m-m_1-1}}(\phi_-(\alpha_{m_1})-\widetilde{\phi}_-(\alpha_{m_1})).$$
Hence
$5^m(\beta_m-\alpha_m)>5^{m_1+1}(\phi_-(\alpha_{m_1})-\widetilde{\phi}_-(\alpha_{m_1}))$
which yields that
$$\lim_{m\rightarrow\infty}5^m(\beta_m-\alpha_m)\geq 5^{m_1+1}(\phi_-(\alpha_{m_1})-\widetilde{\phi}_-(\alpha_{m_1}))>0.$$
Thus $\lim_{m\rightarrow\infty}
5^m\alpha_m<\lim_{m\rightarrow\infty} 5^m\beta_m.$ $\Box$

\textbf{Lemma 5.4.} \emph{
$\lim_{m\rightarrow\infty}5^m\lambda_m>0.$}

\textbf{Remark.} \emph{In $\mathcal{SG}\setminus{V_0}$ case, this is
also a direct result, since $\{5^m\lambda_m\}_{m\geq m_1}$ is then a
monotone increasing sequence.}

\emph{Proof of Lemma 5.4.} Without loss of generality, we assume
that $\lambda_{m_1}$ is the least eigenvalue in
$\mathcal{P}_{m_1}^+$, since Lemma 5.3 says that it suffices  to
prove for this special case. Then $\forall m\geq m_1$, $\lambda_m$
is also the least eigenvalue in $\mathcal{P}_{m}^+$. Note that Lemma
B in Appendix says that $\forall m\geq m_1$, we have $\lambda_m\geq
\phi_-^{(m)}(6).$ Hence
$$\lim_{m\rightarrow\infty}5^m\lambda_m\geq \lim_{m\rightarrow\infty}5^m\phi_-^{(m)}(6)>0,$$
where the existence and positivity  of the second limit are already
shown in $\mathcal{SG}\setminus{V_0}$ case. See $\cite{Fu}$. $\Box$

Now we define
$$\lambda=\frac{3}{2}\lim_{m\rightarrow\infty}5^m\lambda_m.$$ We will
prove $\lambda$ is an primitive Dirichlet eigenvalue of $-\Delta$ on
the fractal domain $\Omega$.

Note that $\forall m\geq m_0$, $\lambda_m\in \mathcal{P}_m^+$, i.e.,
$\lambda_m$ is a root of both $p_m(x)$ and $q_m(x)$ by Lemma 4.5 and
Theorem 3.3. As in Section 4, denote by $\lambda_{i}^{(m)}$ the
successor of $\lambda_m$ of order $(m-i)$ with $2\leq i\leq m$ for
simplicity. Lemma 4.6 and Theorem 3.3 say that the system
$(\ref{9})$ of equations has $1$-dimensional solutions
$(b_1,b_2,\cdots,b_{m-1})$ with $b_1\neq 0$ and $b_{m-1}\neq 0$. We
normalize the solution by requiring $b_1=1$, and write it as
$(b_1^{(m)},b_2^{(m)},\cdots,b_{m-1}^{(m)})$ with $b_1^{(m)}=1$ to
specify its relation to $\lambda_m$. We always denote $b_0^{(m)}=0$
for convenience. As described in Section 4, from
$(b_1^{(m)},b_2^{(m)},\cdots,b_{m-1}^{(m)})$ one can recover the
unique (up to a constant) $\lambda_m$-eigenfunction $u_m$ on
$\Omega_m$ (noticing that $\lambda_{i}^{(m)}\neq 2$ or $5$, $\forall
2\leq i\leq m$). Hence
\begin{eqnarray*}\left\{
\begin{array}{l}
  -\Delta_m u_m=\lambda_m u_m \mbox{ on } \Omega_m,\\
  u_m|_{\partial \Omega_m}=0.
\end{array}\right.
\end{eqnarray*}

For each $m\geq m_0$, we start with the $\lambda_m$-eigenfunction
$u_m$ on $\Omega_m$, and extend $u_m$ to $\Omega$ by successively
using the  eigenfunction extension algorithm $(\ref{3})$
corresponding to the revised eigenvalue sequence
$\{\lambda_m,\phi_-(\lambda_m),\phi_-^{(2)}(\lambda_m),\cdots\}$
(starting from $\lambda_m$, but continued with the standard spectral
decimation eigenvalues) to get a primitive eigenfunction (possessing
the symmetry in each cell $F_w(\mathcal{SG})$ under the reflection
symmetry fixing $F_wq_0$ with word $w$ taking symbols only from
$\{1,2\}$) on $\Omega$. We still denote $u_m$ for this function. Of
cause, $u_m$ may not satisfy the Dirichlet boundary condition on
$L$. $\forall i>m$, we use
$\lambda_{i}^{(m)}=\phi_-^{(i-m)}(\lambda_m)$ to denote the
$i$-level revised eigenvalue. Hence for each $m\geq m_0$, $u_m$ is
an eigenfunction of $-\Delta$ on $\Omega$(not satisfying the
Dirichlet boundary condition), associated to the eigenvalue sequence
$\{\lambda_{i}^{(m)}\}_{i\geq 2}$, where
$\lambda_i^{(m)}=f^{(m-i)}(\lambda_m)$, $\forall 2\leq i\leq m$, and
$\lambda_{i}^{(m)}=\phi_-^{(i-m)}(\lambda_m)$, $\forall i>m$. We use
$b_i^{(m)}$ ($\forall i\geq m$) to denote the value of $u_m$ at
vertex $F_1^iq_0$. Hence $\{b_i^{(m)}\}_{i\geq 0}$ are the values of
$u_m$ on the skeleton of $\Omega$ which conversely determine $u_m$
on $\Omega$. We have the following relationship between
$\{\lambda_i^{(m)}\}_{i\geq 2}$ and $\{b_i^{(m)}\}_{i\geq 0}$.
\begin{equation}\label{512}
(4-\lambda_{i+1}^{(m)})b_i^{(m)}=2b_{i+1}^{(m)}+\frac{(14-3\lambda_{i+1}^{(m)})b_i^{(m)}+(6-\lambda_{i+1}^{(m)})b_{i-1}^{(m)}}{(2-\lambda_{i+1}^{(m)})(5-\lambda_{i+1}^{(m)})},\quad
\forall i\geq 1,
\end{equation}
which follows from the eigenvalue equation at the vertex $F_1^iq_0$.
Note  that  when $1\leq i\leq m-1$, these are exactly the equations
in $(\ref{91})$. Moreover, $u_m$ on $\Omega$ satisfies that
\begin{eqnarray*}\left\{
\begin{array}{l}
  -\Delta u_m=5^m\Phi(\lambda_m) u_m \mbox{ on } \Omega,\\
  u_m(q_0)=0,\\
  u_m|_L=\lim_{i\rightarrow\infty}b_i^{(m)}<\infty,
\end{array}\right.
\end{eqnarray*}
where $\Phi(z)$ is a function defined by
$\Phi(z):=\frac{3}{2}\lim_{k\rightarrow\infty}5^k\phi_{-}^{(k)}(z)$.
The existence  of the limit $\lim_{i\rightarrow\infty}b_i^{(m)}$
will be given later.

It is easy to find that $5^m\Phi(\lambda_m)\rightarrow\lambda$ as
$m$ goes to infinity. Moreover, we have the following lemmas.

\textbf{Lemma 5.5.} \emph{There exists a constant $C_1>0$ depending
only on $m_1$, such that $\forall i\in \mathbb{N}$, $\forall p\in
\mathbb{N}$, we have $|b_{i+p}^{(m)}-b_i^{(m)}|\leq
C_1(\frac{3}{10})^i\| u_m\|_\infty$  uniformly on $m\geq m_1$.}

\emph{Proof.} Without loss of generality, assume $i>m_1$ and
$\lambda_{m_1}$ is not the largest eigenvalue in
$\mathcal{P}_{m_1}^+$. Denote by $\gamma_{m_1}$ the next eigenvalue
of $\lambda_{m_1}$ in $\mathcal{P}_{m_1}^+$. Let
$\{\gamma_m\}_{m\geq m_1}$ be the infinite sequence staring from
$\gamma_{m_1}$ related by
$\gamma_{m+1}=\widetilde{\phi}_-(\gamma_m)$, $\forall m\geq m_1$. We
now show
\begin{equation}\label{513}
\lambda_{i+1}^{(m)}<\gamma_{i+1}<\phi_-(2), \quad \forall m\geq m_1.
\end{equation}

In fact if $m\geq i+1$, then
$$\lambda_{i+1}^{(m)}=f^{(m-i-1)}(\lambda_{m})=f^{(m-i-1)}(\widetilde{\phi}_-^{(m-i-1)}(\lambda_{i+1}))\leq
f^{(m-i-1)}(\phi_-^{(m-i-1)}(\lambda_{i+1}))=
\lambda_{i+1}<\gamma_{i+1}.$$ If $m<i+1$, then
$\lambda_{i+1}^{(m)}=\phi_-^{(i+1-m)}(\lambda_m)<\widetilde{\phi}_-^{(i+1-m)}(\gamma_m)=\gamma_{i+1}$
by using Lemma 5.2. The right inequality of $(\ref{513})$ is
obvious. Hence $(\ref{513})$ always holds.

On the other hand, notice that from $(\ref{512})$,
\begin{eqnarray*}
b_{i+1}^{(m)}-b_{i}^{(m)}&=&\frac{s(\lambda_{i+1}^{(m)})b_{i}^{(m)}-(6-\lambda_{i+1}^{(m)})b_{i-1}^{(m)}}{2(2-\lambda_{i+1}^{(m)})(5-\lambda_{i+1}^{(m)})}-b_i^{(m)}
\\&=&\frac{(6-\lambda_{i+1}^{(m)})(b_i^{(m)}-b_{i-1}^{(m)})-(20\lambda_{i+1}^{(m)}-9(\lambda_{i+1}^{(m)})^2+(\lambda_{i+1}^{(m)})^3)b_i^{(m)}}{2(2-\lambda_{i+1}^{(m)})(5-\lambda_{i+1}^{(m)})}.
\end{eqnarray*}
Hence
$$|b_{i+1}^{(m)}-b_{i}^{(m)}|\leq\frac{|6-\lambda_{i+1}^{(m)}|}{2|2-\lambda_{i+1}^{(m)}|\cdot|5-\lambda_{i+1}^{(m)}|}|b_i^{(m)}-b_{i-1}^{(m)}|+\frac{|20-9\lambda_{i+1}^{(m)}+(\lambda_{i+1}^{(m)})^2|}{2|2-\lambda_{i+1}^{(m)}|\cdot|5-\lambda_{i+1}^{(m)}|}|\lambda_{i+1}^{(m)}|\cdot|b_i^{(m)}|.$$

In the remaining proof, we use $c$ to denote different constants.

By $(\ref{513})$, we have
$$|b_{i+1}^{(m)}-b_{i}^{(m)}|\leq\frac{3}{(2-\gamma_{i+1})(5-\gamma_{i+1})}|b_i^{(m)}-b_{i-1}^{(m)}|+c\gamma_{i+1}|b_i^{(m)}|.$$
Noticing that $\gamma_i=O(\frac{1}{5^i})$ and $|b_i^{(m)}|\leq\|
u_m\|_\infty$, we get
$$|b_{i+1}^{(m)}-b_{i}^{(m)}|\leq(\frac{3}{10}+\frac{c}{5^i})|b_i^{(m)}-b_{i-1}^{(m)}|+\frac{c}{5^i}\|u_m\|_\infty.$$
Hence
$$|b_{i+1}^{(m)}-b_{i}^{(m)}|\leq\frac{3}{10}|b_i^{(m)}-b_{i-1}^{(m)}|+\frac{c}{5^i}\|u_m\|_\infty.$$
Similarly we have the estimates
$$|b_{i}^{(m)}-b_{i-1}^{(m)}|\leq\frac{3}{10}|b_{i-1}^{(m)}-b_{i-2}^{(m)}|+\frac{c}{5^{i-1}}\|u_m\|_\infty$$
till
$$|b_{m_1+2}^{(m)}-b_{m_1+1}^{(m)}|\leq\frac{3}{10}|b_{m_1+1}^{(m)}-b_{m_1}^{(m)}|+\frac{c}{5^{m_1+1}}\|u_m\|_\infty.$$
A routine argument shows that
$$|b_{i+1}^{(m)}-b_{i}^{(m)}|\leq(\frac{3}{10})^{i-m_1}|b_{m_1+1}^{(m)}-b_{m_1}^{(m)}|+(\frac{3}{10})^{i-m_1-1}\frac{c}{5^{m_1+1}}\|u_m\|_\infty.$$
Hence we have proved that
$$|b_{i+1}^{(m)}-b_{i}^{(m)}|\leq
c(\frac{3}{10})^{i}\|u_m\|_\infty$$ where $c$ depends only on $m_1$.

Similarly, we have
$$|b_{i+2}^{(m)}-b_{i+1}^{(m)}|\leq
c(\frac{3}{10})^{i+1}\|u_m\|_\infty,$$ till
$$|b_{i+p}^{(m)}-b_{i+p-1}^{(m)}|\leq
c(\frac{3}{10})^{i+p-1}\|u_m\|_\infty.$$ By adding up the above
estimates, we finally get $|b_{i+p}^{(m)}-b_{i}^{(m)}|\leq
C_1(\frac{3}{10})^{i}\|u_m\|_\infty. \Box$

\textbf{Lemma 5.6.} \emph{For each $m\geq m_1$,
$\lim_{i\rightarrow\infty} b_i^{(m)}$ exists. Moreover, there exists
a constant $C_2>0$ depending only on $m_1$, such that
$|\lim_{i\rightarrow\infty}b_i^{(m)}|\leq
C_2(\frac{3}{10})^m\|u_m\|_\infty$ uniformly on $m\geq m_1$.}

\emph{Proof.} For each $m\geq m_1$, Lemma 5.5 says that each
sequence $\{b_i^{(m)}\}_{i\geq 1}$ is a Cauchy sequence, hence
$\lim_{i\rightarrow\infty}b_i^{(m)}$ exists.

Taking $i=m$, $p=1$ in Lemma 5.5, noticing that $b_m^{(m)}=0$, we
get that $|b_{m+1}^{(m)}|\leq C_1(\frac{3}{10})^m\|u_m\|_\infty.$

On the other hand, $\forall i>m+1$, notice that $|b_i^{(m)}|\leq
|b_i^{(m)}-b_{m+1}^{(m)}|+|b_{m+1}^{(m)}|$. By using Lemma 5.5
again, we have
$$|b_i^{(m)}|\leq C_1(\frac{3}{10})^{m+1}\|u_m\|_\infty+
C_1(\frac{3}{10})^m\|u_m\|_\infty=C_2(\frac{3}{10})^m\|u_m\|_\infty.$$

Letting $i\rightarrow \infty$, we get the desired result. $\Box$

 In the following context,  for each $m\geq m_1$,  let $\theta_m$ denote the limit
$\lim_{i\rightarrow\infty} b_i^{(m)}/\|u_m\|_\infty$. Lemma 5.6
guarantees the existence of this limit, and furthermore,
$|\theta_m|\leq C_2(\frac{3}{10})^m$. Let $v_m:=\frac{u_m}{\|
u_m\|_\infty}$.  Then $v_m$ on $\Omega$ satisfies that
\begin{eqnarray*}\left\{
\begin{array}{l}
  -\Delta v_m=5^m\Phi(\lambda_m) v_m \mbox{ on } \Omega,\\
  v_m(q_0)=0,\\
  v_m|_L=\theta_m.
\end{array}\right.
\end{eqnarray*}

We will prove that $\{v_m\}_{m\geq m_1}$ contains a subsequence
converging  uniformly to a continuous function  on $\Omega$, which
is a Dirichlet eigenfunction associated to $\lambda$.

\textbf{Lemma 5.7.}\emph{ $\{\partial_n v_m(q_0)\}_{m\geq m_1}$ is
uniformly bounded, i.e., there exist a constant $C_3>0$ depending
only on $m_1$, such that $|\partial_n v_m(q_0)|\leq C_3$.}

\emph{Proof.} Let $m\geq m_1$. Choosing a harmonic function $h$ such
that $h(q_0)=1$, $h(F_1q_0)=h(F_2q_0)=0$, the local Gauss-Green
formula on $F_0(\mathcal{SG})$ says that
$$\mathcal{E}_{F_0(\mathcal{SG})}(v_m,h)=\int_{F_0(\mathcal{SG})}(-\Delta v_m)hd\mu+\sum_{\partial F_0(\mathcal{SG})}h\partial_n v_m.$$
Hence $|\partial_n v_m(q_0)|\leq
|\mathcal{E}_{F_0(\mathcal{SG})}(v_m,h)|+|\int_{F_0(\mathcal{SG})}(-\Delta
v_m)hd\mu|.$

 Since $h$ is harmonic on $F_0(\mathcal{SG})$, we have
$\mathcal{E}_{F_0(\mathcal{SG})}(v_m,h)=\frac{5}{3}\mathcal{E}(v_m\circ
F_0, h\circ F_0)=\frac{5}{3}\mathcal{E}_0(v_m\circ F_0, h\circ
F_0).$ Noticing that $h(q_0)=1$, $h(F_1q_0)=h(F_2q_0)=0$, we get
$|\mathcal{E}_{F_0(\mathcal{SG})}(v_m,h)|\leq c_1$, since
$\|v_m\|_\infty=1$.

On the other hand, since $-\Delta v_m=5^m\Phi(\lambda_m)v_m$, we
have $|\int_{F_0(\mathcal{SG})}(-\Delta v_m)hd\mu|\leq
5^m\Phi(\lambda_m)\|v_m\|_\infty\cdot\|h\|_\infty\mu(F_0(\mathcal{SG}))\leq
c_2$, since $5^m\Phi(\lambda_m)\rightarrow \lambda$.

Hence $|\partial_n v_m(q_0)|\leq c_1+c_2\triangleq C_3$. $\Box$

\textbf{Lemma 5.8.} \emph{$\{\mathcal{E}(v_m)\}_{m\geq m_1}$ is
uniformly bounded, i.e., there exists a constant $C_4>0$ depending
only on $m_1$, such that $\mathcal{E}(v_m)\leq C_4$.}

\emph{Proof.} $\forall n\geq m_1$, let $K_n$ be the part of $\Omega$
above $\partial\Omega_{n}\setminus\{q_0\}$. We first prove
$\{\mathcal{E}_{K_n}(v_m)\}_{m\geq m_1}$ is uniformly bounded and
the upper bound is independent on $n$.

Fix $n\geq m_1$, $m\geq m_1$. The Gauss-Green formula says that
$\int_{K_n}\Delta v_md\mu=\sum_{\partial K_n}\partial_n v_m$. From
the symmetry property of $v_m$, $\partial_n v_m$ takes same value
along $\partial K_n\setminus\{q_0\}$. Hence we get
\begin{equation}\label{514}-5^m\Phi(\lambda_m)\int_{K_n}v_md\mu=\partial_n
v_m(q_0)+2^n\partial_nv_m(F_1^n(q_0)).\end{equation}

On the other hand, the Gauss-Green formula also says that
\begin{eqnarray*}
\mathcal{E}_{K_n}(v_m)&=&\int_{K_n}(-\Delta
v_m)v_md\mu+\sum_{\partial K_n}v_m\partial_nv_m\\
&=&5^m\Phi(\lambda_m)\int_{K_n}v_m^2d\mu+2^nv_m(F_1^nq_0)\partial_nv_m(F_1^nq_0),
\end{eqnarray*}
since $v_m(q_0)=0$. Combined with $(\ref{514})$, it follows
$$\mathcal{E}_{K_n}(v_m)=5^m\Phi(\lambda_m)\int_{K_n}v_m^2d\mu+v_m(F_1^nq_0)(-5^m\Phi(\lambda_m)\int_{K_n}v_md\mu-\partial_n v_m(q_0)).$$

Since $5^m\Phi(\lambda_m)\rightarrow\lambda$, there exists a
constant $c>0$, such that $5^m\Phi(\lambda_m)\leq c$. Hence
$$\mathcal{E}_{K_n}(v_m)\leq c\|v_m\|^2_\infty+\|v_m\|_\infty(c\|v_m\|_\infty+|\partial_n v_m(q_0)|).$$
Using Lemma 5.7, we get $\mathcal{E}_{K_n}(v_m)\leq 2c+C_3\triangleq
C_4$. Since the above inequality is independent on $n$, we then get
the desired result by passing $n$ to infinity. $\Box$

Now we come to the main purpose of this section.

\emph{Proof of Theorem 3.6.}

 For each $m\geq m_1$, since $v_m\in \mathcal{F}$, we
have $$|v_m(x)-v_m(y)|\leq \mathcal{E}(v_m)^{1/2}
d(x,y)^{1/2},\quad\forall x,y\in \Omega,$$ where $d(\cdot,\cdot)$ is
the effective resistance metric on $\Omega$. Hence by Lemma 5.8,
$$|v_m(x)-v_m(y)|\leq C^{1/2}_4d(x,y)^{1/2}, \quad \forall x,y\in
\Omega$$ holds uniformly on $m\geq m_1$. Thus $\{v_m\}_{m\geq m_1}$
is equicontinuous. Moreover, notice that $\{v_m\}_{m\geq m_1}$ is
also uniformly bounded. Then  using the Arzel\`{a}-Ascoli theorem,
there exists a subsequence $\{v_{m_k}\}$ of $\{v_m\}$ which
converges uniformly to a continuous function $v$ on $\Omega$.

Let $G_\Omega(x,y)$ denote the Green's function associated to
$\Omega$. See the explicit expression for $G_\Omega(x,y)$ in
$\cite{Guo}$. Then $\forall k$, we have
\begin{equation}\label{515}
v_{m_k}(x)=\int_\Omega
G_\Omega(x,y)5^{m_k}\Phi(\lambda_{m_k})v_{m_k}(y)d\mu(y)+h_{m_k}(x),
\end{equation}
where $h_{m_k}$ is a harmonic function on $\Omega$ taking the same
boundary values as $v_{m_k}$. Namely, $h_{m_k}(q_0)=0$ and
$h_{m_k}|_L=\theta_{m_k}$. If $k\rightarrow\infty$, then
$\theta_{m_k}\rightarrow 0$ and hence $h_{m_k}$ goes to $0$
uniformly on $\Omega$ by the maximum principle. Hence by letting
$k\rightarrow\infty$ on both side of $(\ref{515})$, we get
$$v(x)=\int_{\Omega}G_\Omega(x,y)(\lambda v(y))d\mu(y).$$
Thus we finally get
\begin{eqnarray*}\left\{
\begin{array}{l}
  -\Delta v=\lambda v \mbox{ in } \Omega,\\
  v|_{\partial\Omega}=0.
\end{array}\right.
\end{eqnarray*}
Hence $v$ is a Dirichlet eigenfunction associated to $\lambda$.
$\Box$

Thus for each sequence $\{\lambda_m\}_{m\geq m_0}$, we have proved
that $\lambda=\frac{3}{2}\lim_{m\rightarrow}5^m\lambda_m$ is a
symmetric primitive Dirichlet  eigenvalue of $-\Delta$ on $\Omega$.
We denote by $\mathcal{P}_*^{+}$ the totality of all this kind of
eigenvalues. Of cause, $\mathcal{P}_*^{+}\subset \mathcal{P}^{+}$.
Lemma 5.3 and Lemma 5.4 guarantee that all eigenvalues in
$\mathcal{P}_*^{+}$ are distinct and they are all greater than $0$.
In the next section, we will prove that all eigenvalues in
$\mathcal{P}^{+}$ come in this way. Namely,
$\mathcal{P}_*^{+}=\mathcal{P}^{+}$.

The skew-symmetric case is similar. We denote by $\mathcal{P}_*^{-}$
the set of skew-symmetric eigenvalues generated in this way. Let
$\mathcal{P}_*=\mathcal{P}_*^+\cup\mathcal{P}_*^-$ denote all the
associated primitive eigenvalues. Let $\mathcal{M}_*$ be the set of
miniaturized eigenvalues generated by $\mathcal{P}_*^-$.
Accordingly, $\mathcal{P}_*^{-}\subset \mathcal{P}^{-}$,
$\mathcal{P}_*\subset \mathcal{P}$ and $\mathcal{M}_*\subset
\mathcal{M}$.

\section{Complete Dirichlet spectrum of $-\Delta$}

It is clear that the weak spectral decimation recipe constructs many
primitive eigenvalues (hence also many miniaturized eigenvalues) of
$-\Delta$ on $\Omega$. Recall that the standard spectral decimation
recipe also constructs many localized eigenvalues of $-\Delta$ on
$\Omega$. It is natural to ask do these recipes construct the whole
Dirichlet spectrum $\mathcal{S}$? In this section, we will give an
affirmative  answer of this question.

Till now, for each $m\geq 2$, we have proved that the Dirichlet
spectrum $\mathcal{S}_m$ of the discrete Laplacian $-\Delta_m$ on
$\Omega_m$ consists of $\mathcal{L}_m$, $\mathcal{P}_m$ and
$\mathcal{M}_m$ the three types of eigenvalues. After passing the
approximations to the limit, we have proved that there are at least
three types of eigenvalues $\mathcal{L}$, $\mathcal{P}_*$ and
$\mathcal{M}_*$ in the Dirichlet spectrum $\mathcal{S}$ of
$-\Delta$, which could generated by the (weak) spectral decimation
recipe. Namely, $\mathcal{S}\supset \mathcal{L}\cup
\mathcal{P}_*\cup \mathcal{M}_*$. We call all of the above three
types of eigenvalues \emph{raw eigenvalues}. By the \emph{raw
multiplicity} of the raw eigenvalue $\lambda$, we mean the
multiplicity of the associated eigenvalue $\lambda_{m_0}$ of
$-\Delta_{m_0}$, where $m_0$ is the generation of birth. Since
linearly independent eigenfunctions of $-\Delta_{m_0}$ belonging to
$\lambda_{m_0}$ give rise to linearly independent eigenfunctions of
$-\Delta$, and the fact that all primitive graph eigenvalues have
only raw multiplicity $1$, the raw multiplicity of $\lambda$ is not
greater than the true multiplicity of $\lambda$.

Denote by $\mathcal{S}_*$ the collection of raw eigenvalues of
$-\Delta$, then $\mathcal{S}_*=\mathcal{L}\cup \mathcal{P}_*\cup
\mathcal{M}_*$ and $\mathcal{S}_*\subset \mathcal{S}$. Hence we need
to prove $\mathcal{S}_*=\mathcal{S}$, $ \mathcal{P}_*=\mathcal{P}$
and $\mathcal{M}_*=\mathcal{M}$ and the raw multiplicity of each
element of $\mathcal{S}_*$ coincides with its true multiplicity.

Comparing to the proof of the analogous problem for the standard
$\mathcal{SG}\setminus{V_0}$ case (see details in $\cite{Fu}$), to
prove the above results, the following proposition will play a vital
role.  Recall that
$a_m=\sharp(V_m^\Omega\setminus\partial\Omega_m)=\frac{3^{m+1}-1}{2}-2^{m+1}$.

\textbf{Proposition 6.1.} \emph{Let $0<\kappa_1\leq
\kappa_2\leq\cdots$ be the rearrangement of elements of
$\mathcal{S}_*$ each repeated according to its raw multiplicity. Let
$\{\kappa_{m,i}\}_{1\leq i\leq a_m}$ be the $m$-level graph
eigenvalues of $-\Delta_m$ on $\Omega_m$ including multiplicities.
Then
$$\lim_{m\rightarrow\infty}\sum_{1\leq i\leq a_m}\frac{1}{\frac{3}{2}5^m\kappa_{m,i}}=\sum_{i=1}^{\infty}\frac{1}{\kappa_i}<\infty.$$}

In order to prove this proposition, we first list some notations and
lemmas. It is more convenient to consider the following slightly
different classification of all the raw eigenvalues of $-\Delta$,
$$\mathcal{S}_*=\mathcal{L}\cup \mathcal{P}_*^+\cup
\widetilde{\mathcal{P}}_*^-$$ where
$\widetilde{\mathcal{P}}_*^-=\mathcal{P}_*^-\cup \mathcal{M}_*$,
since miniaturized eigenvalues have the same generation mechanism as
the skew-symmetric primitive eigenvalues. In the following, we
always use $\alpha, \beta,\gamma$ to denote $\mathcal{L}$,
$\mathcal{P}_*^+$, $\widetilde{\mathcal{P}}_*^-$ type eigenvalues
respectively. Accordingly, $\forall m\geq 2$, all the $m$-level
graph eigenvalues are classified into  the three types
$\mathcal{L}_m$, $\mathcal{P}_m^+$ and
$\widetilde{\mathcal{P}}_m^-$, where
$\widetilde{\mathcal{P}}_m^-=\mathcal{P}_m^-\cup \mathcal{M}_m$, and
we always use $\alpha_m, \beta_m,\gamma_m$ to denote eigenvalues in
them respectively. For simplicity, we denote
$A_m=\sharp\mathcal{L}_m$, $B_m=\sharp\mathcal{P}_m^+$ and
$C_m=\sharp\widetilde{\mathcal{P}}_m^-$. Of course,
$a_m=A_m+B_m+C_m$. Moreover, recall that $\rho_m^\Omega(5)$ and
$\rho_m^\Omega(6)$ are the multiplicities of $m$-level initial
eigenvalues $5$ and $6$ respectively.  See the exact values of them
in Section 3.

\textbf{Lemma 6.1.} \emph{$\mathcal{L}=\bigcup_{k=3}^\infty
\mathcal{L}^{k}$  (disjoint union) where $\mathcal{L}^{k}\subset
[5^k\Phi(3), 5^k\Phi(5)]$. }

 \emph{Proof.} $\forall \alpha\in \mathcal{L}$, let $\{\alpha_m\}_{m\geq m_0}$ be the corresponding sequence of eigenvalues with a generation of fixation $m_1$.
Then
$\alpha=\frac{3}{2}\lim_{m\rightarrow\infty}5^m\alpha_m=5^{m_1}\Phi(\alpha_{m_1}).$

If $\alpha_{m_1}$ is an initial eigenvalue, then $\alpha_{m_1}$ can
only be equal to $5$. If $\alpha_{m_1}$ is a continued eigenvalue,
then $\alpha_{m_1}=\phi_+(\alpha_{m_1-1})$, which yields that $3\leq
\alpha_{m_1}\leq 5$. Hence we always have $3\leq \alpha_{m_1}\leq
5$.

Noticing that each localized eigenvalue has generation of birth at
least $3$, denote by $\mathcal{L}^{k}$ the set of eigenvalues with
$m_1=k$, $k=3,4,\cdots$. Then $\mathcal{L}=\bigcup_{k=3}^\infty
\mathcal{L}^{k}$  and $\mathcal{L}^{k}\subset [5^k\Phi(3),
5^k\Phi(5)]$. Since $\phi_-(5)<3$, we have $\Phi(5)<5\Phi(3)$. Hence
 $\mathcal{L}=\bigcup_{k=3}^\infty
\mathcal{L}^{k}$ is a disjoint union.  $\Box$

\textbf{Lemma 6.2.} \emph{$\mathcal{P}_*^+=\bigcup_{k=2}^\infty
\mathcal{P}_*^{+,k}$  (disjoint union) where
$\mathcal{P}_*^{+,2}\subset (0,5^2\Phi(6)]$ and
$\mathcal{P}_*^{+,k}\subset [5^k\Phi(\phi_-(3)), 5^k\Phi(6)]$ for
$k\geq 3$. }

\emph{Proof.} $\forall\beta\in \mathcal{P}_*^+$, let
$\{\beta_m\}_{m\geq m_0}$ be the corresponding sequence of
eigenvalues with a generation of fixation $m_1$. Then
$\beta=\frac{3}{2}\lim_{m\rightarrow\infty}5^m\beta_m=5^{m_1}\lim_{n\rightarrow\infty}\frac{3}{2}5^n\widetilde{\phi}_-^{(n)}(\beta_{m_1}).$

If $\beta_{m_1}$ is a continued eigenvalue (hence $m_1\geq 3$), then
we must have $\beta_{m_1}=\widetilde{\phi}_+(\beta_{m_1-1})$, which
obviously yields that
$\beta_{m_1}>\widetilde{\phi}_-(\beta_{m_1-1}^*)$ where
$\beta_{m_1-1}^*$ denotes the largest  eigenvalue in
$\mathcal{P}_{m_1-1}^+$. If $\beta_{m_1}$ is an initial eigenvalue
with $m_1\geq 3$, then obviously
$\beta_{m_1}>\widetilde{\phi}_-(\beta_{m_1-1}^*)$. Hence we always
have $\beta_{m_1}>\widetilde{\phi}_-(\beta_{m_1-1}^*)$ if $m_1\geq
3$.

Moreover, When $m_1>3$, if we denote $\beta_{m_1-1}^{**}$ the
largest eigenvalue in $\mathcal{P}_{m_1-1}^+$ except for
$\beta_{m_1-1}^{*}$, then we have
$\widetilde{\phi}_-(\beta_{m_1-1}^*)>\phi_-(\beta_{m_1-1}^{**})$. It
is easy to check that $\beta_{m_1-1}^{**}>\phi_+(\beta_{m_1-2}^*)>3$
since $m_1>3$. Thus
$\beta_{m_1}>\widetilde{\phi}_-(\beta_{m_1-1}^*)>\phi_-(3)$. When
$m_1=3$, it can be checked directly that
$\beta_3>\widetilde{\phi}_-(\beta_2^*)\approx1.33>\phi_-(3)$. Hence
we always have
$\beta_{m_1}>\widetilde{\phi}_-(\beta_{m_1-1}^*)>\phi_-(3)$ if
$m_1\geq 3$. By Lemma 5.2, we have
$\widetilde{\phi}_-^{(n)}(\beta_{m_1})>\phi_-^{(n)}(\widetilde{\phi}_-(\beta_{m_1-1}^*)),\forall
n\in \mathbb{N}.$ Hence if $m_1\geq 3$, we have
\begin{eqnarray*}
\beta&=&5^{m_1}\lim_{n\rightarrow\infty}\frac{3}{2}5^n\widetilde{\phi}_-^{(n)}(\beta_{m_1})\\
&\geq&
5^{m_1}\lim_{n\rightarrow\infty}\frac{3}{2}5^n\phi_-^{(n)}(\widetilde{\phi}_-(\beta_{m_1-1}^*))\\
&\geq&5^{m_1}\lim_{n\rightarrow\infty}\frac{3}{2}5^n\phi_-^{(n)}(\phi_-(3))\\
&=& 5^{m_1}\Phi(\phi_-(3)).
\end{eqnarray*}

On the other hand, when $m_1\geq 2$, we always have
$$\beta=5^{m_1}\lim_{n\rightarrow\infty}\frac{3}{2}5^n\widetilde{\phi}_-^{(n)}(\beta_{m_1})\leq
5^{m_1}\lim_{n\rightarrow\infty}\frac{3}{2}5^n\phi_-^{(n)}(6)=5^{m_1}\Phi(6).$$

Denote by $\mathcal{P}_*^{+,k}$ the set of  eigenvalues with
$m_1=k$, $k=2,3,\cdots$. Then
$\mathcal{P}_*^{+}=\bigcup_{k=2}^\infty\mathcal{P}_*^{+,k}$ where
$\mathcal{P}_*^{+,2}\subset (0,5^2\Phi(6)]$ and
$\mathcal{P}_*^{+,k}\subset[5^k\Phi(\phi_-(3)),5^k\Phi(6)]$ for
$k\geq 3$.

Next we need to prove
$\mathcal{P}_*^{+}=\bigcup_{k=2}^\infty\mathcal{P}_*^{+,k}$ is a
disjoint union. $\forall 2\leq k<k'$, take an element $\beta$ in
$\mathcal{P}_*^{+,k}$,  $\beta'$ in $\mathcal{P}_*^{+,k'}$
respectively. Then
$\beta=5^{k}\lim_{n\rightarrow\infty}\frac{3}{2}5^n\widetilde{\phi}_-^{(n)}(\beta_{k})$
for some eigenvalue $\beta_{k}$ in $\mathcal{P}_{k}^+$, and
$\beta'=5^{k'}\lim_{n\rightarrow\infty}\frac{3}{2}5^n\widetilde{\phi}_-^{(n)}(\beta_{k'}')$
for some eigenvalue $\beta_{k'}'$ in $\mathcal{P}_{k'}^+$.

Note that $\widetilde{\phi}_-^{(k'-k)}(\beta_{k})$ and $\beta_{k'}'$
both belong to $\mathcal{P}_{k'}^+$. Since $k'$ is the generation of
fixation of $\beta'$, we can easily get
$\widetilde{\phi}_-^{(k'-k)}(\beta_{k})<\beta_{k'}'$. Then by using
Lemma 5.3, we have $\beta<\beta'$.

From the arbitrariness of $\beta$, $\beta'$ and $k$, $k'$, we
finally get that
$\mathcal{P}_*^{+}=\bigcup_{k=2}^\infty\mathcal{P}_*^{+,k}$ is a
disjoint union. $\Box$

\textbf{Lemma 6.3.}  \emph{Let $0<\alpha_1\leq \alpha_2\leq\cdots$
be the rearrangement of elements of $\mathcal{L}$ each repeated
according to its raw multiplicity. Let $\{\alpha_{m,i}\}_{1\leq
i\leq A_m}$ be the $m$-level localized eigenvalues of $-\Delta_m$ on
$\Omega_m$ including multiplicities. Then
$$\lim_{m\rightarrow\infty}\sum_{1\leq i\leq A_m}\frac{1}{\frac{3}{2}5^m\alpha_{m,i}}=\sum_{i=1}^{\infty}\frac{1}{\alpha_i},$$}
providing $\sum_{i=1}^{\infty}\frac{1}{\alpha_i}<\infty$.

\emph{Proof.} Noticing that
$\lim_{m\rightarrow\infty}\frac{\rho_m^\Omega(6)}{5^m}=0$, it
suffices to show that
\begin{equation}\label{61}
\sum_{1\leq i\leq A_m\atop \alpha_{m,i}\neq
6}\frac{1}{\frac{3}{2}5^m\alpha_{m,i}}-\sum_{i=1}^{A_m-\rho_{m}^\Omega(6)}\frac{1}{\alpha_i},
\end{equation}
converges to $0$ as $m$ goes to infinity.

$\forall m\geq 2$, denote $D_m=A_m-\rho_m^\Omega(6)$. By Lemma 6.1,
$\{\alpha_{1},\alpha_2,\cdots,\alpha_{D_m}\}$ is an arrangement of
elements of $\bigcup_{k=3}^m\mathcal{L}^k$ each being repeated
according to its raw multiplicity. The first sum of $(\ref{61})$ has
also $D_m$ terms, which can be rearranged so that
$$\lim_{n\rightarrow\infty}\frac{3}{2}5^{m+n}\phi_-^{(n)}(\alpha_{m,i})=\alpha_i,\quad 1\leq i\leq D_m.$$

Hence by using Lemma 6.1, $(\ref{61})$ is  equal to
$\sum_{k=3}^m\sum_{\alpha_i\in
\mathcal{L}^k}(\frac{1}{\frac{3}{2}5^m\alpha_{m,i}}-\frac{1}{\alpha_i})$.
If $\alpha_{i}\in \mathcal{L}^k$ ($k=3,\cdots,m$), then
$\alpha_{i}=5^k\Phi(\theta)$ for some $\theta\in [3,5]$ and
accordingly the corresponding $\alpha_{m,i}$ is of the form
$\alpha_{m,i}=\phi_-^{(m-k)}(\theta)$. Hence
$$0<\frac{1}{\frac{3}{2}5^m\alpha_{m,i}}-\frac{1}{\alpha_i}=\frac{1}{5^k}(\frac{1}{\frac{3}{2}5^{m-k}\phi_-^{(m-k)}(\theta)}-\frac{1}{\Phi(\theta)}).$$
Since $\frac{1}{\frac{3}{2}5^n\phi_-^{(n)}(x)}$ converges to
$\frac{1}{\Phi(x)}$ uniformly on $[3,5]$ as $n$ goes to infinity,
$\forall\varepsilon>0$, the last expression is dominated by
$\frac{\varepsilon}{5^k}$ whenever $m-k$ is greater than some number
$N$. When $m-k\leq N$, the same expression is dominated by
$\frac{1}{5^mR}$ for $R=\frac{3}{2}\inf_{3\leq x\leq
5}\phi_-^{(N)}(x)$. The number of $\alpha_i$'s in $\mathcal{L}^k$ is
less than $A_{k-1}+\rho_k^\Omega(5)$, so $(\ref{61})$ is dominated
by
$$\sum_{k=3}^{m-N-1}\frac{A_{k-1}+\rho_{k}^\Omega(5)}{5^k}\varepsilon+\sum_{k=m-N}^m\frac{A_{k-1}+\rho_{k}^\Omega(5)}{5^mR}\leq c_1\varepsilon+c_2(\frac{3}{5})^m\frac{1}{R}$$
for some constants $c_1,c_2>0$. Then let $m$ be large enough,
$(\ref{61})$ can be dominated by $(c_1+c_2)\varepsilon$. Hence we
have proved $\sum_{1\leq i\leq A_m\atop \alpha_{m,i}\neq
6}\frac{1}{\frac{3}{2}5^m\alpha_{m,i}}-\sum_{i=1}^{D_m}\frac{1}{\alpha_i}$
converges to $0$ as $m$ goes to infinity.
 $\Box$

\textbf{Lemma 6.4.} \emph{Let $0<\beta_1< \beta_2<\cdots$ be the
elements of $\mathcal{P}_*^+$ in increasing order. Let
$\{\beta_{m,i}\}_{1\leq i\leq B_m}$ be the $m$-level symmetric
primitive eigenvalues of $-\Delta_m$ on $\Omega_m$. Then
$$\lim_{m\rightarrow\infty}\sum_{1\leq i\leq B_m}\frac{1}{\frac{3}{2}5^m\beta_{m,i}}=\sum_{i=1}^{\infty}\frac{1}{\beta_i},$$}
providing $\sum_{i=1}^{\infty}\frac{1}{\beta_i}<\infty$.

\emph{Proof.} It suffices to prove that
\begin{equation}\label{62}
\sum_{1\leq i\leq
B_m}\frac{1}{\frac{3}{2}5^m\beta_{m,i}}-\sum_{i=1}^{B_m}\frac{1}{\beta_i},
\end{equation}
converges to $0$ as $m$ goes to infinity.

By Lemma 6.2, $\{\beta_{1},\beta_2,\cdots,\beta_{B_m}\}$ is an
arrangement of elements of $\bigcup_{k=2}^m\mathcal{P}_*^{+,k}$. The
first sum of $(\ref{62})$  can be rearranged so that
$$\lim_{n\rightarrow\infty}\frac{3}{2}5^{m+n}\widetilde{\phi}_-^{(n)}(\beta_{m,i})=\beta_i,\quad 1\leq i\leq B_m.$$

Hence by using Lemma 6.2, $(\ref{62})$ is  equal to
$\sum_{k=2}^m\sum_{\beta_i\in
\mathcal{P}_*^{+,k}}(\frac{1}{\frac{3}{2}5^m\beta_{m,i}}-\frac{1}{\beta_i}).$
The $k=2$ term converges to $0$ as $m$ goes to infinity since
$\sharp\mathcal{P}_*^{+,2}=B_2=3$.

Hence we only need to prove
\begin{equation}\label{63}\sum_{k=3}^m\sum_{\beta_i\in
\mathcal{P}_*^{+,k}}|\frac{1}{\frac{3}{2}5^m\beta_{m,i}}-\frac{1}{\beta_i}|.
\end{equation}
converges to $0$ as $m$ goes to infinity.

If $\beta_{i}\in \mathcal{P}_*^{+,k}$ ($k=3,\cdots,m$), then
$\beta_{i}=5^k\lim_{n\rightarrow\infty}\frac{3}{2}5^n\widetilde{\phi}_-^{(n)}(\theta)$
for some $\theta\in \mathcal{P}_k^+$ and accordingly the
corresponding $\beta_{m,i}$ is of the form
$\beta_{m,i}=\widetilde{\phi}_-^{(m-k)}(\theta)$. Hence
\begin{equation}\label{64}
|\frac{1}{\frac{3}{2}5^m\beta_{m,i}}-\frac{1}{\beta_i}|=\frac{1}{5^k}|\frac{1}{\frac{3}{2}5^{m-k}\widetilde{\phi}_-^{(m-k)}(\theta)}-\frac{1}{\lim_{n\rightarrow\infty}\frac{3}{2}5^n\widetilde{\phi}_-^{(n)}(\theta)}|.
\end{equation}

From the proof of Lemma 6.2, we have
$$\frac{3}{2}5^n\phi_-^{(n)}(\phi_-(3))<\frac{3}{2}5^n\widetilde{\phi}_-^{(n)}(\theta)<\frac{3}{2}5^n\phi_-^{(n)}(6).$$
Then by the proof of Lemma 5.1, $\forall\varepsilon>0$, the right
side of formula $(\ref{64})$ is dominated by
$\frac{1}{5^k}\varepsilon$ whenever $m-k$ is greater than some
number $N$. When $m-k\leq N$,
$\frac{1}{\frac{3}{2}5^{m}\widetilde{\phi}_-^{(m-k)}(\theta)}$ is
dominated by $\frac{1}{5^m R}$ for $R=\frac{3}{2}\phi_-^{(N+1)}(3)$.
The number of $\beta_i$'s in $\mathcal{P}_*^{+,k}$ is controlled by
$B_k$, so the sum $(\ref{63})$ is dominated by
\begin{eqnarray}\label{65}
&&\sum_{k=3}^{m-N-1}\frac{B_k}{5^k}\varepsilon+\sum_{k=m-N}^m\frac{B_k}{5^mR}+\sum_{k=m-N}^m\sum_{\beta_i\in\mathcal{P}_*^{+,k}}\frac{1}{\beta_i}\nonumber\\
\leq&&
c_1\varepsilon+c_2(\frac{2}{5})^m\frac{1}{R}+\sum_{k=m-N}^m\sum_{\beta_i\in\mathcal{P}_*^{+,k}}\frac{1}{\beta_i}.
\end{eqnarray}
Noticing that $\sum_{i=1}^\infty\frac{1}{\beta_{i}}<\infty$, the
last term goes to $0$ as $m$ goes to infinity. Hence for large $m$,
$(\ref{65})$ is less than $(c_1+c_2+1)\varepsilon$. Thus we have
proved $\sum_{1\leq i\leq
B_m}\frac{1}{\frac{3}{2}5^m\beta_{m,i}}-\sum_{i=1}^{B_m}\frac{1}{\beta_i}$
converges to $0$ as $m$ goes to infinity.  $\Box$

\textbf{Lemma 6.5.} \emph{Let $0<\gamma_1\leq \gamma_2\leq\cdots$ be
the elements of $\widetilde{\mathcal{P}}_*^-$ in increasing order
repeated according to their raw multiplicities. Let
$\{\gamma_{m,i}\}_{1\leq i\leq C_m}$ be the $m$-level
$\widetilde{\mathcal{P}}_m^-$ type eigenvalues of $-\Delta_m$ on
$\Omega_m$ including multiplicities. Then
$$\lim_{m\rightarrow\infty}\sum_{1\leq i\leq C_m}\frac{1}{\frac{3}{2}5^m\gamma_{m,i}}=\sum_{i=1}^{\infty}\frac{1}{\gamma_i},$$}
providing $\sum_{i=1}^{\infty}\frac{1}{\gamma_i}<\infty$.

The proof is similar to those of Lemma 6.3 and Lemma 6.4.

\emph{Proof of Proposition 6.1.}

Let $0<\widetilde{\kappa}_1\leq \widetilde{\kappa}_2\leq\cdots$ be
the rearrangement of elements of $\mathcal{S}$ each repeated
according to its true multiplicity. Let $\widetilde{v}_1,
\widetilde{v}_2,\cdots$ be the associated eigenfunctions. Let
$G_\Omega(x,y)$ be the Green's function for $\Omega$. Then
$G_\Omega(x,y)$ can be expanded as a uniformly convergence series
$$G_\Omega(x,y)=\sum_{i=1}^\infty\frac{\widetilde{v}_i(x)\widetilde{v}_i(y)}{\widetilde{\kappa}_i},\quad\forall x,y\in \Omega.$$

Since $\mathcal{S}_*\subset \mathcal{S}$ and the raw multiplicity is
not greater than the true one, we get that
$\sum_{i=1}^{\infty}\frac{1}{\kappa_i}<\infty$. Hence
$\sum_{i=1}^{\infty}\frac{1}{\alpha_i}<\infty$,
$\sum_{i=1}^{\infty}\frac{1}{\beta_i}<\infty$, and
$\sum_{i=1}^{\infty}\frac{1}{\gamma_i}<\infty$. The by adding up the
results in Lemma 6.3, Lemma 6.4 and Lemma 6.5, we have
$$\lim_{m\rightarrow\infty}\sum_{1\leq i\leq a_m}\frac{1}{\frac{3}{2}5^m\kappa_{m,i}}=\sum_{i=1}^{\infty}\frac{1}{\kappa_i}<\infty. \Box$$

\textit{Proof of Theorem 3.8.}

Based on Proposition 6.1, following similar argument in $\cite{Fu}$,
we finally get $\mathcal{S}_*=\mathcal{S}$, hence
$\mathcal{S}=\mathcal{L}\cup \mathcal{P}_*^{+}\cup
\widetilde{\mathcal{P}}_*^-$. As a immediate consequence, we have
$\mathcal{P}_*^{+}=\mathcal{P}^+$, $\mathcal{P}_*^{-}=\mathcal{P}^-$
and $\mathcal{M}_*=\mathcal{M}$. Thus we have
$\mathcal{S}=\mathcal{L}\cup \mathcal{P}\cup \mathcal{M}$ where the
union is disjoint. $\Box$

\textit{Proof of Theorem 3.7.}

The proof follows immediately from Lemma 5.3, Theorem 3.6 and
Theorem 3.8. $\Box$

The skew-symmetric analog of Theorem 3.7 is obvious. we omit it.

\textit{Proof of Theorem 3.1.}

We only need to prove $\mathcal{P}^+\cap \mathcal{P}^-=\emptyset.$

Let $\lambda$ belong to both $\mathcal{P}^+$ and $\mathcal{P}^-$.
Then there exist a symmetric primitive eigenfunction $u_1$ and  a
skew-symmetric primitive eigenfunction $u_2$ on $\Omega$, both
associated to $\lambda$. Let $m\geq 1$, consider the subdomain
$F_1^m(\Omega)$ of $\Omega$, which is a $m$-times contraction of
$\Omega$, with the boundary $F_1^mq_0$ and $F_1^m(L)$. From Lemma
4.12, the restriction of $u_1$, $u_2$ to $F_1^m(\Omega)$ should be
linear independent. Notice that on the bottom line segment
$F_1^m(L)$, both $u_1$ and $u_2$ satisfy the Dirichlet  condition.
Let $u$ be a linear composition of $u_1$ and $u_2$ such that $u$
satisfy the Dirichlet boundary condition on the vertex $F_1^mq_0$.
Then obviously we have that $-\Delta u=\lambda u $ on
$F_1^m(\Omega)$ satisfying the Dirichlet boundary condition on
$\partial F_1^m(\Omega)$.
 Hence $u\circ F_1^m$ becomes a
symmetric primitive Dirichlet  eigenfunction on $\Omega$ associated
to an eigenvalue $\frac{1}{5^m}\lambda$ by using the scaling
property of $\Delta$. Thus $\frac{1}{5^m}\lambda\in \mathcal{P}^+$.
However, from Lemma 5.3 and Lemma 5.4, we have a smallest positive
element in $\mathcal{P}^+$. Hence we get a contradiction by choosing
$m$ sufficiently large.
 $\Box$

Hence we have constructed the complete Dirichlet spectrum of
$-\Delta$ on $\Omega$, and the raw multiplicity of each element of
$\mathcal{S}$ coincides with its true multiplicity.

Finally we turn to the Weyl's eigenvalue asymptotics on $\Omega$. As
showed before, we use $\rho^0(x)$ and $\rho^\Omega(x)$ to denote the
Dirichlet eigenvalue counting function with respect to
$\mathcal{SG}\setminus{V_0}$ and $\Omega$, respectively.

\textit{Proof of Theorem 3.10.}

We divide $\rho^\Omega(x)$ into four parts $\rho^{\mathcal{L}}(x)$,
$\rho^{\mathcal{P}^+}(x)$, $\rho^{\mathcal{P}^-}(x)$ and
$\rho^{\mathcal{M}}(x)$ corresponding to different types of
eigenvalues. The exact definitions are:
$\rho^{\mathcal{L}}(x)=\sharp\{\lambda\in \mathcal{L}: \lambda\leq
x\},$ $\rho^{\mathcal{P}^+}(x)=\sharp\{\lambda\in \mathcal{P}^+:
\lambda\leq x\},$ $\rho^{\mathcal{P}^-}(x)=\sharp\{\lambda\in
\mathcal{P}^-: \lambda\leq x\}$ and
$\rho^{\mathcal{M}}(x)=\sharp\{\lambda\in \mathcal{M}: \lambda\leq
x\}.$
 Obviously,
$$\rho^\Omega(x)-\rho^{\mathcal{L}}(x)=\rho^{\mathcal{P}^+}(x)+\rho^{\mathcal{P}^-}(x)+\rho^{\mathcal{M}}(x).$$

For  $\rho^{\mathcal{P}^+}(x)$, denote $\beta_m^*$ the largest
eigenvalue in $\mathcal{P}_m^+$, and $\beta^{(m)}$ the eigenvalue in
$\mathcal{P}^+$ corresponding to the sequence
$\{\widetilde{\phi}_-^{(n)}(\beta_m^*)\}_{n\geq 0}$, i.e.,
$\beta^{(m)}=\lim_{n\rightarrow\infty}\frac{3}{2}5^{n+m}
\widetilde{\phi}_-^{(n)}(\beta_m^*)$. By using Lemma 5.2, it is easy
to check that
\begin{equation}\label{ppp}
c_15^m=\lim_{n\rightarrow\infty}\frac{3}{2}5^{n+m}
{\phi}_-^{(n)}(2)\leq\lim_{n\rightarrow\infty}\frac{3}{2}5^{n+m}
{\phi}_-^{(n)}(\beta_m^{**})\leq\beta^{(m)}\leq\lim_{n\rightarrow\infty}\frac{3}{2}5^{n+m}
{\phi}_-^{(n)}(6)=c_25^m,
\end{equation}
for appropriate constants $c_1,c_2>0$, where $\beta_{m}^{**}$ denote
the largest eigenvalue in $\mathcal{P}_m^+$ except $\beta_m^*$.

Notice that the bottom $r_m$ eigenvalues in $\mathcal{P}^+$ are
generated from eigenvalues in $\mathcal{P}_m^+$ by extending these
eigenvalues by choosing $\widetilde{\phi}_-$ relation for all
$m'>m$. Hence we get
$$\rho^{\mathcal{P}^+}(\beta^{(m)})=r_m, \quad\forall m\geq 2.$$
Using $(\ref{ppp})$, we get $\rho^{\mathcal{P}^+}(c_15^m)\leq r_m$,
and $\rho^{\mathcal{P}^+}(c_25^m)\geq r_m$.

Denote by $k_0$ the least number such that $5^{k_0}c_1\geq c_2$.
$\forall x\geq 25c_2$, choose a number $m$ such that $c_25^m\leq
x<c_25^{m+1}$. Then $c_25^m\leq x< c_15^{m+k_0+1}$. Hence
$$c_3 x^{\log2/\log5}\leq r_m\leq\rho^{\mathcal{P}^+}(c_25^m)\leq\rho^{\mathcal{P}^+}(x)\leq\rho^{\mathcal{P}^+}(c_15^{m+k_0+1})\leq r_{m+k_0+1}\leq c_4 x^{\log2/\log5},$$
for appropriate constants $c_3,c_4>0$. Thus we have proved that for
$x$ large enough,
$$c_3x^{\log2/\log5}\leq\rho^{\mathcal{P}^+}(x)\leq c_4x^{\log2/\log5}. $$

Similar argument yields that for $x$ large enough,
$$c_5x^{\log2/\log5}\leq\rho^{\mathcal{P}^-}(x)\leq c_6x^{\log2/\log5},
$$ for appropriate constants $c_5,c_6>0$.

Now we consider $\rho^\mathcal{M}(x)$. Notice that for each
$\lambda'\in\{\lambda\in \mathcal{M}: \lambda\leq x\}$, there exists
a $k\geq 1$, such that $\lambda'$ has multiplicity $2^k$ in
$\mathcal{M}$, and $\frac{1}{5^k}\lambda'\in \{\lambda\in
\mathcal{P}^-: \lambda\leq \frac{x}{5^k}\}$. Hence
$$\rho^{\mathcal{M}}(x)\leq\sum_{k} 2^k\rho^{\mathcal{P}^-}(\frac{x}{5^k}).$$
Denote $\lambda_*$ the least eigenvalue in $\mathcal{P}^-$. Then
$$\rho^{\mathcal{M}}(x)\leq\sum_{k=1}^{[\log(x/\lambda_*)/\log 5]} 2^k\rho^{\mathcal{P}^-}(\frac{x}{5^k})\leq c_6\cdot\sum_{k=1}^{[\log(x/\lambda_*)/\log 5]} 2^k(\frac{x}{5^k})^{\log2/\log 5}\leq c_7(\log x)x^{\log2/\log5},$$
for  an appropriate constant $c_7>0$.

Taking the above estimates into account, we  get
\begin{equation}\label{x}
\rho^\Omega(x)-\rho^\mathcal{L}(x)=O(x^{\log2/\log5}\log x) \quad
\mbox{ as } x\rightarrow\infty.
\end{equation}

On the other hand, by using the consequences of the usual spectral
decimation for the Dirichlet Laplacian on $\mathcal{SG}\setminus
V_0$, we will also prove that
\begin{equation}\label{y}
\rho^0(x)-\rho^\mathcal{L}(x)=O(x^{\log2/\log5}\log x) \quad \mbox{
as } x\rightarrow\infty.
\end{equation}
(a similar result can be found in Kigami $\cite{Ki5}$.)

In fact, as showed before, $\mathcal{L}$ is also a subset of the
Dirichlet spectrum $\mathcal{D}$ of $-\Delta$ on
$\mathcal{SG}\setminus V_0$. Besides $\mathcal{L}$, there remain
some $2$-series, $5$-series and $6$-series eigenvalues in
$\mathcal{D}$, whose associated eigenfunction having support
touching the line segment $L$. We denote them by $\mathcal{R}^2$,
$\mathcal{R}^5$ and $\mathcal{R}^6$ respectively. Thus we have
$$\rho^ 0(x)-\rho^{\mathcal{L}}(x)=\rho^{\mathcal{R}^2}(x)+\rho^{\mathcal{R}^5}(x)+\rho^{\mathcal{R}^6}(x),$$
where $\rho^{\mathcal{R}^2}(x), \rho^{\mathcal{R}^5}(x),
\rho^{\mathcal{R}^6}(x)$ are the eigenvalue counting functions of
the associated type eigenvalues.

For $\mathcal{R}^5$, we use $\mathcal{R}_m^5$ to denote the
associated total $m$-level graph eigenvalues. Then it is easy to
verify that $\mathcal{R}_m^5$ consists of $1+2^{m-1}$ initial
eigenvalues and $(m+1)2^{m-1}-2$ continued eigenvalues. Notice that
the bottom $\sharp\mathcal{R}_m^5$ eigenvalues in $\mathcal{R}^5$
are generated from eigenvalues in $\mathcal{R}_m^5$ by extending
these eigenvalues by choosing $\phi_-$ relations for all $m'>m$.
Since $5$ is the largest eigenvalue in $\mathcal{R}_m^5$, we get
$$\rho^{\mathcal{R}^5}(c_85^m)=\sharp\mathcal{R}_m^5=(1+2^{m-1})+((m+1)2^{m-1}-2)=(m+2)2^{m-1}-1,$$ for
appropriate constant $c_8>0$. Similarly as the analysis on
$\rho^{\mathcal{P}^+}(x)$, we then get
$$\rho^{\mathcal{R}^5}(x)=O(x^{\log2/\log 5}\log x)\quad \mbox{ as } x\rightarrow\infty.$$

Following similar argument, we also have
$$\rho^{\mathcal{R}^2}(x)=O(x^{\log2/\log 5})\quad \mbox{ as }
x\rightarrow\infty,$$ and
$$\rho^{\mathcal{R}^6}(x)=O(x^{\log2/\log 5}\log x)\quad \mbox{ as } x\rightarrow\infty.$$

Taking these estimates into account, we then get $(\ref{y})$.

Thus Theorem 3.10 follows from $(\ref{x})$ and $(\ref{y})$. $\Box$

\section{The Neumann case}

In this section, we give a brief discussion on the Neumann spectrum
of $-\Delta$ on $\Omega$. Throughout this section, for simplicity,
we omit the terms ``graph" and ``Neumann" without causing any
confusion. The main object  in this section is to prove Theorem 3.12
and Theorem 3.14. We will also give a comment on how to modify the
proof of Theorem 3.6 suitably to prove its Neumann counterpart,
Theorem 3.15, at the end of this section.

As indicated in Section 3, we want to impose a Neumann condition on
the graph $\Omega_m$ by imagining that it is embedded in a larger
graph by reflecting in each boundary vertex and imposing the
$\lambda_m$-eigenvalue equation on the even extension of $u_m$. It
is convenient to allow $m=1$, in which case there are only three
boundary points in $\Omega_1$ and no others. As introduced in
Section 3, $\mathcal{P}_m^{N}$ denotes the totality of primitive
Neumann eigenvalues of the discrete Laplacian $-\Delta_m$ on
$\Omega_m$. Due to the eigenspace dimensional counting argument in
Section 3, this time we need to find out $2^m$ symmetric primitive
eigenvalues and $2^m-1$ skew-symmetric primitive eigenvalues.

We focus our discussion on $\mathcal{P}_{m}^{+,N}$, the symmetric
case, and describe a similar weak spectral decimation which relates
$\mathcal{P}_m^{+,N}$ with $\mathcal{P}_{m+1}^{+,N}$.  Let $u_m$ be
a $\lambda_m$-eigenfunction of $-\Delta_m$ on $\Omega_m$ with
$\lambda_m\in\mathcal{P}_{m}^{+,N}$. Still denote by
$(b_0,b_1,\cdots,b_m)$ the values of $u_m$ on the skeleton of
$\Omega_m$. Write $\lambda_i^{(m)}$ the successor of $\lambda_m$ of
order $(m-i)$ with $1\leq i\leq m$. (This time we begin with
$\lambda_{1}^{(m)}$.) Assume that none of $\lambda_i^{(m)}$'s is
equal to $2$ or $5$ for $2\leq i\leq m$. Then $u_m$ is uniquely
determined by $(b_0,b_1,\cdots,b_m)$. In addition to the eigenvalue
equations at the vertex $F_1q_0,F_1^2q_0,\cdots,F_1^{m-1}q_0$ as
described in Section 4, we impose the equations
\begin{equation}\label{a}
(4-\lambda_1^{(m)})b_0=4b_1
\end{equation}
 at $q_0$ and
\begin{equation}\label{b}
(4-\lambda_m)b_m=2b_{m-1}+2b_m
\end{equation}
 at $F_1^mq_0$
according to the Neumann boundary condition. Hence
$(b_0,b_1,\cdots,b_m)$ can be viewed as a non-zero vector solution
of a system of equations consisting of $m+1$ equations in $m+1$
unknowns, whose determinant is
$$\left|
  \begin{array}{cccccccc}
  4-\lambda_1^{(m)} & -4\\
  l(\lambda_2^{(m)}) & s(\lambda_2^{(m)}) & r(\lambda_2^{(m)})   \\
       & \ddots & \ddots & \ddots &  \\
       &  & l(\lambda_m^{(m)}) & s(\lambda_{m}^{(m)})& r(\lambda_{m}^{(m)}) \\
   &&& -2& 2-\lambda_m^{(m)}
  \end{array}
\right|.$$ Hence $\lambda_{m}$ should be a solution of the following
equation
\begin{equation}\label{16}
q_m^N(x):=\left|
  \begin{array}{cccccccc}
  4-f^{(m-1)}(x) & -4\\
    l(f^{(m-2)}(x)) & s(f^{(m-2)}(x)) & r(f^{(m-2)}(x))  \\
         & \ddots & \ddots & \ddots &  \\
      &  & l(x) & s(x) &r(x)\\
    &&& -2& 2-x\\
  \end{array}
\right|=0.
\end{equation}

Thus if $\lambda_m$ is a root of $q_m^N(x)$  and none of
$f^{(i)}(\lambda_m)$'s with $0\leq i\leq m-2$ is equal to $2$ or
$5$, then $\lambda_m\in\mathcal{P}_m^{+,N}$. We should mention here
that when $m\geq 2$, comparing to $q_m(x)$ in the Dirichlet case,
$q_m^N(x)$ is a $(m+1)\times(m+1)$ tridiagonal determinant,
containing $q_m(x)$ in the center as a $(m-1)\times (m-1)$ minor.
Namely, we can write
$$q_m^N(x)=\left|
  \begin{array}{cccccccc}
  4-f^{(m-1)}(x) & -4 & 0 &\cdots & 0& 0\\
    l(f^{(m-2)}(x)) & &&&& 0  \\
    \vdots &&q_m(x)&&&\vdots\\
    0 &&&&&r(x)\\
    0 & 0 &\cdots & 0 & -2 & 2-x\\
  \end{array}
\right|.
$$
The degree of $q_m^N(x)$ is $3(2^{m-1}-1)+2^{m-1}+1=2^{m+1}-2$,
since the degree of $q_m(x)$ is $3(2^{m-1}-1)$. The analysis on
$q_m^N(x)$ is more complicated than that on $q_m(x)$ since for
$q_m(x)$ we can always use the expansion of $q_m(x)$ along the first
or last row to get a relation between two polynomials in same type
but with smaller degree.

The following lemma is a slight modification of the form of
$q_m^N(x)$ from a $(m+1)\times(m+1)$ determinant to a $m\times m$
determinant.

\textbf{Lemma 7.1.} \emph{Let $m\geq 2$. Then $$q_m^N(x)=(2-x)(x-6)\left|
  \begin{array}{cccccccc}
  4-f^{(m-1)}(x) & -4\\
    l(f^{(m-2)}(x)) & s(f^{(m-2)}(x)) & r(f^{(m-2)}(x))  \\
         & \ddots & \ddots & \ddots &  \\
      &  & l(f(x)) & s(f(x)) &r(f(x))\\
    &&& 1& f(x)-1\\
  \end{array}
\right|.
$$}

\emph{Proof.} Substituting the expression for $r(x)$ into
$(\ref{16})$, we get
\begin{eqnarray*}
q_m^N(x)&=&(2-x)\left|
  \begin{array}{cccccccc}
  4-f^{(m-1)}(x) & -4\\
    l(f^{(m-2)}(x)) & s(f^{(m-2)}(x)) & r(f^{(m-2)}(x))  \\
         & \ddots & \ddots & \ddots &  \\
      &  & l(x) & s(x) &-2(5-x)\\
    &&& -2& 1\\
  \end{array}
\right|\\&=&(2-x)\left|
  \begin{array}{cccccccc}
  4-f^{(m-1)}(x) & -4\\
    l(f^{(m-2)}(x)) & s(f^{(m-2)}(x)) & r(f^{(m-2)}(x))  \\
         & \ddots & \ddots & \ddots &  \\
      &  & l(x) & s(x)-4(5-x) &-2(5-x)\\
    &&& 0& 1\\
  \end{array}
\right|\\&=&(2-x)\left|
  \begin{array}{cccccccc}
  4-f^{(m-1)}(x) & -4\\
    l(f^{(m-2)}(x)) & s(f^{(m-2)}(x)) & r(f^{(m-2)}(x))  \\
         & \ddots & \ddots & \ddots &  \\
      &  & l(x) & s(x)-4(5-x) \\
  \end{array}
\right|.
\end{eqnarray*}

Noticing that $s(x)-4(5-x)=(x-6)(f(x)-1)$ and $l(x)=x-6$, we get the
desired result. $\Box$

The following lemma focuses on the possibility of the roots of
$q_m^N(x)$ satisfying $f^{(i)}(x)=2$ or $5$ for some $0\leq i\leq
m-2$.

\textbf{Lemma 7.2.} \emph{Let $m\geq 2$, and $x$ be a predecessor of
$2$ or $5$ of order $i$ with $0\leq i\leq m-2$. Then $q_m^N(x)=0$. }

\emph{Proof.} Firstly, let $x$ be a predecessor of $2$ of order $i$
with $0\leq i\leq m-2$. Then $f^{(i)}(x)=2$ and $f^{(i+1)}(x)=6$. If
$0\leq i<m-2$, the proof is the same as that of Lemma 4.2. So we
only need to check the $i=m-2$ case.  In this case we have
\begin{eqnarray*}
q_m^N(x)=\left|
  \begin{array}{cccccccc}
  4-f^{(m-1)}(x) & -4\\
    l(f^{(m-2)}(x)) & s(f^{(m-2)}(x)) & r(f^{(m-2)}(x))  \\
         & \ddots & \ddots & \ddots &  \\
       \end{array}
\right|=\left|
  \begin{array}{ccccccc}
         -2 & -4  \\
            -4 & -8 & 0 \\
     &   \ddots &\ddots & \ddots \\
  \end{array}
\right|=0.
\end{eqnarray*}

Secondly, let $x$ be a predecessor of $5$ of order $i$ with $0\leq
i\leq m-2$.   Then $f^{(i)}(x)=5$ and $f^{(i+1)}(x)=\cdots
f^{(m-1)}(x)=0$. Hence we have
\begin{eqnarray*}
q_m^N(x)=\left|
  \begin{array}{ccccccccccc}
      4 &  -4\\
      l(0)& s(0)&r(0) \\
      & \ddots & \ddots &\ddots\\
      && l(0)& s(0)&r(0)\\
   & &&  l(5) & s(5)& r(5)\\
   &&&&\ddots & \ddots &\ddots\\
  \end{array}
\right|=\left|
  \begin{array}{ccccccccccc}
      4 &  -4\\
      -6& 26&-20 \\
      & \ddots & \ddots &\ddots\\
      && -6& 26&-20\\
   & &&  -1 & 1& 0\\
   &&&&\ddots & \ddots &\ddots\\
  \end{array}
\right|=0.
\end{eqnarray*}
Thus we always have $q_m^N(x)=0$. $\Box$

This lemma means that for $m\geq 2$, all the predecessors of $2$ or
$5$ of order $i$ with $0\leq i\leq m-2$ are unwanted roots of
$q_m^N(x)$. To exclude them out, we define
$$p_m^N(x):=\frac{q_m^N(x)}{(x-2)(x-5)\cdots(f^{(m-2)}(x)-2)(f^{(m-2)}(x)-5)}\quad \mbox{ for  } m\geq 2,$$
and $$p_1^N(x):=q_1^N(x).$$ Now we can say if $\lambda_m$ is a root
of the polynomial $p_m^N(x)$, then
$\lambda_m\in\mathcal{P}_m^{+,N}$. It is easy to check that the
degree of $p_m^N(x)$ is $2^m$, since the degree of $q_m^N(x)$ is
$2^{m+1}-2$ and the number of all the unwanted roots of $q_m^N(x)$
is $2(1+2+\cdots 2^{m-2})=2^{m}-2$ for $m\geq 2$ and $0$ for $m=1$.
The following is an easy observation on $p_m^N(x)$.

\textbf{Lemma 7.3.} \emph{For $m\geq 1$, $p_m^N(x)$ always has roots
$0$ and $6$.}

\emph{Proof.} We only need to check $q_m^N(0)=q_m^N(6)=0$. It is
easy to see that
\begin{eqnarray*}
q_m^N(0)=\left|
  \begin{array}{ccccccccccc}
      4 &  -4\\
      l(0)& s(0)&r(0) \\
      & \ddots & \ddots &\ddots\\
      && l(0)& s(0)&r(0)\\
  &&&-2&2\\
  \end{array}
\right|=\left|
  \begin{array}{ccccccccccc}
      4 &  -4\\
      -6& 26&-20 \\
      & \ddots & \ddots &\ddots\\
      && -6& 26&-20\\
   &&& -2&2\\
  \end{array}
\right|=0.
\end{eqnarray*}

$q_m^N(6)=0$ follows from Lemma 7.1 for $m\geq 2$, and from direct
computation for $m=1. \Box$

In order to study the distribution of roots of $p_m^N(x)$, we now
introduce a type of auxiliary polynomials $l_m(x)$ associated to
$p_m^N(x)$. First, $\forall m\geq 1$, let $\widetilde{l}_m(x)$
denote the $m\times m$ minor  located in the upper left corner of
$q_m^N(x)$. Namely, $\widetilde{l}_1(x):=4-x$ and for $m\geq 2$,
$$
\widetilde{l}_m(x):=\left|
  \begin{array}{cccccccc}
  4-f^{(m-1)}(x) & -4\\
    l(f^{(m-2)}(x)) & s(f^{(m-2)}(x)) & r(f^{(m-2)}(x))  \\
         & \ddots & \ddots & \ddots &  \\

      &  & l(x) & s(x) &\\
  \end{array}
\right|.
$$
Note that the $(m-1)\times(m-1)$ minor located in the  bottom right
corner of $\widetilde{l}_m(x)$ is $q_m(x)$. The degree of
$\widetilde{l}_m(x)$ is $2^{m+1}-3$ since it is reduced by $1$
comparing to the degree of $q_m^N(x)$. With similar argument in the
proof of Lemma 7.2, all the predecessors of $2$ or $5$ of order $i$
with $0\leq i\leq m-2$ are roots of $\widetilde{l}_m(x)$. To exclude
them out, we define
$$l_m(x):=\frac{\widetilde{l}_m(x)}{(x-2)(x-5)\cdots(f^{(m-2)}(x)-2)(f^{(m-2)}(x)-5)}\quad \mbox{ for  } m\geq 2,$$
and $$l_1(x):=\widetilde{l}_1(x).$$ It is easy to check that the
degree of $l_m(x)$ is $2^m-1$, since the degree of
$\widetilde{l}_m(x)$ is $2^{m+1}-3$ and the number of all the
unwanted roots of $\widetilde{l}_m(x)$ is $2(1+2+\cdots
2^{m-2})=2^{m}-2$ for $m\geq 2$ and $0$ for $m=1$.

Based on the property
$$\widetilde{l}_m(x)=s(x)\widetilde{l}_{m-1}(f(x))-r(f(x))l(x)\widetilde{l}_{m-2}(f^{(2)}(x)),$$ $l_m(x)$ can be
analyzed in a similar way like $p_m(x)$ or $\widetilde{p}_m(x)$ in
the Dirichlet case. We then have:

 \textbf{Lemma
7.4.} \emph{$l_m(0)>0$ and $l_m(6)<0$, $\forall m\geq 1$.}

\emph{Proof.} $l_m(0)>0$ follows from a similar argument in the
proof of Proposition 4.1 and Proposition 4.2.

To prove $l_m(6)<0$, we only need to prove $\widetilde{l}_m(6)<0$ by
the definition of $l_m(x)$. It can be checked that
$\widetilde{l}_1(6)=-2<0$ and $\widetilde{l}_2(6)=-40<0$ by a direct
computation. For $m\geq 3$, an expansion of $\widetilde{l}_m(6)$
along the first row yields that
$$\widetilde{l}_m(6)=(4-f^{(m-1)}(6))q_m(6)+4(f^{(m-2)}(6)-6)q_{m-1}(6).$$
Recall that in the proof of Proposition 4.1(3), we have proved that
$q_{m}(6)\leq q_{m-1}(6)<0$, $\forall m\geq 3.$ Hence
$$\widetilde{l}_m(6)\leq(4-f^{(m-1)}(6)+4f^{(m-2)}(6)-24)q_m(6)=(f^{(m-2)}(6)-5)(f^{(m-2)}(6)+4)q_m(6)<0,$$
noticing that $f^{(m-2)}(6)\leq -6$ whenever $m\geq 3$. $\Box$

\textbf{Lemma 7.5.} \emph{For each $m\geq 1$, $l_m(x)$ has $2^m-1$
distinct real roots between $0$ and $6$ satisfying
$$0<\beta_{m,1}<\beta_{m,2}<\cdots<\beta_{m,2^m-1}<6.$$
Moreover, \begin{eqnarray*}
0&<&\beta_{m+1,1}<\phi_{-}(\beta_{m,1}),\\
\phi_{-}({\beta_{m,k-1}})&<&\beta_{m+1,k}<\phi_{-}(\beta_{m,k}),
\quad\forall 2\leq k\leq 2^m-1,\\
\phi_{-}({\beta_{m,2^m-1}})&<&\beta_{m+1,2^m}<\phi_{+}(\beta_{m,2^m-1}),\\
\phi_{+}({\beta_{m,2^{m+1}-k}})&<&\beta_{m+1,k}<\phi_{+}(\beta_{m,2^{m+1}-k-1}),
\quad\forall 2^m+1\leq k\leq 2^{m+1}-2,\\
\phi_{+}(\beta_{m,1})&<&\beta_{m+1,2^{m+1}-1}<6.
\end{eqnarray*}}

\emph{Proof.} It follows from a similar argument in the proof of
Lemma 4.4. $\Box$

The following lemma shows a relation between $p_m^N(x)$'s and
$l_m(x)$'s.

\textbf{Lemma 7.6.}\emph{ Let $m\geq 2$. Then
$p_m^N(x)=(2-x)l_m(x)-4l_{m-1}(f(x))$.}

\emph{Proof.} This is easy to get since we have
$$q_m^N(x)=(2-x)\widetilde{l}_m(x)-4(2-x)(5-x)\widetilde{l}_{m-1}(f(x)),\quad \forall m\geq 2,$$
using the expansion  along the last row of $q_m^N(x)$. $\Box$

Now we consider the distribution of roots of $p_m^N(x)$.

\textbf{Lemma 7.7.} \emph{For each $m\geq 1$, $p_m^N(x)$ has $2^m$
distinct roots between $0$ and $6$ (including $0$ and $6$).
Moreover, $p_m^N(0+)<0$ and $p_m^N(6-)<0$.}

\emph{Proof.} When $m=1$, it naturally holds.

Let $m\geq 2$. From Lemma 7.5, $l_m(x)$  has $2^m-1$ distinct real
roots between $0$ and $6$ satisfying
$$0<\beta_{m,1}<\beta_{m,2}<\cdots<\beta_{m,2^m-1}<6.$$
For each $1\leq k\leq 2^m-1$, using Lemma 7.6, we have
$$p_m^N(\beta_{m,k})=(2-\beta_{m,k})l_m(\beta_{m,k})-4l_{m-1}(f(\beta_{m,k}))=-4l_{m-1}(f(\beta_{m,k})).$$

When $k=1$, by Lemma 7.5, $0<\beta_{m,1}<\phi_-(\beta_{m-1,1})$,
hence $0<f(\beta_{m,1})<\beta_{m-1,1}$. Combined with $l_{m-1}(0)>0$
from Lemma 7.4, it follows $l_{m-1}(f(\beta_{m,1}))>0$, hence
$p_m^N(\beta_{m,1})<0$.

When $2\leq k\leq 2^{m-1}-1$, following from Lemma 7.5, we have
$\phi_-(\beta_{m-1,k-1})<\beta_{m,k}<\phi_-(\beta_{m-1,k})$, hence
$\beta_{m-1,k-1}<f(\beta_{m,k})<\beta_{m-1,k}$. Combined with
$l_{m-1}(0)>0$, it follows $l_{m-1}(f(\beta_{m,k}))\sim(-1)^{k-1}$,
hence $p_{m}^N(\beta_{m,k})\sim(-1)^k$.

When $k=2^{m-1}$, following from Lemma 7.5, we have
$\phi_-(\beta_{m-1,2^{m-1}-1})<\beta_{m,2^{m-1}}<\phi_+(\beta_{m-1,2^{m-1}-1})$,
hence $f(\beta_{m,2^{m-1}})>\beta_{m-1,2^{m-1}-1}$. Combined with
$l_{m-1}(0)>0$, it follows $l_{m-1}(f(\beta_{m,2^{m-1}}))<0$, hence
$p_{m}^N(\beta_{m,2^{m-1}})>0$.

When $2^{m-1}+1\leq k\leq 2^{m}-2$, following from Lemma 7.5, we
have
$\phi_+(\beta_{m-1,2^m-k})<\beta_{m,k}<\phi_+(\beta_{m-1,2^m-k-1})$,
hence $\beta_{m-1,2^m-k-1}<f(\beta_{m,k})<\beta_{m-1,2^m-k}$.
Combined with $l_{m-1}(0)>0$, it follows
$l_{m-1}(f(\beta_{m,k}))\sim(-1)^{k-1}$, hence
$p_{m}^N(\beta_{m,k})\sim(-1)^k$.

When $k= 2^{m}-1$, following from Lemma 7.5, we have
$\phi_+(\beta_{m-1,1})<\beta_{m,2^m-1}<6$, hence
$f(\beta_{m,2^m-1})<\beta_{m-1,1}$. Combined with $l_{m-1}(0)>0$, it
follows $l_{m-1}(f(\beta_{m,2^m-1}))>0$, hence
$p_{m}^N(\beta_{m,2^m-1})<0$.

Hence we have proved $p_m^N(\beta_{m,k})\sim(-1)^k, \forall 1\leq
k\leq 2^m-1.$ So there exist at least $2^m-2$ roots of $p_m^N(x)$,
each located strictly between each two consecutive $\beta_{m,k}$'s.
Moreover, Lemma 7.3 says that $0$ and $6$ are also roots of
$p_m^N(x)$. Thus we have found $2^m$ distinct roots of $p_m^N(x)$.
Since the order of $p_m^N(x)$ is also $2^m$, these are the total
roots of $p_m^N(x)$.

Furthermore, from the fact that $p_m^N(\beta_{m,1})<0$ and
$p_m^N(\beta_{m,2^m-1})<0$, we have $p_{m}^N(0+)<0$ and
$p_m^N(6-)<0$. $\Box$

In all that follows, we denote $$\lambda_{m,1}=0<
\lambda_{m,2}<\cdots<\lambda_{m,2^m-1}<\lambda_{m,2^m}=6$$ the $2^m$
distinct roots of $p_m^N(x)$ in increasing order, $\forall m\geq 1$.
In order to study the relation of roots of two consecutive
$p_m^N(x)$'s, we prove the following two lemmas:

\textbf{Lemma 7.8.}\emph{ Let $m\geq 1$ and $2\leq k\leq 2^{m}$,
then
$$p_{m+1}^N(\phi_-(\lambda_{m,k}))=-\frac{\lambda_{m,k}}{2}\frac{\phi_-(\lambda_{m,k})-6}{\phi_-(\lambda_{m,k})-5}\cdot l_m(\lambda_{m,k}),$$
and
$$p_{m+1}^N(\phi_+(\lambda_{m,k}))=-\frac{\lambda_{m,k}}{2}\frac{\phi_+(\lambda_{m,k})-6}{\phi_+(\lambda_{m,k})-5}\cdot l_m(\lambda_{m,k}).$$}

\emph{Proof.} For simplicity we only prove the first equality. The
second will follow from a similar argument. It is easy to see that
$\lambda_{m,k}$ is also a root of $q_m^N(x)$ and none of
$f^{(i)}(\lambda_{m,k})$'s ($0\leq i\leq m-2$) is equal to $2$ or
$5$.

By Lemma 7.1,
$q_{m+1}^N(\phi_-(\lambda_{m,k}))=(2-\phi_-(\lambda_{m,k}))(\phi_-(\lambda_{m,k})-6)\cdot
A$ where
$$
A=\left|
  \begin{array}{cccccccc}
  4-f^{(m-1)}(\lambda_{m,k}) & -4\\
    l(f^{(m-2)}(\lambda_{m,k})) & s(f^{(m-2)}(\lambda_{m,k})) & r(f^{(m-2)}(\lambda_{m,k}))  \\
         & \ddots & \ddots & \ddots &  \\
      &  & l(\lambda_{m,k}) & s(\lambda_{m,k}) &r(\lambda_{m,k})\\
    &&& 1& \lambda_{m,k}-1\\
  \end{array}
\right|.
$$
Noticing that from $q_{m}^N(\lambda_{m,k})=0$, we have
$$
\left|
  \begin{array}{cccccccc}
  4-f^{(m-1)}(\lambda_{m,k}) & -4\\
    l(f^{(m-2)}(\lambda_{m,k})) & s(f^{(m-2)}(\lambda_{m,k})) & r(f^{(m-2)}(\lambda_{m,k}))  \\
         & \ddots & \ddots & \ddots &  \\
      &  & l(\lambda_{m,k}) & s(\lambda_{m,k}) &r(\lambda_{m,k})\\
    &&& -1& 1-\lambda_{m,k}/2\\
  \end{array}
\right|=0.
$$
The summation of the above two determinants yields that
$A=\frac{\lambda_{m,k}}{2}\widetilde{l}_m(\lambda_{m,k})$. Hence
$$q_{m+1}^N(\phi_-(\lambda_{m,k}))=\frac{\lambda_{m,k}}{2}(2-\phi_-(\lambda_{m,k}))(\phi_-(\lambda_{m,k})-6)\cdot\widetilde{l}_m(\lambda_{m,k}),
$$
which yields the desired result. $\Box$

\textbf{Lemma 7.9.}\emph{ Let $m\geq 1$. Then
$(-1)^{k-1}l_m(\lambda_{m,k})>0,\quad \forall 1\leq k\leq 2^m$.}

\emph{Proof.} Let $\beta_{m,1},\beta_{m,2},\cdots\beta_{m,2^m-1}$
denote the $2^{m}-1$ distinct roots of $l_m(x)$ in increasing order
as described in Lemma 7.5. Then by the proof of Lemma 7.7, we have
$$\lambda_{m,1}=0<\beta_{m,1}<\lambda_{m,2}<\beta_{m,2}<\cdots<\lambda_{m,2^m-1}<\beta_{m,2^m-1}<\lambda_{m,2^m}=6.$$
Combined with the fact $l_m(\lambda_{m,1})=l_m(0)>0$  by Lemma 7.4,
it follows the desired result. $\Box$

Now we can prove the following  lemma.

 \textbf{Lemma 7.10.}
\emph{For each $m\geq 1$, $\mathcal{P}_m^{+,N}$ consists of at least
$2^m$ distinct eigenvalues satisfying
$$
\lambda_{m,1}=0<\lambda_{m,2}<\cdots<\lambda_{m,2^m-1}<\lambda_{m,2^m}=6.
$$
Moreover,
\begin{eqnarray}\label{80}
\phi_{-}({\lambda_{m,k-1}})&<&\lambda_{m+1,k}<\phi_{-}(\lambda_{m,k}),
\quad\forall 2\leq k\leq 2^m,\nonumber\\
\phi_{+}({\lambda_{m,2^{m+1}-k+1}})&<&\lambda_{m+1,k}<\phi_{+}(\lambda_{m,2^{m+1}-k}),
\quad\forall 2^m+1\leq k\leq 2^{m+1}-2,\nonumber\\
\phi_{+}(\lambda_{m,2})&<&\lambda_{m+1,2^{m+1}-1}<6.
\end{eqnarray}}

\emph{Proof.} Noticing that each root of $p_m^N(x)$ belongs to
$\mathcal{P}_m^{+,N}$, we only need to prove the results for the
roots of $p_m^N(x)$. The first statement follows from Lemma 7.7. We
now prove the second statement.
 From Lemma 7.8 and Lemma 7.9, we have
$$p_{m+1}^N(\phi_-(\lambda_{m,k}))=-\frac{\lambda_{m,k}}{2}\frac{\phi_-(\lambda_{m,k})-5}{\phi_-(\lambda_{m,k})-6}\cdot l_m(\lambda_{m,k})\sim -l_m(\lambda_{m,k})\sim(-1)^k,\quad\forall 2\leq k\leq 2^m,$$
and similarly,
$$p_{m+1}^N(\phi_+(\lambda_{m,k}))=-\frac{\lambda_{m,k}}{2}\frac{\phi_+(\lambda_{m,k})-5}{\phi_+(\lambda_{m,k})-6}\cdot l_m(\lambda_{m,k})\sim -l_m(\lambda_{m,k})\sim(-1)^k,\quad\forall 2\leq k\leq 2^m.$$
Following the above facts and Lemma 7.7, we can list the signs of
the values of $p_{m+1}^N(x)$ at different point $x$ in the following
table.
\begin{flushleft}
\scriptsize{
\begin{tabular}{cccccccccccccccc}
  \hline
  $x:$  &$0$ & $0+$ & $\phi_-(\lambda_{m,2})$ & $\phi_-(\lambda_{m,3})$ & $\cdots$ & $\phi_-(\lambda_{m,2^m})$ & $\phi_+(\lambda_{m,2^m})$ & $\cdots$ & $\phi_+(\lambda_{m,3})$ & $\phi_+(\lambda_{m,2})$ & $6-$&$6$\\
  $p_{m+1}^N(x):$  & $0$ &$-$ & $+$ & $-$ & $\cdots$ & $+$ & $+$ & $\cdots$ & $-$ & $+$ & $-$&$0$ \\
\end{tabular}}
\end{flushleft}
Hence there exist at least $2^{m+1}$ distinct roots of
$p_{m+1}^N(x)$ satisfying $(\ref{80})$. Moreover, these are the
totality of the roots of $p_{m+1}^N(x)$ since the degree of
$p_{m+1}^N(x)$ is also $2^{m+1}$. Hence we get the desired
distribution of roots of $p_{m+1}^N(x)$. $\Box$

The estimate $\phi_{+}(\lambda_{m,2})<\lambda_{m+1,2^{m+1}-1}<6$ in
Lemma 7.10 can be refined into
\begin{equation}\label{13}
\phi_{+}(\lambda_{m,2})<\lambda_{m+1,2^{m+1}-1}<5
\end{equation}
 by using the
following lemma.

\textbf{Lemma 7.11.} \emph{For $m\geq 2$, let
$\lambda_{m,1}=0,\lambda_{m,2},\cdots,\lambda_{m,2^m-1},\lambda_{m,2^m}=6$
be the $2^m$ distinct roots of $p_m^N(x)$ in increasing order. Then
$$\lambda_{m,k}+\lambda_{m,2^m-k+1}=5, \quad\forall 2\leq k\leq 2^m-1.$$}

\emph{Proof.} From Lemma 7.1, it is easy to see that if $q_m^N(x)=0$
and $x\neq 2$ or $6$, then $q_m^N(5-x)=0$. Obviously, each
$\lambda_{m,k}$ ($2\leq k\leq 2^m-1$) satisfies this property.
$\Box$

Hence we have the following Neumann analog of Lemma 4.5.

 \textbf{Lemma 7.12.} \emph{For each $m\geq 1$,
$\mathcal{P}_m^{+,N}$ consists of at least $2^m$ distinct
eigenvalues satisfying
$$\lambda_{m,1}=0<\lambda_{m,2}<\cdots<\lambda_{m,2^m-1}<5<\lambda_{m,2^m}=6.$$
Moreover,\begin{eqnarray}\label{71}
\phi_{-}({\lambda_{m,k-1}})&<&\lambda_{m+1,k}<\phi_{-}(\lambda_{m,k}),
\quad\forall 2\leq k\leq 2^m,\nonumber\\
\phi_{+}({\lambda_{m,2^{m+1}-k+1}})&<&\lambda_{m+1,k}<\phi_{+}(\lambda_{m,2^{m+1}-k}),
\quad\forall 2^m+1\leq k\leq 2^{m+1}-1.
\end{eqnarray} }

\emph{Proof.} It follows from Lemma 7.10 and Lemma 7.11. $\Box$

The following is a Neumann analog of Lemma 4.6.

\textbf{Lemma 7.13.} \emph{Let $\lambda_m$ be a root of $p_m^N(x)$,
$u_m$  a primitive $\lambda_m$-eigenfunction on $\Omega_m$, and
$(b_0,b_1,\cdots,b_m)$  the values of $u_m$ on the skeleton of
$\Omega_m$. Then $b_0\neq 0$ and $b_{m}\neq 0$.}

\emph{Proof.} Without loss of generality, assume $m\geq 3$. We still
use $\lambda_{i}^{(m)}$ to denote the successor of $\lambda_m$ of
order $(m-i)$ with $1\leq i\leq m$. From the definition of
$p_m^N(x)$, $\lambda_{i}^{(m)}\neq 2$ or $5$, for each $2\leq i\leq
m$. Vector $(b_0,b_1,\cdots,b_{m})$ can be viewed as a non-zero
vector solution of system $(\ref{9})$ of equations and in addition
the two Neumann boundary eigenvalue equations $(\ref{a})$ and
$(\ref{b})$.

Suppose $b_{m}=0$. Then from $(\ref{b})$, $b_{m-1}=0$. It is easy to
check that the determinant of the remaining equations in $m-1$
unknowns $(b_0,b_1,\cdots,b_{m-2})$ is
$\widetilde{l}_{m-1}(\lambda_{m-1}^{(m)})$. Since
$(b_0,b_1,\cdots,b_{m-2})$ should be a non-zero vector, we have
$\widetilde{l}_{m-1}(\lambda_{m-1}^{(m)})=0$, hence
$l_{m-1}(\lambda_{m-1}^{(m)})=0$. Noticing that from Lemma 7.6, we
have
$p_m^N(\lambda_m)=(2-\lambda_m)l_m(\lambda_m)-4l_{m-1}(\lambda_{m-1}^{(m)})$.
Hence we get that $l_m(\lambda_m)=0$ since both
$l_{m-1}(\lambda_{m-1}^{(m)})$ and $p_m^N(\lambda_m)$ are equal to
$0$. But this is impossible, since Lemma 7.5 says that  if
$l_{m-1}(\lambda_{m-1}^{(m)})=0$ then $l_m(\lambda_m)$ could not
equal to $0$. Hence $b_m\neq 0$.

On the other hand, if $b_0=0$, then by substituting it into the
system, noticing that none of $\lambda_{i}^{(m)}$'s is equal to $2$
or $5$, we can get $b_1=0,\cdots, b_{m}=0$ successively, which
contradicts to $b_{m}\neq 0$. Hence $b_0\neq 0$. $\Box$

This is the whole story of the symmetric case. The skew-symmetric
case is slightly different but very similar. The result is shown in
Section 3, but the proof is omitted.

\emph{Proof of Theorem 3.12 and Theorem 3.14.}

The results follows by using Lemma 7.12, Lemma 7.13 and their
skew-symmetric analogs, and the eigenspace dimension counting
formula $(\ref{123})$, following a similar argument for the
Dirichlet case. $\Box$

We should remark that Lemma 7.13 and its skew-symmetric analog show
that there is no primitive eigenfunction (or miniaturized
eigenfunction) that is simultaneously Dirichlet and Neumann ($D-N$).
Hence the only possible $D-N$ eigenfunctions are localized
eigenfunctions. This is same as the $\mathcal{SG}\setminus{V_0}$
case.

Before closing this section, we will make a comment on how to prove
Theorem 3.15, by a suitable modification of the proof of Theorem
3.6.

\emph{Proof of Theorem 3.15.}

Let $\{\lambda_m\}_{m\geq m_0}$ be a sequence of symmetric primitive
graph Neumann eigenvalues related by $\widetilde{\phi}_{\pm}$
relations, with all but a finite number of $\widetilde{\phi}_{-}$'s.
An argument similar to Lemma 5.1 says that the limit
$\lambda:=\frac{3}{2}\lim_{m\rightarrow\infty}5^m\lambda_m$ exists.
We will prove $\lambda\in \mathcal{P}^{+,N}$.

Let $m_1$ be the generation of fixation of $\{\lambda_m\}_{m\geq
m_0}$. For $m\geq m_1$, we still use $u_m$ to denote the associated
$\lambda_m$-eigenfunction on $\Omega_m$, and extend it to the whole
domain $\Omega$ using a similar recipe as for the Dirichlet case.
Then $u_m$ is also defined on $\Omega$, satisfying
\begin{eqnarray*}\left\{
\begin{array}{l}
  -\Delta u_m=5^m\Phi(\lambda_m) u_m \mbox{ on } \Omega,\\
  \partial_n u_m|_{\partial\Omega_m}=0,
\end{array}\right.
\end{eqnarray*}
with $5^m\Phi(\lambda_m)\rightarrow\lambda$ as $m$ goes to infinity.

As for the Dirichlet case,  defined $v_m:=\frac{u_m}{\parallel
u_m\parallel_{\infty}}$. With suitable modification of Lemma 5.7 and
Lemma 5.8, we still have the equicontinuity  of the sequence
$\{v_m\}_{m\geq m_1}$. Hence a similar argument still yields that
there exists a subsequence $\{v_{m_k}\}$ of $\{v_m\}$ which
converges uniformly to a continuous function $v$ on $\Omega$. Hence
we only need to prove that $v$ is the associated Neumann
eigenfunction of $\lambda$. We will use the notations defined in
Definition 3.11.

Fix an integer $n\geq m_1$, we still use $K_n$ to denote the part of
$\Omega$ above $\partial\Omega_n\setminus\{q_0\}$. Let
$G_{K_n}(x,y)$ denote the Green's function associate to the simple
domain $K_n$.

Then $\forall k$, we have
$$v_{m_k}(x)=\int_{K_n}G_{K_n}(x,y)5^{m_k}\Phi(\lambda_{m_k})v_{m_k}(y)d\mu(y)+h_{m_k}^{(n)}(x)\mbox{ on } K_n,$$
where $h_{m_k}^{(n)}$ is a harmonic function on $K_n$, taking the
same boundary values as $v_{m_k}$ on $\partial K_n$. Noticing that
the sequence $\{h_{m_k}^{(n)}\}$ converges uniformly on $K_n$ as $k$
goes to infinity. Let $k\rightarrow\infty$, we then get
$$v(x)=\lambda\int_{K_n}G_{K_n}(x,y)v(y)d\mu(y)+h^{(n)}\mbox{ on }
K_n,$$ where $h^{(n)}$ is the uniform limit of $\{h_{m_k}^{(n)}\}$
on $K_n$. Thus $-\Delta v=\lambda v$ on $K_n$.

Let $\varphi$ be a test function in $\widetilde{\mathcal{F}}$. Now
we calculate $\mathcal{E}_{K_n}(v,\varphi)$.

The Gauss-Green formula says that
$$\mathcal{E}_{K_n}(v,\varphi)=\lambda\int_{K_n}v\varphi
d\mu+\sum_{\partial K_n}\varphi\partial_n v.$$ It is easy to check
that $\partial_n v(q_0)=0$ since $\partial_n v_{m_k}(q_0)=0$ for
each $m_k$, by imaging that $v_{m_k}$ is defined on a large domain
by reflecting in $q_0$, being the even extension of itself. Hence we
have
\begin{equation}\label{xyz}
\mathcal{E}_{K_n}(v,\varphi)=\lambda\int_{K_n}v\varphi
d\mu+\partial_n v(F_1^nq_0)\cdot\sum_{\partial
K_n\setminus\{q_0\}}\varphi,
\end{equation}
 since obviously $v$ takes same value
along $\partial K_n\setminus\{q_0\}$.

On the other hand, the Gauss-Green formula also says that
$\int_{K_n}\Delta v d\mu=\sum_{\partial K_n}\partial_n v$, hence
$\partial_n v(F_1^nq_0)=-\frac{1}{2^n}\lambda\int_{K_n}v d\mu$.
Substituting it into $(\ref{xyz})$, we then have
\begin{eqnarray*}
\mathcal{E}_{K_n}(v,\varphi)&=&\lambda\int_{K_n}v\varphi
d\mu-\lambda\int_{K_n}v
d\mu\cdot(\frac{1}{2^n}\sum_{K_n\setminus\{q_0\}}\varphi)\\
&=&\lambda\int_{K_n}v\varphi d\mu-\lambda\int_{K_n}(v-v_{n})
d\mu\cdot(\frac{1}{2^n}\sum_{K_n\setminus\{q_0\}}\varphi).
\end{eqnarray*}
The last equality follows from the fact that $\int_{K_n}\Delta
v_nd\mu=\sum_{\partial K_n}\partial_nv_n=0$ and $-\Delta
v_n=5^n\Phi(\lambda_n)v_n$ on $K_n$.

Taking $n=m_k$, we then get $\forall k$,
$$\mathcal{E}_{K_{m_k}}(v,\varphi)
=\lambda\int_{K_{m_k}}v\varphi d\mu-\lambda\int_{K_{m_k}}(v-v_{m_k})
d\mu\cdot(\frac{1}{2^{m_k}}\sum_{K_{m_k}\setminus\{q_0\}}\varphi).
$$

Letting $k$ goes to infinity, we get
$$\widetilde{\mathcal{E}}(v,\varphi)
=\lambda\int_{\widetilde{\Omega}}v\varphi d\mu,$$ since $v_{m_k}$
converges uniformly to $v$ and
$|\frac{1}{2^{m_k}}\sum_{K_{m_k}\setminus\{q_0\}}\varphi|$ is
controlled by $\parallel\varphi\parallel_\infty$.

Hence $v$ is an eigenfunction of the Neumann Laplacian $-\Delta_N$
associated to $\lambda$ from the arbitrariness of the test function
$\varphi$. $\Box$

\section{Spectral asymptotics, ratio gaps and clusters}

In this section, we list some unproved conjectures related to the
structure  of the spectrum of $-\Delta$ on $\Omega$. For simplicity,
we only discuss the Dirichlet spectrum $\mathcal{S}$ on $\Omega$.

In Tables 8.1, 8.2, 8.3 and 8.4 we present the eigenvalues and their
multiplicities in $\mathcal{S}_m$ for level $m=2,3,4,5$, where we
use $\lambda_{m,k}^+$, $\lambda_{m,k}^-$, $\lambda_{m,k}$ to denote
the $k$'th $\mathcal{P}_m^+$, $\mathcal{P}_m^-$, $\mathcal{L}_m$
type eigenvalues respectively, and use
$\mathcal{M}_m(\lambda_{m',k}^-)$ to denote the miniaturized
eigenvalue generated from $\lambda_{m',k}^-$.

The following conjectures list some interesting phenomena we
observed from the data.

\textbf{Conjecture 8.1.} \emph{Let $\rho_m^\Omega(x)$ denote the
eigenvalue counting function of $\mathcal{S}_m$, i.e.,
$\rho_m^\Omega(x)=\sharp\{\lambda_m\in \mathcal{S}_m: \lambda_m\leq
x\}$. Then $\rho_m^\Omega(\phi_-^{(m-k)}(5))=3^k-2^k$ for $k<m$.}

\textbf{Remark.}\emph{ Here $3^k-2^k$ is the difference between
$a_k$ and $a_{k-1}$.}

This conjecture suggests that the bottom $3^k-2^k$ eigenvalues of
the Dirichlet spectrum of $\Omega$ should be generated from the
bottom $3^k-2^k$ eigenvalues in $\mathcal{S}_m$ and the largest of
these eigenvalues should be
$\lim_{n\rightarrow\infty}\frac{3}{2}5^n\phi_-^{(n-k)}(5)=c 5^k$ for
the appropriate choice of $c$. For the  eigenvalue counting function
$\rho^\Omega(x)$ on $\Omega$, we then have
$\rho^\Omega(c5^k)=3^k-2^k.$ Roughly this suggests an asymptotic
growth rate $\rho^\Omega(x)\sim x^{\log 3/\log 5}$ as
$x\rightarrow\infty$, which is similar to the
$\mathcal{SG}\setminus{V_0}$ case. But more precisely, this also
implies that
$$\rho^\Omega(x)=c_1x^{\log3/\log5}-c_2x^{\log2/\log 5}$$ along the sequence
$x=c5^k$ for some appropriate constants $c_1$ and $c_2$. Hence, in
analogy with the $\mathcal{SG}\setminus{V_0}$ case, we hopefully
believe the following more precise conjecture.

\textbf{Conjecture 8.2.} \emph{There exist two periodic functions
$g_1(t)$ and $g_2(t)$ of period $\log 5$, which are bounded above,
bounded away from zero, and necessarily discontinuous at the value
$\log c$, such that
\begin{equation}\label{aaa}
\rho^\Omega(x)=g_1(\log x)x^{\log 3/\log 5}+g_2(\log x)x^{\log2/\log
5}+o(x^{\log2/\log 5}).
\end{equation}}

Here, comparing to formula $(\ref{aaaa})$, besides the leading term
$g_1(\log x)x^{\log 3/\log 5}$, the asymptotic second term of the
eigenvalue counting function appears. This is very analogous to the
conjectures of Weyl and Berry.

\textbf{Conjecture 8.3.}\emph{ There exist gaps in the ratios of
eigenvalues from the Dirichlet spectrum $\mathcal{S}$ of $-\Delta$.
That is, we can find infinitely many pairs of consecutive
eigenvalues $\lambda$, $\lambda'$ with $\frac{\lambda'}{\lambda}\geq
c$ for some constant $c>1$.}

\textbf{Remark.} \emph{In fact, in the discrete spectrum
$\mathcal{S}_m$, one can observe that gap appears above  each
$\phi_-^{(m-k)}(5)$ for $k<m$. Moreover, there are also smaller gaps
below miniaturized eigenvalues.}

In $\cite{Bo}$ it was shown that on $\mathcal{SG}\setminus{V_0}$
there exist gaps in the ratios of eigenvalues. The existence of gaps
is an interesting phenomenon in itself, but it also has important
applications to analysis on fractals. See details in $\cite{Bo},
\cite{Ki5},\cite{Str5}$. Thus it is of great interest to know
whether similar phenomenon exists for fractals other than
$\mathcal{SG}$. In fact $\cite{Zhou}$ shows that this is the case
for Vicsek set. Also $\cite{Dren,Hare,Zhou1}$ investigates this
question for a variant of the $\mathcal{SG}$ type fractal.

\textbf{Conjecture 8.4.} \emph{In the spectrum $\mathcal{S}_m$,
between consecutive $5$ and $6$ type localized eigenvalues, there is
exactly one $\mathcal{P}^+$ and one $\mathcal{P}^-$ type eigenvalue
(except the case that the two consecutive eigenvalues are
$\phi_-(5)$ and $\phi_+(6)=3$, where there is nothing in between).}

\textbf{Conjecture 8.5.} \emph{In the spectrum $\mathcal{S}_m$, the
number of distinct eigenvalues between $5-\varepsilon$ and $5$ goes
to $\infty$ as $m\rightarrow\infty$ for any $\varepsilon>0$.}

\textbf{Remark.} \emph{This means in $\mathcal{S}$ there exist
eigenvalue clusters, that is, arbitrarily many distinct eigenvalues
in an arbitrarily small interval.}

We say the spectrum $\mathcal{S}$ exhibits spectral clustering.
Clustering  does not occur on the $\mathcal{SG}\setminus{V_0}$ case.
Experimental evidence suggests that it does occur on the pentagasket
$\cite{Ada}$ and on the Julia sets $\cite{Flo}$. It is proved that
in $\cite{Con}$ it also does occur on Vicsek set.

\begin{flushleft}
\begin{center}
\scriptsize{
\begin{tabular}{c|cccc|c|cccc}
  \hline
$m=2$ & eigenvalue $\lambda_m$ & $\frac{3}{2}5^m\lambda_m$ & multi & type & $m=2$ & eigenvalue $\lambda_m$ & $\frac{3}{2}5^m\lambda_m$ & multi & type\\
  1 & $\lambda^+_{2,1}$=1.064568 & 39.92& 1 & $\mathcal{P}^{+}$ & 4 & $\lambda^+_{2,3}$=5.472834 & 205.23& 1 & $\mathcal{P}^+$ \\
  2 & $\lambda^-_{2,1}$=3.381966 &126.82& 1 & $\mathcal{P}^-$  &5 & $\lambda^-_{2,2}$=5.618034 & 210.68&1 & $\mathcal{P}^-$ \\
  3 & $\lambda^+_{2,2}$=4.462598 &167.35& 1 & $\mathcal{P}^+$ && \\
 \hline
\end{tabular}}
\begin{center}
\textbf{Table 8.1.} \small{The $2$-level eigenvalues in
$\mathcal{S}_2$ in increasing order.}
\end{center}
\end{center}
\end{flushleft}
\begin{flushleft}
\begin{center}
 \scriptsize{
\begin{tabular}{c|cccc|c|cccc}
  \hline
$m=3$ & eigenvalue $\lambda_m$ & $\frac{3}{2}5^m\lambda_m$ & multi & type &$m=3$ & eigenvalue $\lambda_m$ & $\frac{3}{2}5^m\lambda_m$ & multi & type\\
  1 & $\lambda^+_{3,1}$=0.187518 & 35.16& 1 & $\mathcal{P}^{+}$  & 11 & $\lambda^-_{3,4}$=3.902230 &731.67& 1 & $\mathcal{P}^-$ \\
  2 & $\lambda^-_{3,1}$=0.558733 &104.76& 1 & $\mathcal{P}^-$ &   12 & $\lambda^+_{3,6}$=4.517231 & 846.98&1 & $\mathcal{P}^+$ \\
  3 & $\lambda^+_{3,2}$=0.805532 & 151.04& 1 & $\mathcal{P}^+$ &13 & $\lambda^-_{3,5}$=4.803115 &900.58& 1 & $\mathcal{P}^-$ \\
  4 & $\lambda^-_{3,2}$=1.247636 & 233.93& 1 & $\mathcal{P}^-$ & 14 & $\lambda^+_{3,7}$=4.946726 &927.51& 1 & $\mathcal{P}^+$ \\
  5 & $\lambda^+_{3,3}$=1.329287 & 249.24& 1 & $\mathcal{P}^+$&  15 & $\lambda_{3,1}$=5 &937.50& 1 & $\mathcal{L}$ \\
  6 & $\lambda^-_{3,3}$=3.059152 &573.59& 1 & $\mathcal{P}^-$& 16 & $\lambda^+_{3,8}$=5.424059 &1017.01& 1 & $\mathcal{P}^+$ \\
  7 & $\lambda^+_{3,4}$=3.075910 & 576.73&1 & $\mathcal{P}^+$ & 17 & $\lambda^-_{3,6}$=5.429135 &1017.96& 1 & $\mathcal{P}^-$ \\
  8,9 & $\mathcal{M}_3(\lambda^-_{2,1})$=3.381966 &634.12& 2 &  $\mathcal{M}$& 18,19& $\mathcal{M}_3(\lambda^-_{2,2})$=5.618034 &1053.38& 2 & $\mathcal{M}$ \\
  10 & $\lambda^+_{3,5}$=3.713736 &696.33& 1 & $\mathcal{P}^+$ &20--24 & $\lambda_{3,2}$=6 &1125.00& 5 & $\mathcal{L}$\\
   \hline
\end{tabular}}
\begin{center}
\textbf{Table 8.2.} \small{The $3$-level eigenvalues in
$\mathcal{S}_3$ in increasing order.}
\end{center}
\end{center}
\end{flushleft}
\begin{flushleft}
\begin{center} \scriptsize{
\begin{tabular}{c|cccc|c|cccc}
  \hline
$m=4$ & eigenvalue $\lambda_m$ & $\frac{3}{2}5^m\lambda_m$ & multi & type & $m=4$ & eigenvalue $\lambda_m$ & $\frac{3}{2}5^m\lambda_m$ & multi & type\\
  1 & $\lambda^+_{4,1}$=0.035755 &33.52& 1 & $\mathcal{P}^{+}$&34 & $\lambda^+_{4,10}$=3.631877 &3404.88& 1 & $\mathcal{P}^+$ \\
  2 & $\lambda^-_{4,1}$=0.100554 &94.27& 1 & $\mathcal{P}^-$& 35 & $\lambda^-_{4,8}$=3.656967 &3428.41& 1 & $\mathcal{P}^-$ \\
  3 & $\lambda^+_{4,2}$=0.146945 &137.76& 1 & $\mathcal{P}^+$&36 & $\lambda^+_{4,11}$=3.760496 &3525.46& 1 & $\mathcal{P}^+$ \\
  4 & $\lambda^-_{4,2}$=0.249495 & 233.90&1 & $\mathcal{P}^-$&37,38 & $\mathcal{M}_4(\lambda^-_{3,4})$=3.902230 &3658.34& 2 & $\mathcal{M}$ \\
  5 & $\lambda^+_{4,3}$=0.277423 & 260.08&1 & $\mathcal{P}^+$& 39 & $\lambda^-_{4,9}$=3.982762 &3733.84& 1 & $\mathcal{P}^-$ \\
  6,7 & $\mathcal{M}_4(\lambda^-_{3,1})$=0.558733 & 523.81&2 &  $\mathcal{M}$&40 & $\lambda^+_{4,12}$=4.074531 &3819.87& 1 & $\mathcal{P}^+$ \\
  8 & $\lambda^+_{4,4}$=0.645454 &605.11& 1 & $\mathcal{P}^+$&41 & $\lambda^+_{4,13}$=4.223191 & 3959.24&1 & $\mathcal{P}^+$ \\
  9 & $\lambda^-_{4,3}$=0.652593 &611.81& 1 & $\mathcal{P}^-$&42 & $\lambda^-_{4,10}$=4.241362 &3976.28& 1 & $\mathcal{P}^-$ \\
  10 & $\lambda^-_{4,4}$=0.843591 &790.87& 1 & $\mathcal{P}^-$&43 & $\lambda^-_{4,11}$=4.573615 &4287.76& 1 & $\mathcal{P}^-$ \\
  11 & $\lambda^+_{4,5}$=0.857718 &804.11& 1 & $\mathcal{P}^+$& 44 & $\lambda^+_{4,14}$=4.586787 &4300.11& 1 & $\mathcal{P}^+$ \\
  12 & $\lambda^+_{4,6}$=0.965805 &905.44& 1 & $\mathcal{P}^+$&45 & $\lambda^+_{4,15}$=4.735683 & 4439.70&1 & $\mathcal{P}^+$ \\
  13 & $\lambda^-_{4,5}$=1.065699 &999.09& 1 & $\mathcal{P}^-$&46 & $\lambda^-_{4,12}$=4.793032 &4493.47& 1 & $\mathcal{P}^-$ \\
  14,15 & $\mathcal{M}_4(\lambda^-_{3,2})$=1.247636 &1169.66& 2 &  $\mathcal{M}$& 47,48 & $\mathcal{M}_4(\lambda^-_{3,5})$=4.803115 &4502.92& 2 & $\mathcal{M}$ \\
  16 & $\lambda^+_{4,7}$=1.263652 &1184.67& 1 & $\mathcal{P}^+$& 49 & $\lambda^+_{4,16}$=4.926848 &4618.92& 1 & $\mathcal{P}^+$ \\
  17 & $\lambda^-_{4,6}$=1.358256 &1273.37& 1 & $\mathcal{P}^-$& 50 & $\lambda^-_{4,13}$=4.979948 &4668.70& 1 & $\mathcal{P}^-$ \\
  18 & $\lambda^+_{4,8}$=1.372367 &1286.59& 1 & $\mathcal{P}^+$&51 & $\lambda^+_{4,17}$=4.993259 &4681.18& 1 & $\mathcal{P}^+$ \\
  19 & $\lambda_{4,1}$=1.381966 &1295.59& 1 & $\mathcal{L}$&52--57 & $\lambda_{4,4}$=5 &4687.50& 6 & $\mathcal{L}$ \\
  20--24 & $\lambda_{4,2}$=3 &2812.50& 5 & $\mathcal{L}$&58 & $\lambda^+_{4,18}$=5.423778 &5084.79& 1 & $\mathcal{P}^+$ \\
  25,26 & $\mathcal{M}_4(\lambda^-_{3,3})$=3.059152 &2867.96& 2 &  $\mathcal{M}$&59 & $\lambda^-_{4,14}$=5.423779 &5084.79& 1 & $\mathcal{P}^-$ \\
  27 & $\lambda^-_{4,7}$=3.078348 &2885.95& 1 & $\mathcal{P}^-$ & 60,61 & $\mathcal{M}_4(\lambda^-_{3,6})$=5.429135 & 5089.81&2 & $\mathcal{M}$ \\
  28 & $\lambda^+_{4,9}$=3.078431 &2886.03& 1 & $\mathcal{P}^+$&62--65 & $\mathcal{M}_4(\lambda^-_{2,2})$=5.618034 &5266.91& 4 & $\mathcal{M}$ \\
  29--32 & $\mathcal{M}_4(\lambda^-_{2,1})$=3.381966 &3170.59& 4 & $\mathcal{M}$ &66--89 & $\lambda_{4,5}$=6 &5625.00& 24 & $\mathcal{L}$\\
   33 & $\lambda_{4,3}$=3.618034 &3391.91& 1 & $\mathcal{L}$ &&\\
  \hline
\end{tabular}}
\begin{center}
\textbf{Table 8.3.} \small{The $4$-level eigenvalues in
$\mathcal{S}_4$ in increasing order.}
\end{center}
\end{center}
\end{flushleft}
\begin{flushleft}
\begin{center}
\scriptsize{
\begin{tabular}{c|cccc|c|cccc}
  \hline
$m=5$ & eigenvalue $\lambda_m$ & $\frac{3}{2}5^m\lambda_m$ & multi & type & $m=5$ & eigenvalue $\lambda_m$ & $\frac{3}{2}5^m\lambda_m$ & multi & type\\
 1 & $\lambda^+_{5,1}$=0.007039 & 33.00& 1 & $\mathcal{P}^+$& 112 & $\lambda^+_{5,20}$=3.620288 &16970.1& 1 & $\mathcal{P}^+$\\
 2 & $\lambda^-_{5,1}$=0.019385 & 90.87& 1 & $\mathcal{P}^-$& 113 & $\lambda^-_{5,16}$=3.623927 & 16987.2& 1 & $\mathcal{P}^-$\\
3 & $\lambda^+_{5,2}$=0.028430 & 133.27& 1 & $\mathcal{P}^+$&114 & $\lambda^+_{5,21}$=3.644882 & 17085.4& 1 & $\mathcal{P}^+$\\
 4 & $\lambda^-_{5,2}$=0.049571 & 232.36& 1 & $\mathcal{P}^-$&  115,116 & $\mathcal{M}_5(\lambda^-_{4,8})$=3.656967 & 17142.0& 2 & $\mathcal{M}$\\
5 & $\lambda^+_{5,3}$=0.055860 & 261.84& 1 & $\mathcal{P}^+$& 117 & $\lambda^-_{5,17}$=3.694772 & 17319.2& 1 & $\mathcal{P}^-$\\
 6,7 & $\mathcal{M}_5(\lambda^-_{4,1})$=0.100554 & 471.35& 2 & $\mathcal{M}$&118 & $\lambda^+_{5,22}$=3.720985 & 17442.1& 1 & $\mathcal{P}^+$\\
8 & $\lambda^+_{5,4}$=0.123515 & 578.98& 1 & $\mathcal{P}^+$&119 & $\lambda^+_{5,23}$=3.749413 & 17575.4& 1 & $\mathcal{P}^+$\\
 9 & $\lambda^-_{5,3}$=0.125398 &587.80& 1 & $\mathcal{P}^-$&  120 & $\lambda^-_{5,18}$=3.753145 & 17592.9& 1 & $\mathcal{P}^-$\\
  10 & $\lambda^-_{5,4}$=0.166319 & 779.62& 1 & $\mathcal{P}^-$& 121--124 & $\mathcal{M}_5(\lambda^-_{3,4})$=3.902230 & 18291.7& 4 & $\mathcal{M}$\\
11 & $\lambda^+_{5,5}$=0.170850 & 800.86& 1 & $\mathcal{P}^+$&125 & $\lambda^-_{5,19}$=3.908588 & 18321.5& 1 & $\mathcal{P}^-$\\
12 & $\lambda^+_{5,6}$=0.196017 & 918.83& 1 & $\mathcal{P}^+$&126 & $\lambda^+_{5,24}$=3.912510 & 18339.9& 1 & $\mathcal{P}^+$\\
 13 & $\lambda^-_{5,5}$=0.217665 & 1020.30& 1 & $\mathcal{P}^-$&127 & $\lambda^+_{5,25}$=3.971467 & 18616.3& 1 & $\mathcal{P}^+$\\
 14,15 & $\mathcal{M}_5(\lambda^-_{4,2})$=0.249495 &1169.51& 2 & $\mathcal{M}$&  128,129 & $\mathcal{M}_5(\lambda^-_{4,9})$=3.982762 & 18669.2& 2 & $\mathcal{M}$\\
16 & $\lambda^+_{5,7}$=0.264441 & 1239.57& 1 & $\mathcal{P}^+$& 130 & $\lambda^-_{5,20}$=3.997137 & 18736.6& 1 & $\mathcal{P}^-$\\
 17 & $\lambda^-_{5,6}$=0.286684 & 1343.83& 1 & $\mathcal{P}^-$& 131 & $\lambda^+_{5,26}$=4.069518 & 19075.9& 1 & $\mathcal{P}^+$\\
18 & $\lambda^+_{5,8}$=0.290993 & 1364.03& 1 & $\mathcal{P}^+$& 132 & $\lambda^-_{5,21}$=4.103862 & 19236.9& 1 & $\mathcal{P}^-$\\
  19 & $\lambda_{5,1}$=0.293638 & 1376.43& 1 & $\mathcal{L}$&133 & $\lambda^+_{5,27}$=4.116582 & 19296.5& 1 &  $\mathcal{P}^+$\\
 20--23 & $\mathcal{M}_5(\lambda^-_{3,1})$=0.558733 & 2619.06& 4 & $\mathcal{M}$&134 & $\lambda_{5,7}$=4.122334 &19323.4& 1 & $\mathcal{L}$\\
24 & $\lambda^+_{5,9}$=0.644676 & 3021.92& 1 & $\mathcal{P}^+$&135 & $\lambda^+_{5,28}$=4.219041 &19776.8& 1 & $\mathcal{P}^+$\\
 25 & $\lambda^-_{5,7}$=0.644693 & 3022.00& 1 & $\mathcal{P}^-$& 136 & $\lambda^-_{5,22}$=4.219295 & 19777.9& 1 & $\mathcal{P}^-$\\
 26,27 & $\mathcal{M}_5(\lambda^-_{4,3})$=0.652593 & 3059.03& 2 & $\mathcal{M}$& 137,138 & $\mathcal{M}_5(\lambda^-_{4,10})$=4.241362 &19881.4& 2 & $\mathcal{M}$\\
28--32 & $\lambda_{5,2}$=0.697224 & 3268.24& 5 & $\mathcal{L}$&139--143 & $\lambda_{5,8}$=4.302776 & 20169.3& 5 & $\mathcal{L}$\\
 33,34 & $\mathcal{M}_5(\lambda^-_{4,4})$=0.843591 & 3954.33& 2 & $\mathcal{M}$&  144,145 & $\mathcal{M}_5(\lambda^-_{4,11})$=4.573615 & 21438.8& 2 & $\mathcal{M}$\\
 35 & $\lambda^-_{5,8}$=0.864034 &4050.16& 1 & $\mathcal{P}^-$&146 & $\lambda^-_{5,23}$=4.588806 & 21510.0& 1 & $\mathcal{P}^-$\\
36 & $\lambda^+_{5,10}$=0.866936 &4063.76& 1 & $\mathcal{P}^+$&147 & $\lambda^+_{5,29}$=4.588882 & 21510.4& 1 & $\mathcal{P}^+$\\
37 & $\lambda_{5,3}$=0.877666 & 4114.06& 1 & $\mathcal{L}$&148 & $\lambda_{5,9}$=4.706362 & 22061.1& 1 & $\mathcal{L}$\\
38 & $\lambda^+_{5,11}$=0.890579 & 4174.59& 1 & $\mathcal{P}^+$&149 & $\lambda^+_{5,30}$=4.710126 & 22078.7& 1 & $\mathcal{P}^+$\\
 39 & $\lambda^-_{5,9}$=0.921042 & 4317.38& 1 & $\mathcal{P}^-$& 150 & $\lambda^-_{5,24}$=4.717827 & 22114.8& 1 & $\mathcal{P}^-$\\
 40 & $\lambda^+_{5,12}$=0.951360 &4459.50& 1 & $\mathcal{P}^+$&151 & $\lambda^+_{5,31}$=4.742035 & 22228.3& 1 & $\mathcal{P}^+$\\
 41 & $\lambda^-_{5,10}$=1.013289 & 4749.79& 1 & $\mathcal{P}^-$& 152 & $\lambda^-_{5,25}$=4.791572 & 22460.5& 1 & $\mathcal{P}^-$\\
42 & $\lambda^+_{5,13}$=1.031636 & 4835.79& 1 & $\mathcal{P}^+$&   153,154 & $\mathcal{M}_5(\lambda^-_{4,12})$=4.793032 & 22467.3& 2 & $\mathcal{M}$\\
43,44 & $\mathcal{M}_5(\lambda^-_{4,5})$=1.065699 & 4995.46& 2 & $\mathcal{M}$&   155--158 & $\mathcal{M}_5(\lambda^-_{3,5})$=4.803115 & 22514.6& 4 & $\mathcal{M}$\\
45 & $\lambda^+_{5,14}$=1.095777 & 5136.45& 1 & $\mathcal{P}^+$&159 & $\lambda^+_{5,32}$=4.809185 & 22543.1& 1 & $\mathcal{P}^+$\\
 46 & $\lambda^-_{5,11}$=1.097686 & 5145.40& 1 & $\mathcal{P}^-$&160 & $\lambda^+_{5,33}$=4.844770 & 22709.9& 1 & $\mathcal{P}^+$\\
  47--50 & $\mathcal{M}_5(\lambda^-_{3,2})$=1.247636 &5848.29& 4 &  $\mathcal{M}$& 161 & $\lambda^-_{5,26}$=4.847489 & 22722.6& 1 & $\mathcal{P}^-$\\
51 & $\lambda^+_{5,15}$=1.259109 & 5902.07& 1 & $\mathcal{P}^+$&  162 & $\lambda^-_{5,27}$=4.932207 & 23119.7& 1 & $\mathcal{P}^-$\\
 52 & $\lambda^-_{5,12}$=1.260744 & 5909.74& 1 & $\mathcal{P}^-$&163 & $\lambda^+_{5,34}$=4.934639 & 23131.1& 1 & $\mathcal{P}^+$\\
53 & $\lambda^+_{5,16}$=1.291565 &6054.21& 1 & $\mathcal{P}^+$&164 & $\lambda^+_{5,35}$=4.950036 & 23203.3& 1 & $\mathcal{P}^+$\\
 54 & $\lambda^-_{5,13}$=1.314754 & 6162.91& 1 & $\mathcal{P}^-$& 165 & $\lambda^-_{5,28}$=4.963126 & 23264.7& 1 & $\mathcal{P}^-$\\
55 & $\lambda^+_{5,17}$=1.358055 &6365.88& 1 & $\mathcal{P}^+$&    166,167 & $\mathcal{M}_5(\lambda^-_{4,13})$=4.979948 &23343.5& 2 & $\mathcal{M}$\\
  56,57 & $\mathcal{M}_5(\lambda^-_{4,6})$=1.358256 & 6366.83& 2 &  $\mathcal{M}$&168 & $\lambda^+_{5,36}$=4.987488 &23378.9& 1 & $\mathcal{P}^+$\\
 58 & $\lambda^-_{5,14}$=1.377582 & 6457.42& 1 & $\mathcal{P}^-$& 169 & $\lambda^-_{5,29}$=4.997193 &23424.3& 1 & $\mathcal{P}^-$\\
59 & $\lambda^+_{5,18}$=1.380161 & 6469.50& 1 & $\mathcal{P}^+$&170 & $\lambda^+_{5,37}$=4.998947 &23432.6& 1 & $\mathcal{P}^+$\\
60--65 & $\lambda_{5,4}$=1.381966 &6477.97& 6 & $\mathcal{L}$&171--195 & $\lambda_{5,10}$=5 & 23437.5& 25 & $\mathcal{L}$\\
66--89 & $\lambda_{5,5}$=3 & 14063.0& 24 & $\mathcal{L}$&196 & $\lambda^+_{5,38}$=5.423778 & 25424.0& 1 & $\mathcal{P}^+$\\
   90--93 & $\mathcal{M}_5(\lambda^-_{3,3})$=3.059152 & 14339.8& 4 &   $\mathcal{M}$& 197 & $\lambda^-_{5,30}$=5.423778 & 25424.0& 1 & $\mathcal{P}^-$\\
   94,95 & $\mathcal{M}_5(\lambda^-_{4,7})$=3.078348 & 14429.8& 2 &   $\mathcal{M}$&  198,199 & $\mathcal{M}_5(\lambda^-_{4,14})$=5.423779 & 25424.0& 2 & $\mathcal{M}$\\
96 & $\lambda^+_{5,19}$=3.078432 & 14430.2& 1 & $\mathcal{P}^+$&    200--203 & $\mathcal{M}_5(\lambda^-_{3,6})$=5.429135 & 25449.1& 4 & $\mathcal{M}$\\
 97 & $\lambda^-_{5,15}$=3.078432 & 14430.2& 1 & $\mathcal{P}^-$&  204--211 & $\mathcal{M}_5(\lambda^-_{2,2})$=5.618034 & 26334.5& 8 & $\mathcal{M}$\\
  98--105 & $\mathcal{M}_5(\lambda^-_{2,1})$=3.381966 & 15853.0& 8 &  $\mathcal{M}$&212--300 & $\lambda_{5,11}$=6 & 28125.0& 89 & $\mathcal{L}$\\
106--111 & $\lambda_{5,6}$=3.618034 & 16959.5& 6 & $\mathcal{L}$&&\\
 \hline
\end{tabular}}
\begin{center}
\textbf{Table 8.4.} \small{The $5$-level eigenvalues in
$\mathcal{S}_5$ in increasing order.}
\end{center}
\end{center}
\end{flushleft}

\section{Further discussion}
In this section, we discuss to what extent our method can be
extended to other domains in $\mathcal{SG}$. In particular, we will
focus on $\Omega_x$ ($0<x<1$).  It seems that we can analyze the
spectrum of $-\Delta$ on $\Omega_x$ case by case following the
similar recipe for the $\Omega_1$ case. However, it is hardly to
develop a general method which is suitable for all cases, although
we believe that we could make clear the structures of the spectra.
We let $L_x$ denote the bottom boundary of $\Omega_x$. Thus $L_x$
will be a Cantor set for generic $x$, and an union of intervals if
$x$ is a dyadic rational. We may assume without loss of generality
that $\frac{1}{2}<x<1$, for if not we may first solve the problem
for $\Omega_{2x}$, and then simply dilate the solution to
$\Omega_x$.

For simplicity, we only discuss the Dirichlet spectrum of $-\Delta$.
Obviously, it will suffice to describe the discrete Dirichlet
spectra of $-\Delta_m$'s for all $m$. Hence the first problem is how
to define the graph approximations. Similarly to  $\Omega_1$, the
fractal domain $\Omega_x$ can be realized as the limit of a sequence
of graphs $\Omega_{x,m}$. More precisely, $\forall m\geq 1$, let
$V_m^{\Omega_x}$ be a subset of $V_m$ with all vertices lying along
or under $L_x$  removed. Let $\Omega_{x,m}$ be the subgraph of
$\Gamma_m$ restricted to $V_m^{\Omega_x}$. Denote by
$\partial\Omega_{x,m}$ the boundary of the finite graph
$\Omega_{x,m}$. It is easy to find that
$V_m^{\Omega_x}\setminus\partial\Omega_{x,m}$,
$\partial\Omega_{x,m}$ approximate to $\Omega_x$ and
$\partial\Omega_x$ as $m$ goes to infinity respectively. See Fig.
9.1 and Fig. 9.2 for $\Omega_x$ and $\Omega_{x,m}$ where $x=3/4$.

On $\Omega_{x,m}$ the Dirichlet $\lambda_m$-eigenvalue equations
consists of exactly $\sharp
(V_m^{\Omega_x}\setminus\partial\Omega_{x,m})$ equations in $\sharp
(V_m^{\Omega_x}\setminus\partial\Omega_{x,m})$ unknowns. We denote
by $\mathcal{S}_m(x)$ the spectrum of $-\Delta_m$ on $\Omega_{x,m}$
for each $m\geq 1$. Accordingly, $\mathcal{S}_m(x)$ should consists
 of (at least) three types of eigenvalues, denoted by
$\mathcal{L}_m(x)$, $\mathcal{P}_m(x)$ and $\mathcal{M}_m(x)$
respectively. $\mathcal{P}_m(x)$ can also be split into symmetric
part $\mathcal{P}_m^+(x)$ and skew-symmetric part
$\mathcal{P}_m^-(x)$. We omit the precise definitions since they are
obvious. To ensure that there is no other eigenvalue in
$\mathcal{S}_m(x)$, the following eigenspace dimensional counting
formula is hoped to be held,
$$\sharp
(V_m^{\Omega_x}\setminus\partial\Omega_{x,m})=\sharp\mathcal{L}_m(x)+\sharp\mathcal{P}_m(x)+\sharp\mathcal{M}_m(x).$$
\begin{figure}[ht]
\begin{center}
\includegraphics[width=6cm,totalheight=5.3cm]{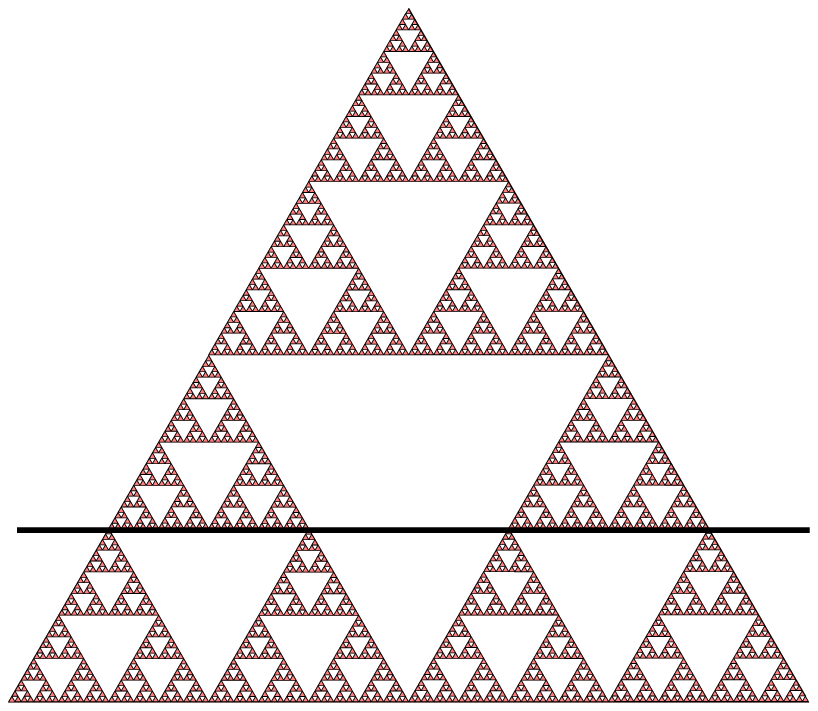}
\begin{center}
\textbf{Fig. 9.1.} \small{$\Omega_{3/4}$.}
\end{center}
\end{center}
\end{figure}

\begin{figure}[ht]
\begin{center}
\includegraphics[width=14.5cm,totalheight=4.5cm]{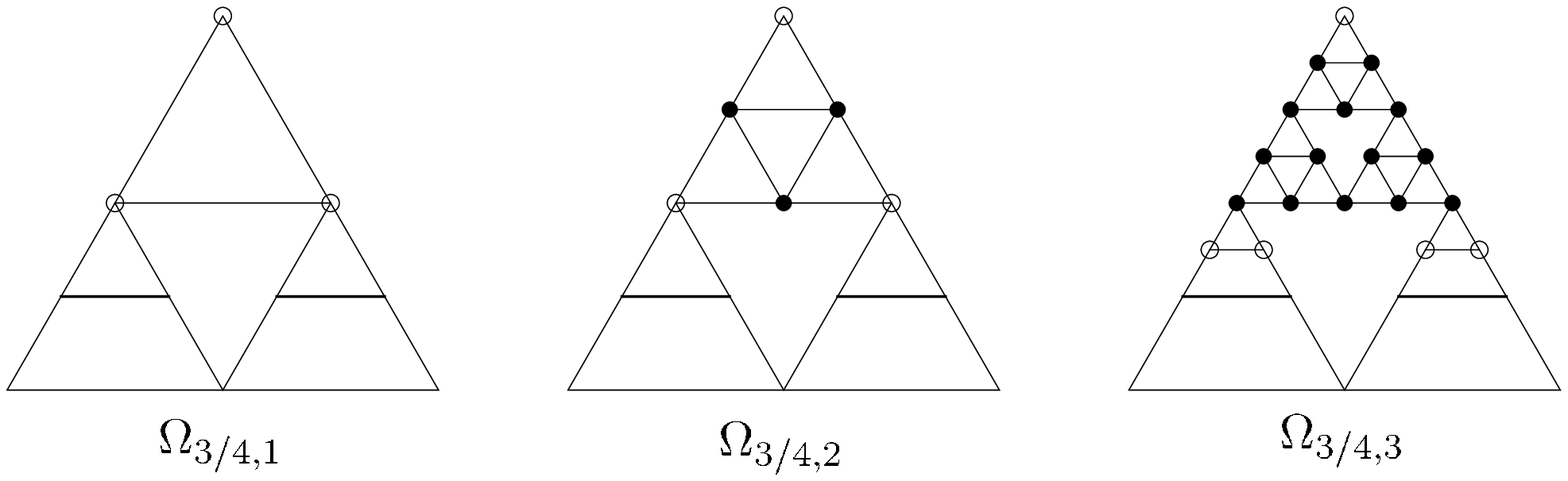}
\begin{center}
\textbf{Fig. 9.2.} \small{The first $3$ graphs, $\Omega_{3/4,1},
\Omega_{3/4,2}, \Omega_{3/4,3}$ in the approximations to
$\Omega_{3/4}$ with inside points and boundary points represented by
dots and circles respectively.}
\end{center}
\end{center}
\end{figure}

 We now focus on a particular example $\Omega_{3/4}$ to illustrate how
to extend the recipe for $\Omega_1$ to general case. We should be
particular interested in the primitive eigenvalues. We begin with
$\mathcal{P}_m^{+}(3/4)$, the symmetric case. It is convenient to
define the skeleton of $\Omega_m(3/4)$ by
$(q_0,F_1q_0,F_{10}F_1q_0,\cdots,F_{10}F_1^{m-2}q_0)$ for $m\geq 3$
and $(q_0,F_1q_0)$ for $m=1$ or $2$. Let $u_m$ be a
$\lambda_m$-eigenfunction of $-\Delta_m$ with
$\lambda_m\in\mathcal{P}_m^{+}(3/4)$. Denote by
$(b_0,b_1,b_2,\cdots,b_m)$ the values of $u_m$ on the skeleton of
$\Omega_{3/4,m}$ where $b_0=b_m=0$ by the Dirichlet boundary
condition. It is easy to observe that when $i\geq 2$, the eigenvalue
equation at the vertex $F_{10}F_1^{i-1}q_0$ is exactly same as that
of $\Omega_1$ case with suitably reindexed. Hence the generation
mechanism of symmetric primitive  eigenvalues is quite similar to
the $\Omega_1$ case. Based on this observation, one can easily find
that $\sharp\mathcal{P}_m^+(3/4)=2^m-2$ for $m\geq 2$ by still using
the weak spectral decimation method. A similar argument yields that
$\sharp\mathcal{P}_m^-(3/4)=2^m-2^{m-2}-2$ for $m\geq 2$.

To verify the eigenspace dimensional counting formula, we only look
at the first $4$ levels of approximations since the continued
process is similar.

When $m=1$, the result is trivial since there is no inside point in
$\Omega_{3/4,1}$. Hence $\sharp\mathcal{S}_1(3/4)=0=\sharp
V_1^{\Omega_{3/4}}\setminus\partial\Omega_{3/4,1}$.

When $m=2$, it is easy to check that there are only primitive
eigenvalues. Hence
$\sharp\mathcal{S}_2(3/4)=\sharp\mathcal{P}_2^+(3/4)+\sharp\mathcal{P}_2^-(3/4)=2+1=\sharp
V_2^{\Omega_{3/4}}\setminus\partial\Omega_{3/4,2}$.

When $m=3$, it is easy to check that there are $4$ initial localized
eigenvalues, i.e., $5$ with multiplicity $1$  and $6$ with
multiplicity $3$; there are $6$ symmetric primitive  eigenvalues and
$4$ skew-symmetric primitive eigenvalues; there is no miniaturized
eigenvalues. Hence
$\sharp\mathcal{S}_3(3/4)=\sharp\mathcal{L}_3(3/4)+\sharp\mathcal{P}_3^+(3/4)+\sharp\mathcal{P}_3^-(3/4)=4+6+4=\sharp
V_3^{\Omega_{3/4}}\setminus\partial\Omega_{3/4,3}$.

When $m=4$, it is easy to check that besides $1\cdot 2+3\cdot 1=5$
continued localized eigenvalues, there are $18$ initial localized
eigenvalues, i.e., $5$ with multiplicity $4$ and $6$ with
multiplicity $14$. Hence $\sharp\mathcal{L}_4(3/4)=5+18=23.$ There
are $14$ symmetric primitive  eigenvalues and $10$ skew-symmetric
primitive eigenvalues. Hence $\sharp\mathcal{P}_4(3/4)=14+10=24$.
Moreover, there are some miniaturized eigenvalues which come from
the miniaturizations of eigenvalues in $\mathcal{P}_2^-(1)$. Hence
$\sharp\mathcal{M}_4(3/4)=2\cdot\mathcal{P}_2^-(1)=2\cdot 2=4$. Thus
$\sharp\mathcal{S}_4(3/4)=23+24+4=\sharp
V_4^{\Omega_{3/4}}\setminus\partial\Omega_{3/4,4}$.

It is easy to verify the general formula for general $m$. We will
not attempt to list the details here. However, a more important fact
should be pointed out is that for $\Omega_{3/4}$ case, the
miniaturized eigenvalues in $\mathcal{M}_m(3/4)$ are generated not
from
 those in $\mathcal{P}_k^-(3/4)$ but from those in
$\mathcal{P}_k^-(1)$ for $k\leq m-2$. This means to study
$\mathcal{S}_m(3/4)$, one should first make clear
$\mathcal{S}_m(1)$. Things will be more complicated for general
$\Omega_x$.

Next we briefly present another observation. Still consider a domain
$\Omega_x$ with a series of graph approximations $\{\Omega_{x,m}\}$.
Notice that there are only two possible patterns when passing from
the $m$-level graph approximation to its next level. One is that the
boundary $\partial\Omega_{x,m+1}$ remains unchanged, i.e.,
$\partial\Omega_{x,m+1}=\partial\Omega_{x,m}$, the other  is that
$\partial\Omega_{x,m}\setminus\{q_0\}$ becomes a collection of
inside points of $\Omega_{x,m+1}$, i.e., each point in
$\partial\Omega_{x,m}\setminus\{q_0\}$ is connected with two new
$(m+1)$-level points in $\partial\Omega_{x,m+1}$. In fact, for the
$\mathcal{SG}\setminus{V_0}$ case, when passing from one level to
the next level, the boundaries of graphs are always $V_0$, keeping
unchanged. This is also the reason why spectral decimation can work
for $2$-series eigenvalues (which should be considered as the
primitive eigenvalues in $\mathcal{SG}\setminus{V_0}$ case). As for
the $\Omega_1$ case, when passing from one level to the next level,
the boundaries always change. Due to this phenomenon, the spectral
decimation recipe should be replaced by the weak spectral decimation
recipe for primitive or miniaturized eigenvalues since their
supports always touch the boundaries. For general $\Omega_x$
($0<x<1$), these two possible patterns can both exist. It is natural
to expect that under the first pattern, the two levels of primitive
eigenvalues are related by the spectral decimation (it is obviously
true), while under the second pattern, they are related by a weak
spectral decimation instead. Thus we expect:

\textbf{Conjecture 9.1.} \emph{For a domain $\Omega_x$ ($0<x<1$)
with a series of graph approximations $\{\Omega_{x,m}\}$, if the
boundaries change when passing from $m$-level to $(m+1)$-level, then
there is a weak spectral decimation relating the two levels of
 symmetric (or skew-symmetric) primitive eigenvalues. }

\section{Appendix}
\textbf{\quad Theorem A.} \emph{For each $m\geq 2$, let
$\lambda_{m,1},\lambda_{m,2},\cdots,\lambda_{m,r_m}$ be the $r_m$
distinct eigenvalues in $\mathcal{P}_m^+$ in increasing order. Then
$\lambda_{m+1,r_m+1}>2.$}

To prove this theorem, we need the following lemma:

\textbf{Lemma A.} $p_2(2)<0$, $p_3(2)>0$ and $(-1)^mp_m(2)>0$,
$\forall m\geq 4$.

\emph{Proof.} It is easy to check that $p_2(2)=-8<0$ and
$p_3(2)=68>0$.

Let $m\geq 4$. Then
\begin{eqnarray*}
p_m(x)&=&\frac{q_m(x)}{(x-2)(f(x)-2)\cdots
(f^{(m-3)}(x)-2)}\\
&=&\frac{s(f^{(m-2)}(x))q_{m-1}(x)-l(f^{(m-3)}(x))r(f^{(m-2)}(x))q_{m-2}(x)}{(x-2)(f(x)-2)\cdots
(f^{(m-3)}(x)-2)}\\
&=&\frac{s(f^{(m-2)}(x))}{f^{(m-3)}(x)-2}p_{m-1}(x)+\frac{-l(f^{(m-3)}(x))r(f^{(m-2)}(x))}{(f^{(m-4)}(x)-2)(f^{(m-3)}(x)-2)}p_{m-2}(x).
\end{eqnarray*}

Noticing that
$l(f^{(m-3)}(x))=f^{(m-3)}(x)-6=(f^{(m-4)}(x)-2)(3-f^{(m-4)}(x))$
and choosing $x=2$, we have
\begin{equation}\label{92}
p_m(2)=\frac{s(f^{(m-2)}(2))p_{m-1}(2)+2(2-f^{(m-2)}(2))(5-f^{(m-2)}(2))(3-f^{(m-4)}(2))p_{m-2}(2)}{f^{(m-3)}(2)-2}.
\end{equation}

We will prove the following stronger result than that stated in
Lemma A.
\begin{equation}\label{93}p_m(2)\sim(-1)^m \mbox{
and } p_{m+1}(2)+p_m(2)\sim (-1)^{m+1}, \forall m\geq
4.\end{equation}

Using $(\ref{92})$, it is easy to check that $p_4(2)=14064>0$ and
$p_5(2)=-593514756<0$ by a direct computation. Hence $(\ref{93})$
holds for $m=4$. In order to use the induction, we assume
$(\ref{93})$ holds for $m$ and will prove it for $m+1$.

First, it is easy to get that $p_{m+1}(2)\sim(-1)^{m+1}$, since
otherwise $p_{m+1}(2)+p_m(2)\sim(-1)^m$, which contradicts to the
induction assumption. Hence we only need to prove
$p_{m+2}(2)+p_{m+1}(2)\sim(-1)^m$.

Note that from $(\ref{92})$,
\begin{eqnarray*}
&&p_{m+2}(2)+p_{m+1}(2)\\=&&\frac{(s(f^{(m)}(2))+f^{(m-1)}(2)-2)p_{m+1}(2)+2(2-f^{(m)}(2))(5-f^{(m)}(2))(3-f^{(m-2)}(2))p_{m}(2)}{f^{(m-1)}(2)-2}\\
=&&a_mp_{m+1}(2)+b_m(p_{m+1}(2)+p_m(2)),
\end{eqnarray*}
where
$$a_m=\frac{s(f^{(m)}(2))+f^{(m-1)}(2)-2-2(2-f^{(m)}(2))(5-f^{(m)}(2))(3-f^{(m-2)}(2))}{f^{(m-1)}(2)-2}$$
and
$$b_m=\frac{2(2-f^{(m)}(2))(5-f^{(m)}(2))(3-f^{(m-2)}(2))}{f^{(m-1)}(2)-2}.$$

It is easy to check that $b_m<0$, since
$f^{(m)}(2)<f^{(m-1)}(2)<f^{(m-2)}(2)<0$ noticing that
$f^{(2)}(2)=-6$ and $m\geq 4$. We will prove that $a_m<0$ also. In
fact, the numerator of $a_m$ is
$s(\gamma)+f(\beta)-2-2(2-\gamma)(5-\gamma)(3-\beta)$, where
$\gamma:=f^{(m)}(2)$ and $\beta:=f^{(m-2)}(2)$ for simplicity. By
using $\gamma<f(\beta)<\beta\leq -6$, it is easy to get
\begin{eqnarray*}
s(\gamma)+f(\beta)-2&=&(2-\gamma)(4-\gamma)(5-\gamma)-14+3\gamma+f(\beta)-2\\
&>& (2-\gamma)(4-\gamma)(5-\gamma)-16+4\gamma\\
&>& (2-\gamma)(4-\gamma)(5-\gamma)-(2-\gamma)(5-\gamma)\\
&=& (3-\gamma)(2-\gamma)(5-\gamma),
\end{eqnarray*}
and
$$3-\gamma>3-f(\beta)=3-\beta(5-\beta)>3-5\beta>2(3-\beta).$$
Hence we have $s(\gamma)+f(\beta)-2>2(2-\gamma)(5-\gamma)(3-\beta)$.
Thus the numerator of $a_m$ is positive. Since the denominator of
$a_m$ is obviously negative, we get $a_m<0$.

Hence since $p_{m+1}(2)\sim(-1)^{m+1}$ we have proved before, and
$p_{m+1}(2)+p_m(2)\sim(-1)^{m+1}$ by the induction assumption, we
finally get $p_{m+2}(2)+p_{m+1}(2)\sim(-1)^m$. $\Box$

\emph{Proof of Theorem A.} Recall that in Lemma 4.4, we have proved
that $p_{m+1}(\phi_-(\lambda_{m,r_m}))\sim(-1)^{m+r_m-1}$ and
$p_{m+1}(\phi_+(\lambda_{m,r_m}))\sim(-1)^{m+r_m}$. Furthermore,
$\lambda_{m+1,r_m+1}$ is the only root of $p_{m+1}(x)$ between
$\phi_-(\lambda_{m,r_m})$ and $\phi_+(\lambda_{m,r_m})$.

When $m=2$, we have $p_3(\phi_-(\lambda_{2,r_2}))>0$ and
$p_3(\phi_+(\lambda_{2,r_2}))<0$ since $r_2$ is odd. By Lemma A, we
have $p_3(2)>0$. Since $\lambda_{3,r_2+1}$ is the only root between
$\phi_-(\lambda_{2,r_2})$ and $\phi_+(\lambda_{2,r_2})$, we get
$\lambda_{3,r_2+1}>2$.

When $m\geq 3$, we have
$p_{m+1}(\phi_-(\lambda_{m,r_m}))\sim(-1)^{m-1}$ and
$p_{m+1}(\phi_+(\lambda_{m,r_m}))\sim(-1)^m$ since $r_m$ is always
even. Still by Lemma A, we have $p_{m+1}(2)\sim(-1)^{m-1}$. Since
$\lambda_{m+1,r_m+1}$ is the only root between
$\phi_-(\lambda_{m,r_m})$ and $\phi_+(\lambda_{m,r_m})$, we get
$\lambda_{m+1,r_m+1}>2$. $\Box$

\textbf{Remark.} \emph{This theorem says that when $ m\geq 3$, the
first $m$-level initial eigenvalue is always greater than $2$. }

\textbf{Lemma B.} \emph{Let $m\geq 2$. Then $q_m(x)>0$ whenever
$0<x<\phi_-^{(m)}(6)$.}

\emph{Proof.}  Define $\theta_m(z)=q_m(\phi_-^{(m)}(z))$ on $0<z<6$,
$\forall m\geq 2$.

When $m=2$, $\theta_2(z)=q_2(\phi_-^{(2)}(z))$. Noticing that
$q_2(x)=s(x)$ and $q_2'(x)=-3x^2+22x-35$, an easy calculus shows
that $q_2(x)$ is monotone decreasing when $0<x<\phi_-^{(2)}(6)$.
Hence $\forall 0<x<\phi_-^{(2)}(6)$, we have
$q_2(0)=26>q_2(x)>q_2(\phi_-^{(2)}(6))\approx 12.68$. Thus
\begin{equation}\label{a1}
26>\theta_2(z)>12.68, \forall 0<z<6.
\end{equation}

When $m=3$,
$\theta_3(z)=q_3(\phi_-^{(3)}(z))=s(\phi_-^{(3)}(z))\theta_2(z)-l(\phi_-^{(3)}(z))r(\phi_-^{(2)}(z))$
on $0<z<6$. Noticing that $s(\phi_-^{(3)}(z))=q_2(\phi_-^{(3)}(z))$
and $q_2(x)$ is monotone decreasing when $0<x<\phi_-^{(2)}(6)$, we
have
$$s(0)=26>s(\phi_-^{(3)}(z))>s(\phi_-^{(3)}(6))\approx 22.96.$$
The monotone property of $-l(\phi_-^{(3)}(z))r(\phi_-^{(2)}(z))$ on
$0<z<6$ implies that
$$-84.21>-l(\phi_-^{(3)}(z))r(\phi_-^{(2)}(z))>-120.$$
Hence by using $(\ref{a1})$, we get
$$26\cdot 26-84.21=591.80>\theta_3(z)>22.96\cdot 12.68-120=171.16, \forall 0<z<6.$$

Hence $\theta_3(z)\geq 6 \theta_2(z)>0$ on $0<z<6$.

We now use induction to prove:
\begin{equation}\label{a2}\theta_{m+1}(z)\geq 6\theta_m(z)>0 \mbox{ on
} 0<z<6, \quad \forall m\geq 2.
\end{equation}

Of course, it holds for $m=2$. To use the induction, Assuming
$\theta_{m+1}(z)\geq 6\theta_m(z)>0$ on $0<z<6$, we will prove
$\theta_{m+2}(z)\geq 6\theta_{m+1}(z)>0$ on $0<z<6$.

Consider a polynomial
$g(x)=s(x)-\frac{1}{6}l(x)r(f(x))=6+\frac{115}{3}x-\frac{194}{3}x^2+\frac{89}{3}x^3-\frac{16}{3}x^4+\frac{1}{3}x^5$.
It is easy to compute that
$$g'(x)=\frac{115}{3}-\frac{388}{3}x+89x^2-\frac{64}{3}x^3+\frac{5}{3}x^4\geq \frac{115}{3}-\frac{388}{3}\phi_-^{(4)}(6)-\frac{64}{3}(\phi_-^{(4)}(6))^3\approx36.02>0$$
on $0<x<\phi_-^{(4)}(6).$ Hence $g(x)$ is a monotone increasing
function on the interval $[0,\phi_-^{(4)}(6)]$. So $g(x)\geq g(0)=6$
on $0<x<\phi_-^{(4)}(6)$.

By using an expansion along the last row of
$\theta_{m+2}(z)=q_{m+2}(\phi_-^{(m+2)}(z))$, we have
$$\theta_{m+2}(z)=s(\phi_-^{(m+2)}(z))\theta_{m+1}(z)-\frac{1}{6}l(\phi_-^{(m+2)}(z))r(\phi_-^{(m+1)}(z)) \cdot 6\theta_m(z).$$

By the induction assumption and the fact that
$\phi_-^{(m+2)}(z)<\phi_-^{(m+1)}(z)<2$, we have
\begin{eqnarray*}\theta_{m+2}(z)&\geq&
s(\phi_-^{(m+2)}(z))\theta_{m+1}(z)-\frac{1}{6}l(\phi_-^{(m+2)}(z))r(\phi_-^{(m+1)}(z))
 \theta_{m+1}(z) \\
&=&g(\phi_-^{m+2}(z))\theta_{m+1}(z).
\end{eqnarray*}

Since $0<\phi_-^{(m+2)}(z)<\phi_-^{(4)}(6)$ on $0<z<6$ when $m\geq
2$, we have $g(\phi_-^{(m+2)}(z))\geq 6$. Hence
$$\theta_{m+2}(z)\geq 6\theta_{m+1}(z)>0.$$

Hence we have proved $(\ref{a2})$ holds for $m+1$. From
$(\ref{a2})$, we get the desired result. $\Box$

\textbf{Acknowledgements.} This problem was originally considered by
Professor Robert S. Strichartz. I am grateful to him for addressing
me this problem and many illuminating discussions leading up to the
writing of this work. This work was done while I was visiting the
Department of Mathematics, Cornell University. Portions of this work
were presented at the 2012 Cornell Analysis Seminar. I express my
sincere gratitude to the department for its hospitality. I would
also like to thank the anonymous referee for several important
suggestions, especially for the refinement of Theorem 3.10, which
led to the improvement of the manuscript.

(Hua Qiu) DEPARTMENT OF MATHEMATICS, NANJING UNIVERISITY, NANJING,
210093, CHINA

\emph{E-mail address}: huatony@gmail.com

\end{document}